\documentclass[11pt]{article}
\usepackage[tbtags]{amsmath}
\usepackage{amssymb}
\usepackage{amsthm}
\usepackage[misc]{ifsym}
\usepackage{cases}
\usepackage{mathrsfs}
\usepackage{color}
\usepackage{graphicx}
\usepackage{subfigure}
\numberwithin{equation}{section}
\normalsize
\title{\bf Linear Quadratic Leader-follower Stochastic Differential Games: Closed-Loop Solvability
\thanks{This work is supported by National Key R\&D Program of China (Grant No. 2018YFB1305400), National Natural Science Foundations of China (Grant Nos. 11971266, 11831010, 11571205), and Shandong Provincial Natural Science Foundations (Grant Nos. ZR2020ZD24, ZR2019ZD42).}}
\author{\normalsize Zixuan Li\thanks{\it School of Mathematics, Shandong University, Jinan 250100, P.R. China, E-mail: 201812064@mail.sdu.edu.cn},\quad Jingtao Shi\thanks{\it Corresponding author. School of Mathematics, Shandong University, Jinan 250100, P.R. China, E-mail: shijingtao@sdu.edu.cn}}
\newtheorem{mypro}{Proposition}[section]
\newtheorem{mythm}{Theorem}[section]
\newtheorem{mydef}{Definition}[section]

\newtheorem{Remark}{Remark}[section]
\newtheorem{Example}{Example}[section]

\begin{document}
\maketitle

\noindent{\bf Abstract:}\quad In this paper, a leader-follower stochastic differential game is studied for a linear stochastic differential equation with a quadratic cost functional. The coefficients in the state equation and the weighting matrices in the cost functionals are all deterministic. Closed-loop strategies are introduced, which require to be independent of initial states; and such a nature makes it very useful and convenient in applications. The follower first solves a stochastic linear quadratic optimal control problem, and his optimal closed-loop strategy is characterized by a Riccati equation, together with an adapted solution to a linear backward stochastic differential equation. Then the leader turns to solve a stochastic linear quadratic optimal control problem of a forward-backward stochastic differential equation, necessary conditions for the existence of optimal closed-loop strategies for the leader is given by the existence of a Riccati equation. Some examples are also given.
\vspace{2mm}


\vspace{2mm}


\section{Introduction}

Let us first introduce some notations which will be used throughout the paper.

Let $T>0$ be a finite time duration. Let $\mathbb{R}^{n \times m}$ be the set of all $(n \times m)$ matrices, and let $\mathbb{S}^n$ be the set of all $(n \times n)$ matrices. For any Banach space $H$ (for example, $H=\mathbb{R}^n,\mathbb{R}^{n \times m},\mathbb{S}^n$), let $L^p(0,T;H)\,(1 \leqslant p \leqslant \infty)$ be the space of all $H$-valued functions that are $L^p$-integrable on $[0,T]$, and let $C([0,T];H)$ be the space of all $H$-valued continuous functions on $[0,T]$.

Let $(\Omega,\mathcal{F},\mathbb{F},\mathbb{P})$ be a completed filtered probability space on which a standard one-dimensional Brownian motion $W=\{W(t);0 \leqslant t < \infty \}$ is defined, where $\mathbb{F}=\{\mathcal{F}_t\}_{t\geqslant 0}$ is the natural filtration of $W$ augmented by all the $\mathbb{P}$-null sets in $\mathcal{F}$. We denote
\begin{equation*}
\begin{aligned}
&L^2_{\mathcal{F}_t}(\Omega;H)=\Big\{\xi:\Omega \to H\,|\,\xi\,\,\mbox{is}\,\,\mathcal{F}_t\mbox{-measurable},\,\,\,\mathbb{E}|\xi|^2<\infty \Big\},\,\, t\in(0,T],\\
&L^2_{\mathbb{F}}(0,T;H)=\Big\{f(\cdot):[0,T]\times\Omega\to H\,\big|\,f(\cdot)\,\,\mbox{is}\,\,\mathbb{F}\mbox{-progressively measurable},\,\,\mathbb{E}\int_0^T|f(s)|^2ds<\infty \Big\}.
\end{aligned}
\end{equation*}

We consider the following controlled linear {\it stochastic differential equation} (SDE for short):
\begin{equation}\label{state}
\left\{
\begin{aligned}
dx(t)&=\big[A(s)x(s)+B_1(s)u_1(s)+B_2(s)u_2(s)+b(s)\big]ds\\
     &\quad+\big[C(s)x(s)+D_1(s)u_1(s)+D_2(s)u_2(s)+\sigma(s)\big]dW(s),\quad s \in [0,T],\\
 x(0)&=x,
\end{aligned}
\right.
\end{equation}
where $x\in\mathbb{R}^n$, $A(\cdot),B_i(\cdot),C(\cdot),D_i(\cdot),\,\,i=1,2$ are given deterministic matrix-valued functions of proper dimensions; $b(\cdot),\,\sigma(\cdot)$ are vector-valued $\mathbb{F}$-progressively measurable processes. In the above, $x(\cdot)$ is the \textit{state process} with values in $\mathbb{R}^n$, and $u_1(\cdot),u_2(\cdot)$ are \textit{control processes} taken by the two players in the games, labeled 1 and 2, with values in $\mathbb{R}^{m_1}$ and $\mathbb{R}^{m_2}$, respectively. We introduce the following Hilbert spaces:
\begin{equation}\label{ac}
\begin{aligned}
\mathcal{U}_i[0,T]=\bigg\{u_i:\,&[0,T] \times \Omega \to \mathbb{R}^{m_i} \,\,\big|\,\,u_i(\cdot)\,\, \mbox{is}\, \,\mathbb{F}\mbox{-progressively measurable}, \\
                                &\mathbb{E}\int_{0}^{T} |u_i(s)|^2 ds < \infty \bigg\},\,\quad \,i=1,2.
\end{aligned}
\end{equation}
The control processes $u_1(\cdot)\in\mathcal{U}_1[0,T]$ and $u_2(\cdot)\in\mathcal{U}_2[0,T]$ are called \textit{admissible controls}.

Under some mild conditions on the coefficients, for any $(x,u_1(\cdot),u_2(\cdot)) \in \mathbb{R}^n \times \mathcal{U}_1[0,T] \times \mathcal{U}_2[0,T]$, there exists a unique (strong) solution $x(\cdot)\equiv x(\cdot;x,u_1(\cdot),u_2(\cdot)) \in L^2_{\mathcal{F}}(0,T;\mathbb{R}^n)$ to (\ref{state}). Thus, we can define the \textit{cost functionals} for the players as follows. For $i=1,2$,
\begin{equation}\label{cfi}
\begin{aligned}
&J_i(x;u_1(\cdot),u_2(\cdot))=\mathbb{E}
\Bigg\{\int_0^T \bigg[\bigg\langle
\left( \begin{array}{ccc} Q^i(s) & S^i_1(s)^\top  & S^i_2(s)^\top \\
S^i_1(s) & R^i_{11}(s) & R^i_{12}(s)\\
S^i_2(s) & R^i_{21}(s) & R^i_{22}(s)\\
\end{array} \right)
\left( \begin{array}{c} x(s) \\ u_1(s) \\ u_2(s) \end{array} \right),
\left( \begin{array}{c} x(s) \\ u_1(s) \\ u_2(s) \end{array}\right)\bigg\rangle \\
&\qquad +2\bigg\langle
\left( \begin{array}{c} q^i(s) \\ \rho^i_1(s) \\ \rho^i_2(s) \end{array} \right),
\left( \begin{array}{c} x(s) \\ u_1(s) \\ u_2(s) \end{array}\right)\bigg\rangle \bigg]ds+\big\langle G^ix(T),x(T)\big\rangle+2\big\langle g^i ,x(T)\big\rangle \Bigg\},
\end{aligned}
\end{equation}
where, for $i=1,2$, $G^i$ are symmetric matrices and  $Q^i(\cdot),S^i_1(\cdot),S^i_2(\cdot),R^i_{11}(\cdot),R^i_{12}(\cdot),R^i_{21}(\cdot)$ and $R^i_{22}(\cdot)$ are deterministic matrix-valued functions of proper dimensions with
\begin{equation}\nonumber
Q^i(\cdot)^\top =Q^i(\cdot),\,\quad\,R^i_{jj}(\cdot)^\top =R^i_{jj}(\cdot),\,\quad\,R^i_{12}(\cdot)^\top =R^i_{21}(\cdot),\qquad\,\, i,j=1,2,
\end{equation}
$g^i$ are $\mathcal{F}_T$-measurable random vectors; $q^i(\cdot),\,\,\rho^i_1(\cdot),\,\,\rho^i_2(\cdot)$ are vector-valued $\mathbb{F}$-progressively measurable processes.

In the Stackelberg game (also known as leader-follower game) framework, Player 1 is the follower and Player 2 is the leader. For any choice $u_2(\cdot) \in \mathcal{U}_2[0,T]$ of Player 2 and a fixed initial state $x \in \mathbb{R}^n$, Player 1 would like to choose a $\bar{u}_1(\cdot) \in \mathcal{U}_1[0,T]$ such that $J_1(x;\bar{u}_1(\cdot),u_2(\cdot))$ is the minimum of $J_1(x;u_1(\cdot),u_2(\cdot))$ over $u_1(\cdot) \in \mathcal{U}_1[0,T]$. Knowing Player 1 would take such an optimal control $\bar{u}_1(\cdot)$, Player 2 would like to choose some $\bar{u}_2(\cdot) \in \mathcal{U}_2[0,T]$ to minimize $J_2(x;\bar{u}_1(\cdot),u_2(\cdot))$ over $u_2(\cdot) \in \mathcal{U}_2[0,T]$. We refer to such a problem as a {\it linear quadratic (LQ for short) leader-follower (Stackelberg) stochastic differential game}.

In a more rigorous way, Player 1 wants to find a mapping $\bar{\mu}_1:\mathcal{U}_2[0,T] \times \mathbb{R}^n \to \mathcal{U}_1[0,T] $ and Player 2 want to find a $\bar{u}_2(\cdot) \in \mathcal{U}_2[0,T]$ such that
\begin{equation}\label{problem}
\begin{cases}
J_1(x;\bar{\mu}_1[u_2(\cdot),x](\cdot),u_2(\cdot))=\mathop{\min}\limits_{u_1(\cdot) \in \mathcal{U}_1[0,T]}J_1(x;u_1(\cdot),u_2(\cdot)),\quad \,\,\forall u_2(\cdot) \in \mathcal{U}_2[0,T],\\
J_2(x;\bar{\mu}_1[\bar{u}_2(\cdot),x](\cdot),\bar{u}_2(\cdot))=\mathop{\min}\limits_{u_2(\cdot) \in \mathcal{U}_2[0,T]}J_2(x;\bar{\mu}_1[u_2(\cdot),x](\cdot),u_2(\cdot)).
\end{cases}
\end{equation}
If the above pair $(\bar{\mu}_1[\cdot,x],\bar{u}_2(\cdot))$ exists, we refer to it as an \textit{open-loop solution} or \textit{open-loop Stackelberg equilibrium} to the above LQ leader-follower stochastic differential game, for $x\in\mathbb{R}^n$.

\par The theory of leader-follower game can be traced back to Stackelberg \cite{Stackelberg1934}, who put forward the concept of Stackelberg equilibrium in static competitive economics with a hierarchical structure. Simann and Cruz \cite{SC1973-1,SC1973-2} studied the multi-stages and dynamic LQ leader-follower differential games, where feedback Stackelberg equilibria are introduced. Castanon and Athans \cite{CA1976} considered feedback Stackelberg strategies for the two person linear multi-stages game with quadratic performance criteria and noisy measurements and gave an explicit solution when the information sets are nested in a stochastic case. Bagchi and Ba\c{s}ar \cite{BB1981} investigated the LQ leader-follower stochastic differential game, where the diffusion term of the state equation does not contain the state and control variables. Yong \cite{Yong2002} extended the LQ leader-follower stochastic differential game to random and state-control dependent coefficients, and obtained the feedback representation of the open-loop equilibrium via some stochastic Riccati equations. In the past decades, there have been a great deal of works on this issue, for jump diffusions see \O ksendal et al. \cite{OSU2013}, for different information structures see Ba\c{s}ar and Olsder \cite{BO1998}, Bensoussan et al. \cite{BCS2015}, for time-delayed systems see Xu and Zhang \cite{XZ2016}, Xu et al. \cite{XSZ2018}, for mean field's type models related with multiple followers and large populations see Mukaidani and Xu \cite{MX2015}, Moon and Ba\c{s}ar \cite{MB2018}, Li and Yu \cite{LiYu2018}, Lin et al. \cite{LJZ2019}, Wang and Zhang \cite{WZ2020}, Huang et al. \cite{HSW2020}, for partial/asymmetric/overlapping information see Shi et al. \cite{SWX2016,SWX2017,SWX2020}, for backward stochastic systems see Du and Wu \cite{DW2019}, Zheng and Shi \cite{ZS2020}.

Our interest in this paper lies in the {\it closed-loop solution} or the {\it closed-loop Stackelberg equilibrium} to the above LQ leader-follower stochastic differential game. To our best knowledge, this topic has not been studied in the literature yet, except for our recent conference paper \cite{LiShi2021} for a special case where the stochastic system is homogeneous and the diffusion is independent of the control. In 2014, Sun and Yong \cite{SY2014} introduce the notions of open-loop and closed-loop solvabilities for the LQ stochastic optimal control problem, which is a special case when only one player/controller is considered for open-loop and closed-loop saddle points for an LQ two-person zero-sum stochastic differential game. Sun et al. \cite{SLY2016} further gives more detailed characterizations of the closed-loop solvability for the LQ stochastic optimal control problem. Sun and Yong \cite{SY2019} is devoted to the open-loop and closed-loop Nash equilibria for the LQ two-person nonzero-sum stochastic differential game. The existence of an optimal closed-loop strategy for an LQ mean-field optimal control problem is studied in Li et al. \cite{LSY2016}. Sun and Yong \cite{SY2018} obtained the equivalence of open-loop and closed-loop solvabilities for the LQ stochastic optimal control problem in an infinite horizon. See also their paragraph \cite{SunYong2020}. Very recently, Li et al. \cite{LSY2020} extended the previous results to LQ mean-field two-person zero-sum and nonzero sum stochastic differential games in an infinite horizon.

In this paper, comparing with our former paper \cite{LiShi2021}, we consider an LQ leader-follower stochastic differential game in a more general framework, with control-state dependent diffusion term in the state equation. We first solve the follower's stochastic optimal control problem, and his optimal closed-loop strategy is characterized by a Riccati equation, together with an adapted solution to a linear \textit{backward stochastic differential equation} (BSDE for short). Then we solve the leader's problem, which is a stochastic optimal control problem of a \textit{forward-backward stochastic differential equation} (FBSDE for short). We will give the definition of the the optimal closed-loop strategy for the leader, and necessary condition for the existence of it is given by some new Riccati equations.

The rest of this paper is organized as follows. Section 2 gives some preliminaries, to introduce closed-loop Stackelberg equilibria for the LQ Stackelberg stochastic differential game. Section 3 is devoted to solve the optimization problem of the follower. With the aid of a Riccati equation, the sufficient and necessary conditions of the closed-loop solvability for the follower's problem are given. In Section 4, necessary conditions for the closed-loop solvability for the leader's problem is obtained. In Section 5, the relationship between open-loop solvability and closed-loop solvability is illustrated by some examples. Finally, in Section 6 some concluding remarks are given.

%Problem(SLQ)
\section{Preliminaries}

First of all, we recall the open-loop and closed-loop solvabilities for the LQ stochastic optimal control problem (see \cite{SLY2016}). Consider the linear state equation
\begin{equation}\label{SLQ}
\left\{\begin{aligned}
dX(s)&=\big[A(s)X(s)+B(s)u(s)+b(s)\big]ds\\&\qquad +\big[C(s)X(s)+D(s)u(s)+\sigma(s)\big]dW(s),\ \ s \in [0,T],\\
 X(0)&=x,
\end{aligned}\right.
\end{equation}
and the quadratic cost functional:
\begin{equation}\label{cf}
\begin{split}
J(t,&x;u(\cdot))=\mathbb{E} \bigg\{
\int_t^T \bigg[\bigg\langle
\left( \begin{array}{cc} Q(s) & S(s)^\top \\ S(s) & R(s)\end{array} \right)
\left( \begin{array}{c} X(s) \\ u(s)\end{array} \right),
\left( \begin{array}{c} X(s) \\ u(s)\end{array}\right)\bigg\rangle \\
&+2\bigg\langle
\left( \begin{array}{c} q(s) \\ \rho(s) \end{array} \right),
\left( \begin{array}{c} X(s) \\ u(s)\end{array}\right)\bigg\rangle \bigg]ds
+\big\langle GX(T),X(T)\big\rangle+2\big\langle g,X(T)\big\rangle \bigg\}.
\end{split}
\end{equation}
We adopt the following assumptions.
\par \textbf{(S1)} The coefficients of the state equation (\ref{SLQ}) satisfy the following:
\begin{equation}\nonumber
\begin{cases}
A(\cdot) \in L^1(0,T;\mathbb{R}^{n \times n}),\,\,B(\cdot) \in L^2(0,T;\mathbb{R}^{n \times m}),\,\,b(\cdot) \in L^2_\mathbb{F}(\Omega;L^1(0,T;\mathbb{R}^n)),\\
C(\cdot) \in L^2(0,T;\mathbb{R}^{n \times n}),\,\,D(\cdot) \in L^\infty(0,T;\mathbb{R}^{n \times m}),\,\,\sigma(\cdot) \in L^2_\mathbb{F}(0,T;\mathbb{R}^n).
\end{cases}
\end{equation}
\par \textbf{(S2)} The weighting coefficients of the cost functional (\ref{cf}) satisfy the following:
\begin{equation}\nonumber
\begin{cases}
Q(\cdot) \in L^1(0,T;\mathbb{S}^n),\,\,S(\cdot) \in L^2(0,T;\mathbb{R}^{n \times m}),\,\,R(\cdot) \in L^{\infty}(0,T;\mathbb{S}^{n})\\
q(\cdot) \in L^2_{\mathbb{F}}(\Omega;L^1(0,T;\mathbb{R}^{n})),\,\,\rho(\cdot) \in L^2(0,T;\mathbb{R}^m),\,\,g \in L^2_{\mathcal{F}_T}(\Omega;\mathbb{R}^{n}),\,\,G \in \mathbb{S}^n.\\
\end{cases}
\end{equation}

Under (S1) and (S2), for any $x \in \mathbb{R}^n$ and $u(\cdot) \in \mathcal{U}[0,T] \equiv L^2_{\mathbb{F}}(0,T;\mathbb{R}^m)$, the state equation (\ref{SLQ}) admits a unique strong solution and the cost functional is well-defined. Therefore, the following problem is meaningful.

\textbf{Problem (SLQ)}. For any initial state $x \in \mathbb{R}^n$, find a $\bar{u}(\cdot) \in \mathcal{U}[0,T]$ such that
\begin{equation}\label{PLQ}
J(x;\bar{u}(\cdot))=\mathop{\inf}\limits_{u(\cdot) \in \mathcal{U}[0,T]}J(x;u(\cdot)) \equiv V(x).
\end{equation}
Any $\bar{u}(\cdot) \in \mathcal{U}[0,T]$ satisfying (\ref{PLQ}) is called \textit{an open-loop optimal control} of Problem (SLQ) for $x$, the corresponding $\bar{X}(\cdot)\equiv X(\cdot;x,\bar{u}(\cdot))$ is called \textit{an open-loop optimal state process} and $(\bar{X}(\cdot),\bar{u}(\cdot))$ is called \textit{an open-loop optimal pair}. The map $V(\cdot)$ is called \textit{the value function} of Problem (SLQ).

\begin{mydef}\label{def2.1}
Let $x\in \mathbb{R}^n$. If there exists a (unique) $\bar{u}(\cdot) \in \mathcal{U}[0,T]$ such that (\ref{PLQ}) holds, then we say that Problem (SLQ) is (uniquely) open-loop solvable at $x$. If Problem (SLQ) is (uniquely) open-loop solvable for every $x \in \mathbb{R}^n$, then we say that Problem (SLQ) is (uniquely) open-loop solvable on $\mathbb{R}^n$.
\end{mydef}

The following result is concerned with open-loop solvability of Problem (SLQ) for a given initial state, whose proof can be found in \cite{SLY2016} (see also \cite{SY2014}).

% SLQ open-loop optimal solvability
\begin{mypro}\label{opsn}
Let (S1)-(S2) hold. For an initial state $x \in \mathbb{R}^n$, a state-control pair $(\bar{X}(\cdot),\bar{u}(\cdot))$ is an open-loop optimal pair of Problem (SLQ) if and only if the following hold:\\
(\romannumeral 1) The stationarity condition holds:
\begin{equation}\label{sc}
B(s)^\top\bar{Y}(s)+D(s)^\top\bar{Z}(s)+S(s)\bar{X}(s)+R(s)\bar{u}(s)+\rho(s)=0,\quad a.e.\ s \in[0,T],\ \mathbb{P}\mbox{-}a.s.,
\end{equation}
where $(\bar{Y}(\cdot),\bar{Z}(\cdot)) \in L^2_{\mathbb{F}}(0,T;\mathbb{R}^n) \times L^2_{\mathbb{F}}(0,T;\mathbb{R}^n)$ is the adapted solution to the following BSDE:
\begin{equation}\left\{
\begin{aligned}
d\bar{Y}(s)&=-\big[A(s)^\top\bar{Y}(s)+C(s)^\top\bar{Z}(s)+Q(s)\bar{X}(s)+S(s)^\top\bar{u}(s)+q(s)\big]ds\\
           &\qquad\quad  +\bar{Z}(s)dW(s),\,\,\quad\,s \in [0,T],\\
 \bar{Y}(T)&=G\bar{X}(T)+g.
\end{aligned}\right.
\end{equation}
(\romannumeral 2) The map $u(\cdot) \to J(0;u(\cdot))$ is convex.
\end{mypro}
Next, take $\Theta(\cdot) \in L^2(0,T;\mathbb{R}^{m \times n}) \equiv \mathcal{Q}[0,T]$ and $v(\cdot) \in \mathcal{U}[0,T]$. For any $x \in \mathbb{R}^n$, let us consider the following equation on $[0,T]$:
\begin{equation}\label{cl}
\left\{
\begin{aligned}
dX(s)&=\big\{[A(s)+B(s)\Theta(s)]X(s)+B(s)v(s)+b(s)\big\}ds\\
     &\qquad +\big\{[C(s)+D(s)\Theta(s)]X(s)+D(s)v(s)+\sigma(s)\big\}dW(s),\quad   s \in [0,T],\\
 X(0)&=x,
\end{aligned}\right.
\end{equation}
which admits a unique solution $X(\cdot) \equiv X(\cdot;x,\Theta(\cdot),v(\cdot))$, depending on the $\Theta(\cdot)$ and $v(\cdot)$. The above equation (\ref{cl}) is called a {\it closed-loop system} of the original state equation (\ref{SLQ}) under a {\it closed-loop strategy} $(\Theta(\cdot),v(\cdot))$. We point out that $(\Theta(\cdot),v(\cdot))$ is independent of the initial state $x \in \mathbb{R}^n$. With the above solution $X(\cdot)$, we define
\begin{equation}\label{ccf}
\begin{aligned}
J(t,&x;\Theta(\cdot)X(\cdot)+v(\cdot))=\mathbb{E} \bigg\{
\int_t^T \bigg[\bigg\langle
\left( \begin{array}{cc} Q & S^T \\ S & R\end{array} \right)
\left( \begin{array}{c}  X \\ \Theta X+v\end{array} \right),
\left( \begin{array}{c}  X \\ \Theta X+v\end{array}\right)\bigg\rangle\\
& +2\bigg\langle
\left( \begin{array}{c} q \\ \rho \end{array} \right),
\left( \begin{array}{c} X \\ \Theta X +v \end{array}\right)\bigg\rangle \bigg]ds+\big\langle GX(T),X(T)\big\rangle+2\big\langle g,X(T)\big\rangle \bigg\},
\end{aligned}
\end{equation}
and recall the following definition.
\begin{mydef}\label{def2.2}
A pair $(\bar{\Theta}(\cdot),\bar{v}(\cdot)) \in \mathcal{Q}[0,T] \times \mathcal{U}[0,T]$ is called a \textit{closed-loop optimal strategy} of Problem (SLQ) on $[0,T]$ if
\begin{equation}
\begin{split}
J(x;\bar{\Theta}(\cdot)\bar{X}(\cdot)+\bar{v}(\cdot)) \leqslant J(x;\Theta(\cdot)X(\cdot)+v(\cdot)),\\
\forall x \in \mathbb{R}^n,\,\,\forall  (\Theta(\cdot),v(\cdot)) \in \mathcal{Q}[0,T] \times \mathcal{U}[0,T],
\end{split}
\end{equation}
where $\bar{X}(\cdot)\equiv X(\cdot;x,\bar{\Theta}(\cdot),\bar{v}(\cdot))$, and $X(\cdot)\equiv X(\cdot;x,\Theta(\cdot),v(\cdot))$.
\end{mydef}
We emphasize that the pair $(\bar{\Theta}(\cdot),\bar{v}(\cdot))$ is required to be independent of the initial state $x \in \mathbb{R}^n$. The following result is also from \cite{SLY2016}.

\begin{mypro}\label{relation}
Let (S1)-(S2) hold and let $(\bar{\Theta}(\cdot),\bar{v}(\cdot)) \in \mathcal{Q}[0,T] \times \mathcal{U}[0,T]$. Then the following statements are equivalent:\\
(\romannumeral 1) $(\bar{\Theta}(\cdot),\bar{v}(\cdot))$ is a closed-loop optimal strategy of Problem (SLQ) on $[0,T]$.\\
(\romannumeral 2) For any $x \in \mathbb{R}^n$ and $v(\cdot) \in \mathcal{U}[0,T]$,
\begin{equation}\nonumber
J(x;\bar{\Theta}(\cdot)\bar{X}(\cdot)+\bar{v}(\cdot)) \leqslant J(x;\bar{\Theta}(\cdot)X(\cdot)+v(\cdot)),
\end{equation}
\par where $\bar{X}(\cdot)\equiv X(\cdot;x,\bar{\Theta}(\cdot),\bar{v}(\cdot))$ and $X(\cdot)\equiv X(\cdot;x,\bar{\Theta}(\cdot),v(\cdot))$.\\
(\romannumeral 3) For any $x \in \mathbb{R}^n$ and $u(\cdot) \in \mathcal{U}[0,T]$,
\begin{equation}\label{open and close relation}
J(x;\bar{\Theta}(\cdot)\bar{X}(\cdot)+\bar{v}(\cdot)) \leqslant J(x;u(\cdot)),
\end{equation}
\par where $\bar{X}(\cdot)\equiv X(\cdot;x,\bar{\Theta}(\cdot),\bar{v}(\cdot))$.	
\end{mypro}

From the above result, we see that if $(\bar{\Theta}(\cdot),\bar{v}(\cdot))$ is a closed-loop optimal strategy of Problem (SLQ) on $[0,T]$, then for any fixed initial state $x \in \mathbb{R}^n$, with $\bar{X}(\cdot)$ denoting the state process corresponding to $x$ and $(\bar{\Theta}(\cdot),\bar{v}(\cdot))$, (\ref{open and close relation}) implies that the {\it outcome}
\begin{equation}\nonumber
\bar{u}(\cdot)=\bar{\Theta}(\cdot)\bar{X}(\cdot)+\bar{v}(\cdot) \in \mathcal{U}[0,T]
\end{equation}
is an open-loop optimal control of Problem (SLQ) for $x$. Therefore, for Problem (SLQ), the existence of closed-loop optimal strategies on $[0,T]$ implies the existence of open-loop optimal controls for any $x \in \mathbb{R}^n$.

\vspace{1mm}

We now return to our LQ leader-follower stochastic differential game (\ref{state})-(\ref{problem}). We denote $L^2(0,T;\mathbb{R}^{m_i \times n}) \equiv \mathcal{Q}_i[0,T]$ for $i=1,2$.

First, for any $u_2(\cdot)\in \mathcal{U}_2[0,T]$, take $\Theta_1(\cdot) \in \mathcal{Q}_1[0,T]$ and $v_1(\cdot)\in \mathcal{U}_1[0,T]$. For any $x \in \mathbb{R}^n$, let us consider the following equation on $[0,T]$:
\begin{equation}\label{cse}\left\{
\begin{aligned}
dx^{u_2}(s)&=\Big\{\big[A(s)+B_1(s)\Theta_1(s)\big]x^{u_2}(s)+B_1(s)v_1(s)+B_2(s)u_2(s)+b(s)\Big\}ds\\
           &\quad
           +\Big\{[C(s)+D_1(s)\Theta_1(s)]x^{u_2}(s)+D_1(s)v_1(s)+D_2(s)v_2(s)+\sigma(s)\Big\}dW(s),\ s \in [0,T],\\
 x^{u_2}(0)&=x,
\end{aligned}\right.
\end{equation}
which admits a unique solution $x^{u_2}(\cdot) \equiv x(\cdot;x,\Theta_1(\cdot),v_1(\cdot),u_2(\cdot))$, depending on $\Theta_1(\cdot)$ and $v_1(\cdot)$. The above is called a {\it closed-loop system} of the original state equation (\ref{state}) under the {\it closed-loop strategy} $(\Theta_1(\cdot),v_1(\cdot))$ for the follower. We point out that $(\Theta_1(\cdot),v_1(\cdot))$ is independent of the initial state $x$. With the above solution $x^{u_2}(\cdot)$, we define (Some time variables are omitted if there is no ambiguity.)
\begin{equation}\label{fccf}
\begin{aligned}
&J_1(x;\Theta_1x^{u_2}+v_1,u_2)\\
&=\mathbb{E}\bigg\{\int_0^T \Big[\big\langle (Q^1+\Theta_1^\top S^1_1+S^{1\top}_1\Theta_1+\Theta_1^\top R^1_{11}\Theta_1)x^{u_2},x^{u_2} \big\rangle +2\big\langle (S^1_1+R^1_{11}\Theta_1)x^{u_2},v_1 \big\rangle\\
&\qquad\quad +\big\langle R^1_{11}v_1,v_1 \big\rangle +2\big\langle S^{1\top}_2u_2+q^1+\Theta_1^\top R^1_{12}u_2+\Theta_1^\top \rho^1_1,x^{u_2} \big\rangle +2\big\langle R^1_{12}u_2 +\rho^1_1,v_1 \big\rangle\\
&\qquad\quad +2\big\langle \rho^1_2,u_2 \big\rangle +\big\langle R^1_{22}u_2,u_2 \big\rangle \Big]ds +\big\langle G^1x^{u_2}(T),x^{u_2}(T)\big\rangle +2\big\langle g^1 ,x^{u_2}(T)\big\rangle \bigg\},
\end{aligned}
\end{equation}
and introduce the following notion.
\begin{mydef}\label{def2.3}
A 4-tuple $(\bar{\Theta}_1(\cdot),\bar{v}_1(\cdot),\bar{\Theta}_2(\cdot),\bar{v}_2(\cdot)) \in \mathcal{Q}_1[0,T] \times \mathcal{U}_1[0,T] \times \mathcal{Q}_2[0,T] \times \mathcal{U}_2[0,T] $ is called a (unique) closed-loop Stackelberg equilibrium of our LQ leader-follower stochastic differential game on $[0,T]$ if\\
(\romannumeral 1) For any $x \in \mathbb{R}^n$ and given $u_2(\cdot) \in \mathcal{U}_2[0,T]$, Player 1 could find two maps: $\bar{\Theta}_1:\mathcal{U}_2[0,T]\rightarrow\mathcal{Q}_1[0,T]$ and $\bar{v}_1:\mathcal{U}_2[0,T]\rightarrow\mathcal{U}_1[0,T]$ such that
\begin{equation}
\begin{aligned}
J_1(x;\bar{\Theta}_1[u_2](\cdot)\bar{x}^{u_2}(\cdot)+\bar{v}_1[u_2](\cdot),u_2(\cdot))\leqslant&\ J_1(x;\Theta_1[u_2](\cdot)x^{u_2}(\cdot)+v_1[u_2](\cdot),u_2(\cdot)),\\
&\qquad \forall\, (\Theta_1(\cdot),v_1(\cdot)) \in \mathcal{Q}_1[0,T] \times \mathcal{U}_1[0,T],
\end{aligned}
\end{equation}
where $\bar{x}^{u_2}(\cdot)\equiv x(\cdot;x,\bar{\Theta}_1[u_2](\cdot),\bar{v}_1[u_2](\cdot),u_2(\cdot))$ and $x^{u_2}(\cdot)\equiv x(\cdot;x,\Theta_1[u_2](\cdot),v_1[u_2](\cdot),u_2(\cdot))$.\\
(\romannumeral 2) There exist a (unique) pair $(\bar{\Theta}_2(\cdot), \bar{v}_2(\cdot))\in \mathcal{Q}_2[0,T] \times \mathcal{U}_2[0,T]$  such that
\begin{equation}
\begin{aligned}
&J_2\big(x;\bar{\Theta}_1[\bar{\Theta}_2\bar{x}+\bar{v}_2](\cdot)\bar{x}(\cdot)+\bar{v}_1[\bar{\Theta}_2\bar{x}+\bar{v}_2](\cdot),\bar{\Theta}_2(\cdot)\bar{x}(\cdot)+\bar{v}_2(\cdot)\big)\\
&\leqslant J_2\big(x;\bar{\Theta}_1[\Theta_2\bar{x}^{\Theta_2,v_2}+v_2](\cdot)\bar{x}^{\Theta_2,v_2}(\cdot)+\bar{v}_1[\Theta_2\bar{x}^{\Theta_2,v_2}+v_2](\cdot),\Theta_2(\cdot)\bar{x}^{\Theta_2,v_2}(\cdot)+v_2(\cdot)\big),\\
&\,\,\quad \quad\qquad\qquad\qquad\qquad\qquad\qquad\qquad \forall\, (\Theta_2(\cdot),v_2(\cdot))\in \mathcal{Q}_2[0,T] \times \mathcal{U}_2[0,T],
\end{aligned}
\end{equation}
where $\bar{x}(\cdot)\equiv\bar{x}^{\bar{\Theta}_2,\bar{v}_2}(\cdot)$ with $\bar{x}^{\Theta_2,v_2}(\cdot)\equiv x(\cdot;x,\bar{\Theta}_1[\Theta_2\bar{x}^{\Theta_2,v_2}+v_2](\cdot),\bar{v}_1[\Theta_2\bar{x}^{\Theta_2,v_2}+v_2](\cdot),\Theta_2(\cdot),v_2(\cdot))$ being the solution to the {\it closed-loop system} under the {\it closed-loop strategy} $(\Theta_2(\cdot),v_2(\cdot))$ for the leader.

\end{mydef}

\begin{Remark}
We can easily obtain the equation for $\bar{x}^{\Theta_2,v_2}(\cdot)$ in the above definition, by substituting $u_2(\cdot)$ with $\Theta_2(\cdot)\bar{x}^{\Theta_2,v_2}(\cdot)+v_2(\cdot)$ in (\ref{cse}), noting the dependence of $\bar{\Theta}_1(\cdot)$ and $\bar{v}_1(\cdot)$ on $(\Theta_2(\cdot),v_2(\cdot))$. We will give the details in Section 4, when dealing with the problem of the leader.
\end{Remark}

\section{LQ problem of the follower}

Let us introduce the following assumptions, which will be in force throughout this paper.
\par \textbf{(H1)} The coefficients of the state equation (\ref{state}) satisfy the following:
\begin{equation}\nonumber
\begin{cases}
A(\cdot) \in L^1(0,T;\mathbb{R}^{n \times n}),\,\,B_i(\cdot) \in L^2(0,T;\mathbb{R}^{n \times m_i}),\,\,\,\,\,b(\cdot) \in L^2_\mathbb{F}(\Omega;L^1(0,T;\mathbb{R}^n)),\\
C(\cdot) \in L^2(0,T;\mathbb{R}^{n \times n}),\,\,D_i(\cdot) \in L^{\infty}(0,T;\mathbb{R}^{n \times m_i}),\,\,\sigma(\cdot) \in L^2_\mathbb{F}(0,T;\mathbb{R}^n),\,\,\,\,\,\,\,i=1,2.
\end{cases}
\end{equation}
\par \textbf{(H2)} The weighting coefficients in the cost functional (\ref{cfi}) satisfy the following:
\begin{equation}\nonumber
\begin{cases}
Q^i(\cdot) \in L^1(0,T;\mathbb{S}^n),\,\,S^i_j(\cdot) \in L^2(0,T;\mathbb{R}^{m_i \times n}),\,\,R^i_{jj}(\cdot) \in L^{\infty}(0,T;\mathbb{S}^{m_j}),\\R^i_{12}(\cdot) \in L^{\infty}(0,T;\mathbb{S}^{m_1 \times m_2}),\,\,
R^i_{21}(\cdot) \in L^{\infty}(0,T;\mathbb{S}^{m_2 \times m_1}),\,\,q^i(\cdot)\in L^2_\mathbb{F}(\Omega;L^1(0,T;\mathbb{R}^n)),\\
\rho^i_j(\cdot)\in L^2_\mathbb{F}(0,T;\mathbb{R}^{m_j}),\,\,g^i \in L^2_{\mathcal{F}_T}(\Omega;\mathbb{R}^n),\,\,G^i \in \mathbb{S}^n,\quad i,j=1,2.
\end{cases}
\end{equation}

\textbf{Problem (SLQ)$_f$}. For any $x \in \mathbb{R}^n$, and given $u_2(\cdot) \in \mathcal{U}_2[0,T]$, find $\bar{u}_1(\cdot) \in \mathcal{U}_1[0,T]$ such that
\begin{equation}\label{follower problem}
J_1(x;\bar{u}_1(\cdot),u_2(\cdot))=\underset{u_1(\cdot)\in \mathcal{U}_1[0,T]} {\min}J_1(x;u_1(\cdot),u_2(\cdot)) \equiv V_1(x;u_2(\cdot)).
\end{equation}

It is worth noting that the strategy of the follower depends on the choice of the leader, that is, $u_1(\cdot)$ is related to $u_2(\cdot)$. If there exist a (unique) $\bar{u}_1(\cdot) \in \mathcal{U}_1[0,T]$ such that (\ref{follower problem}) holds, then we say that Problem (SLQ)$_f$ is (uniquely) \textit{open-loop solvable}. $\bar{u}_1(\cdot)$ is called an \textit{open-loop optimal control}, and $(\bar{x}^{u_2}(\cdot),\bar{u}_1(\cdot))\equiv(x(\cdot;x,\bar{u}_1(\cdot),u_2(\cdot)),\bar{u}_1(\cdot))$ is called an \textit{open-loop optimal pair}, of Problem (SLQ)$_f$.

First, using the idea of Proposition \ref{opsn}, we are able to obtain the following result.

%follower's open-loop optimal
\begin{mypro}
Let (H1)-(H2) hold. For a given $x \in \mathbb{R}^n$ and $u_2(\cdot) \in \mathcal{U}_2[0,T]$, a state-control pair $(\bar{x}^{u_2}(\cdot),\bar{u}_1(\cdot))$ is an open-loop optimal pair of Problem (SLQ)$_f$ if and only if the following holds:
\begin{equation}\label{fsc}
B_1^\top \bar{y}+D^\top_1\bar{z}+S^1_1\bar{x}^{u_2}+R^1_{11}\bar{u}_1+R^1_{12}u_2+\rho^1_1=0,\quad a.e.\ s\in[0,T],\ \mathbb{P}\mbox{-}a.s.,
\end{equation}
where $(\bar{y}(\cdot),\bar{z}(\cdot)) \in L^2_{\mathbb{F}}(0,T;\mathbb{R}^n) \times L^2_{\mathbb{F}}(0,T;\mathbb{R}^n)$ is the adapted solution to the following BSDE:
\begin{equation}\left\{
\begin{aligned}
d\bar{y}(s)&=-\big[A(s)^\top \bar{y}(s)+C(s)^\top\bar{z}(s)+Q^1(s)\bar{x}^{u_2}(s)+S^1_1(s)^\top\bar{u}_1(s)\\
           &\qquad +S^1_2(s)^\top u_2(s)+q^1(s)\big]ds+\bar{z}(s)dW(s),\,\,\quad\,s \in [0,T].\\
 \bar{y}(T)&=G^1\bar{x}^{u_2}(T)+g^1,
\end{aligned}\right.
\end{equation}
and the following convexity condition holds:
\begin{equation}
\mathbb{E}\bigg\{\int_0^T \Big[\big\langle Q^1x_0,x_0 \big\rangle +2\big\langle S^1_1x_0,u_1 \big\rangle +\big\langle R^1_{11}u_1,u_1 \big\rangle \Big]ds +\big\langle G^1x_0(T),x_0(T)\big\rangle\bigg\} \geqslant 0,
\end{equation}
where $x_0(\cdot) \in L^2_{\mathbb{F}}(0,T;\mathbb{R}^n)$ is the adapted solution to the following SDE:
\begin{equation}\left\{
\begin{aligned}
dx_0(s)&=\big[A(s)x_0(s)+B_1(s)u_1(s)\big]ds+\big[C(s)x_0(s)+D_1(s)u_1(s)\big]dW(s),\quad s \in [0,T],\\
 x_0(0)&=0.
\end{aligned}\right.
\end{equation}
\end{mypro}

Next, take $\Theta_1(\cdot) \in \mathcal{Q}_1[0,T]$ and $v_1(\cdot)\in \mathcal{U}_1[0,T]$. For any $x \in \mathbb{R}^n$ and $u_2(\cdot)\in \mathcal{U}_2[0,T]$, let us consider the closed-loop system (\ref{cse}) and the  corresponding cost functional (\ref{fccf}). The following result characterizes the closed-loop solvability of Problem (SLQ)$_f$.

%follower closed-loop optimal
\begin{mythm}\label{Th-cl-f}
Let (H1)-(H2) hold. Then Problem (SLQ)$_f$ admits a closed-loop optimal strategy on $[0,T]$ if and only if the following Riccati equation admits a solution $P^1(\cdot) \in C([0,T];\mathbb{S}^n)$:
\begin{equation}\label{follower Riccati}
\begin{cases}
\dot{P}^1+P^1A+A^\top P^1+C^\top P^1C+Q^1\\
\quad -(P^1B_1+C^\top P^1D_1+S^{1\top}_1)(R^1_{11}+D^\top_1P^1D_1)^{\dagger}(B_1^\top P^1 +D^\top_1P^1C+S^1_1)=0,\\
\mathcal{R}(B_1^\top P^1 +D^\top_1P^1C+S^1_1) \subseteq \mathcal{R}(R^1_{11}+D^\top_1P^1D_1),\\
R^1_{11}+D^\top_1P^1D_1 \geqslant 0,\\
P^1(T)=G^1,\\
\end{cases}
\end{equation}	
such that
\begin{equation}\label{range1}
(R^1_{11}+D^\top_1P^1D_1)^{\dagger}(B_1^\top P^1 +D^\top_1P^1C+S^1_1) \in L^2(0,T;\mathbb{R}^{m_1 \times n}),
\end{equation}
and the following BSDE admits an adapted solution $(\eta^1(\cdot),\zeta^1(\cdot)) \in L^2_{\mathbb{F}}(0,T;\mathbb{R}^n) \times L^2_{\mathbb{F}}(0,T;\mathbb{R}^n)$:
\begin{equation}\label{follower BSDE}\left\{
\begin{aligned}
d\eta^1&=-\Big\{ \big[A^\top-(P^1B_1+C^\top P^1D_1+S^{1\top}_1)(R^1_{11}+D^\top_1P^1D_1)^{\dagger}B_1^\top\big]\eta^1\\
&\qquad+\big[C^\top-(P^1B_1+C^\top P^1D_1+S^{1\top}_1)(R^1_{11}+D^\top_1P^1D_1)^{\dagger}D_1^\top\big]\zeta^1\\
&\qquad+\big[C^\top P^1D_2+S^{1\top}_2+P^1B_2-(P^1B_1+C^\top P^1D_1+S^{1\top}_1)(R^1_{11}+D^\top_1P^1D_1)^{\dagger}\\
&\qquad \times(R^1_{12}+D^\top_1P^1 D_2)\big]u_2-(P^1B_1+C^\top P^1D_1+S^{1\top}_1)(R^1_{11}+D^\top_1P^1D_1)^{\dagger}\\
&\qquad \times(D_1^\top P^1\sigma+\rho^1_1)+C^\top P^1\sigma+q^1+P^1 b \Big\}ds+\zeta^1dW,\quad s \in [0,T],\\
&\hspace{-8mm}B^\top_1\eta^1+D^\top_1\zeta^1+(R^1_{12}+D^\top_1P^1D_2)u_2+D^\top_1P^1\sigma+\rho^1_1\\
&\qquad \in \mathcal{R}(R^1_{11}+D^\top_1P^1D_1),\qquad a.e.\,\,s \in [0,T],\,\,\mathbb{P}\mbox{-}a.s.,\\
\eta^1(T)&=g^1,
\end{aligned}\right.
\end{equation}	
which satisfies
\begin{equation}\label{range2}
[R^1_{11}+D^\top_1P^1D_1]^{\dagger}\big[B^\top_1\eta^1+D^\top_1\zeta^1+(R^1_{12}+D^\top_1P^1D_2)u_2+D^\top_1P^1\sigma+\rho^1_1\big] \in L^2_\mathbb{F}(0,T;\mathbb{R}^{m_1}).
\end{equation}
In this case, the closed-loop optimal strategy $(\bar{\Theta}_1(\cdot),\bar{v}_1(\cdot))$ of Problem (SLQ)$_f$ admits the following representation:
\begin{equation}\label{follower closed-loop optimal}\left\{
\begin{aligned}
\bar{\Theta}_1&=-(R^1_{11}+D^\top_1P^1D_1)^{\dagger}(B_1^\top P^1 +D^\top_1P^1C+S^1_1)\\
              &\qquad\quad  +[I-(R^1_{11}+D^\top_1P^1D_1)^\dagger(R^1_{11}+D^\top_1P^1D_1)]\theta_1,\\
     \bar{v}_1&=-(R^1_{11}+D^\top_1P^1D_1)^{\dagger}(B^\top_1\eta^1+D^\top_1\zeta^1+(R^1_{12}+D^\top_1P^1D_2)u_2+D^\top_1P^1\sigma+\rho^1_1)\\
              &\qquad\quad  +[I-(R^1_{11}+D^\top_1P^1D_1)^\dagger(R^1_{11}+D^\top_1P^1D_1)]v,
\end{aligned}\right.
\end{equation}
for some $\theta_1(\cdot) \in \mathcal{Q}_1[0,T]$ and $v(\cdot)\in \mathcal{U}_1[0,T]$. Further, the value function $V_1(\cdot;u_2(\cdot))$ is given by
\begin{equation}\label{value-1}
\begin{split}
&V_1(x;u_2(\cdot))=\underset{u_1(\cdot)\in \mathcal{U}_1[0,T]} {\inf}J_1(x;u_1(\cdot),u_2(\cdot))=\mathbb{E}\biggl\{ \big\langle P^1(0)x,x\big\rangle +2\big\langle \eta^1(0),x\big\rangle \\
      & =\mathbb{E}\biggl\{ \big\langle P^1(0)x,x\big\rangle +2\big\langle \eta^1(0),x\big\rangle +\int_0^T\Big[\big\langle (R^1_{22}+D^\top_2P^1D_2)u_2,u_2\big\rangle\\
      &\qquad +2\big\langle \rho^1_2+D^\top_2P^1\sigma+B^\top_2\eta^1+D^\top_2\zeta^1,u_2 \big\rangle+2\langle \eta^1,b\rangle+2\langle \zeta^1,\sigma \rangle+\langle P^1\sigma,\sigma\rangle\\
      &\qquad-\big|\big[(R^1_{11}+D^\top_1P^1D_1)^\dagger\big]^{\frac{1}{2}}\big[B^\top_1\eta^1+D^\top_1\zeta^1+(R^1_{12}+D^\top_1P^1D_2)u_2+D^\top_1P^1\sigma+\rho^1_1\big]\big|^2\Big]ds \biggr\}.
\end{split}
\end{equation}
\end{mythm}
The proof here is similar to that in \cite{SLY2016,SY2014}, but for the sake of the integrity of the article, we still give the proof.

{\it Proof.} We first prove the necessity. Given $u_2(\cdot)\in\mathcal{U}_2[0,T]$. Let $(\bar{\Theta}_1(\cdot),\bar{v}_1(\cdot)) \in \mathcal{Q}_1[0,T] \times \mathcal{U}_1[0,T]$ be a closed-loop optimal strategy of Problem (SLQ)$_f$ over $[0,T]$. Then, by Proposition \ref{relation}, $\bar{v}_1(\cdot)$ is an open-loop optimal control of the following linear stochastic system:
\begin{equation*}\left\{
\begin{aligned}
dx^{u_2}(s)&=\big\{\big[A(s)+B_1(s)\bar{\Theta}_1(s)\big]x^{u_2}(s)+B_1(s)v_1(s)+B_2(s)u_2(s)+b(s)\big\}ds\\
           &\quad+\big\{\big[C(s)+D_1(s)\bar{\Theta}_1(s)\big]x^{u_2}(s)+D_1(s)v_1(s)+D_2(s)u_2(s)+\sigma(s)\big\}dW(s),\\
           &\hspace{7cm} s \in [0,T],\\
 x^{u_2}(0)&=x,
\end{aligned}\right.
\end{equation*}
with the cost functional
\begin{equation*}
\begin{aligned}
&J_1(x;\bar{\Theta}x^{u_2}+v_1,u_2)=\mathbb{E}\bigg\{\int_0^T \Big[\big\langle (Q^1+\bar{\Theta}_1^\top S^1_1+S^{1\top}_1\bar{\Theta}_1+\bar{\Theta}_1^\top R^1_{11}\bar{\Theta}_1)x^{u_2},x^{u_2} \big\rangle\\
   &\quad +2\big\langle (S^1_1+R^1_{11}\bar{\Theta}_1)x^{u_2},v_1 \big\rangle+\langle R^1_{11}v_1,v_1 \rangle +2\big\langle S^{1\top}_2u_2+q^1+\bar{\Theta}_1^\top R^1_{12}u_2+\bar{\Theta}_1^\top \rho^1_1,x^{u_2} \big\rangle\\
   &\quad +2\big\langle R^1_{12}u_2 +\rho^1_1,v_1 \big\rangle+2\langle \rho^1_2,u_2 \rangle+\langle R^1_{22}u_2,u_2 \rangle \Big]ds
    +\big\langle G^1x^{u_2}(T),x^{u_2}(T)\big\rangle +2\big\langle g^1 ,x^{u_2}(T)\big\rangle \bigg\}.
\end{aligned}
\end{equation*}
Hence, by Proposition \ref{opsn}, for any $x \in \mathbb{R}^n$, the following FBSDE admits an adapted solution $(\bar{x}^{u_2}(\cdot),\bar{y}^{u_2}(\cdot),\bar{z}^{u_2}(\cdot)) \in L^2_{\mathbb{F}}(0,T;\mathbb{R}^n)\times L^2_{\mathbb{F}}(0,T;\mathbb{R}^n)\times L^2_{\mathbb{F}}(0,T;\mathbb{R}^n)$:
\begin{equation}\label{follower optimal FBSDE}
\left\{
\begin{aligned}
  d\bar{x}^{u_2}&=\big[(A+B_1\bar{\Theta}_1)\bar{x}^{u_2}+B_1\bar{v}_1+B_2u_2+b\big]ds\\
                &\qquad+\big[(C+D_1\bar{\Theta}_1)\bar{x}^{u_2}+D_1\bar{v}_1+D_2u_2+\sigma\big]dW(s),\\
 -d\bar{y}^{u_2}&=\big[(A+B_1\bar{\Theta}_1)^\top\bar{y}^{u_2}+(C+D_1\bar{\Theta}_1)^\top\bar{z}^{u_2}+(Q^1+\bar{\Theta}^\top_1S^1_1+S^{1\top}_1\bar{\Theta}_1+\bar{\Theta}_1^\top R^1_{11}\bar{\Theta}_1)\bar{x}^{u_2}\\
                &\qquad +(S^1_1+R^1_{11}\bar{\Theta}_1)^\top\bar{v}_1+S^{1\top}_2u_2+q^1+\bar{\Theta}^\top_1R^1_{12}u_2+\bar{\Theta}^\top_1\rho^1_1\big]ds-\bar{z}^{u_2}dW,\\
\bar{x}^{u_2}(0)&=x,\quad\bar{y}^{u_2}(T)=G^1\bar{x}^{u_2}(T)+g^1,
\end{aligned}\right.
\end{equation}
and the following stationarity condition holds:
\begin{equation}\label{ follower stationarity condition}
B_1^\top \bar{y}+D^\top_1\bar{z}+(S^1_1+R^1_{11}\bar{\Theta}_1)\bar{x}+R^1_{11}\bar{v}_1+R^1_{12}u_2+\rho^1_1=0,\quad a.e.,\,\,\mathbb{P}\mbox{-}a.s.
\end{equation}
Making use of (\ref{ follower stationarity condition}), we may rewrite the BSDE in (\ref{follower optimal FBSDE}) as follows:
\begin{equation}\nonumber
\begin{aligned}
d\bar{y}^{u_2}&=-\big[A^\top\bar{y}+C^\top\bar{z}^{u_2}+(Q^1+S^{1\top}_1\bar{\Theta}_1)\bar{x}^{u_2}+S^{1\top}_1\bar{v}_1+S^{1\top}_2u_2+q^1\\
        &\qquad\quad+\bar{\Theta}_1^\top[B_1^\top\bar{y}^{u_2}+D^\top_1\bar{z}^{u_2}+(S^1_1+R^1_{11}\bar{\Theta}_1)\bar{x}^{u_2}+R^1_{11}\bar{v}_1+R^1_{12}u_2+\rho^1_1]\big]ds+\bar{z}^{u_2}dW\\
        &=-\big[A^\top\bar{y}^{u_2}+C^\top\bar{z}^{u_2}+(Q^1+S^{1\top}_1\bar{\Theta}_1)\bar{x}^{u_2}+S^{1\top}_1\bar{v}_1+S^{1\top}_2u_2+q^1\big]ds+\bar{z}^{u_2}dW.
\end{aligned}
\end{equation}
Thus, we obtain
\begin{equation}\label{follower simple FBSDE}
\left\{
\begin{aligned}
  d\bar{x}^{u_2}&=\big[(A+B_1\bar{\Theta}_1)\bar{x}^{u_2}+B_1\bar{v}_1+B_2u_2+b\big]ds\\
                &\qquad+\big[(C+D_1\bar{\Theta}_1)\bar{x}^{u_2}+D_1\bar{v}_1+D_2u_2+\sigma\big]dW(s),\\
  d\bar{y}^{u_2}&=-\big[A^\top\bar{y}^{u_2}+C^\top\bar{z}^{u_2}+(Q^1+S^{1\top}_1\bar{\Theta}_1)\bar{x}^{u_2}+S^{1\top}_1\bar{v}_1+S^{1\top}_2u_2+q^1\big]ds+\bar{z}^{u_2}dW,\\
\bar{x}^{u_2}(0)&=x,\quad\bar{y}^{u_2}(T)=G^1\bar{x}^{u_2}(T)+g^1,\\
                &\hspace{-8mm}B_1^\top\bar{y}+D^\top_1\bar{z}+(S^1_1+R^1_{11}\bar{\Theta}_1)\bar{x}+R^1_{11}\bar{v}_1+R^1_{12}u_2+\rho^1_1=0,\quad a.e.,\,\,\,\mathbb{P}\mbox{-}a.s.
\end{aligned}\right.
\end{equation}
Since the above admits a solution for each $x \in \mathbb{R}^n$, and $(\bar{\Theta}_1(\cdot),\bar{v}_1(\cdot))$ is independent of $x$, by subtraction solutions corresponding $x$ and $0$, the later from the former, we see that for any $x \in \mathbb{R}^n$, the following FBSDE admits an adapted solution $(x(\cdot),y(\cdot),z(\cdot)) \in L^2_{\mathbb{F}}(0,T;\mathbb{R}^n)\times L^2_{\mathbb{F}}(0,T;\mathbb{R}^n)\times L^2_{\mathbb{F}}(0,T;\mathbb{R}^n)$:
\begin{equation}\left\{
\begin{aligned}
dx(s)&=\big[A(s)+B_1(s)\bar{\Theta}_1(s)\big]x(s)ds+\big[C(s)+D_1(s)\bar{\Theta}_1(s)\big]x(s)dW(s),\\
dy(s)&=-\big\{A^\top(s)y(s)+C^\top(s)z(s)+\big[Q^1(s)+S^{1\top}_1(s)\bar{\Theta}_1(s)\big]x(s)\big\}ds+z(s)dW,\\
 x(0)&=x,\quad y(T)=G^1x(T),\\
     &\hspace{-8mm}B_1^\top(s)y(s)+D^\top_1(s)z(s)+\big[S^1_1(s)+R^1_{11}(s)\bar{\Theta}_1(s)\big]x(s)=0,\ a.e.\,s \in [0,T],\,\,a.s.\,\,\mathbb{P}\mbox{-}a.s.
\end{aligned}\right.
\end{equation}
Now, we let
\begin{equation*}\left\{
\begin{aligned}
d\mathbb{X}(s)&=\big[A(s)+B_1(s)\bar{\Theta}_1(s)\big]\mathbb{X}(s)ds+\big[C(s)+D_1(s)\bar{\Theta}_1(s)\big]\mathbb{X}(s)dW(s),\quad s \in [0,T],\\
\mathbb{X}(0)&=I_{n\times n},
\end{aligned}\right.
\end{equation*}
and let
\begin{equation*}\left\{
\begin{aligned}
d\mathbb{Y}(s)&=-\big\{A^\top(s)\mathbb{Y}(s)+C^\top(s)\mathbb{Z}(s)+\big[Q^1(s)+S^{1\top}_1(s)\bar{\Theta}_1(s)\big]\mathbb{X}(s)\big\}ds\\
              &\qquad +\mathbb{Z}(s)dW,\quad s \in [0,T],\\
\mathbb{Y}(T)&=G^1\mathbb{X}(T).
\end{aligned}\right.
\end{equation*}
Clearly, $\mathbb{X}(\cdot),\,\mathbb{Y}(\cdot),\,\mathbb{Z}(\cdot)$ are all well-defined $\mathbb{S}^n$-matrix valued processes. Further,
\begin{equation}\label{ssc}
B_1^\top(s) \mathbb{Y}(s)+D^\top_1(s)\mathbb{Z}(s)+\big[S^1_1(s)+R^1_{11}(s)\bar{\Theta}_1(s)\big]\mathbb{X}(s)=0,\quad a.e.\,s \in [0,T],\,\,\mathbb{P}\mbox{-}a.s.
\end{equation}
And $\mathbb{X}(\cdot)^{-1}$ exists, which satisfies the following SDE:
\begin{equation}\left\{
\begin{aligned}
d\mathbb{X}^{-1}(s)&=\mathbb{X}^{-1}(s)\big\{\big[C(s)+D_1(s)\bar{\Theta}_1(s)\big]^2-A(s)-B_1(s)\bar{\Theta}_1(s)\big\}ds\\
                   &\quad -\mathbb{X}^{-1}(s)\big[C(s)+D_1(s)\bar{\Theta}_1(s)\big]dW(s),\quad s \in [0,T],\\
 \mathbb{X}^{-1}(0)&=I_{n\times n}.
\end{aligned}\right.
\end{equation}
We define
\begin{equation}\nonumber
P^1(\cdot)=\mathbb{Y}(\cdot)\mathbb{X}(\cdot)^{-1},\qquad \Pi^1(\cdot)=\mathbb{Z}(\cdot)\mathbb{X}(\cdot)^{-1}.
\end{equation}
Then (\ref{ssc}) implies
\begin{equation}\label{Theta condition}
B_1^\top P^1+D^\top_1\Pi^1+(S^1_1+R^1_{11}\bar{\Theta}_1)=0,\quad a.e.,\quad \mathbb{P}\mbox{-}a.s.
\end{equation}
Also, by It\^{o}'s formula, we get
\begin{equation}\nonumber
\begin{aligned}
dP^1&=d\mathbb{Y}\mathbb{X}^{-1}=d\mathbb{Y}\cdot\mathbb{X}^{-1}+\mathbb{Y}d\mathbb{X}^{-1}+d\mathbb{Y}\cdot d\mathbb{X}^{-1}\\
    &=\bigl\{-A^\top P^1-C^\top \Pi^1-Q^1-S^{1\top}_1\bar{\Theta}_1+P^1\big[(C+D_1\bar{\Theta}_1)^2-A-B_1\bar{\Theta}_1\big]\\
    &\qquad -\Pi^1(C+D_1\bar{\Theta}_1)\bigr\}ds+\big[\Pi^1-P^1(C+D_1\bar{\Theta}_1)\big]dW.
\end{aligned}
\end{equation}
Let $\Lambda=\Pi^1-P^1(C+D_1\bar{\Theta}_1)$, then
\begin{equation}\nonumber
\begin{aligned}
dP^1&=\bigl\{-A^\top P^1-C^\top\big[\Lambda+P^1(C+D_1\bar{\Theta}_1)\big]-Q^1-S^{1\top}_1\bar{\Theta}_1+P^1\big[(C+D_1\bar{\Theta}_1)^2\\
    &\qquad -A-B_1\bar{\Theta}_1\big]-\big[\Lambda+P^1(C+D_1\bar{\Theta}_1)\big](C+D_1\bar{\Theta}_1)\bigr\}ds+\Lambda dW\\
    &=-\bigl\{A^\top P^1+P^1A+C^\top\Lambda+\Lambda C+C^\top P^1C+Q^1\\
    &\qquad +(P^1B_1+C^\top P^1D_1+S^{1\top}_1+\Lambda D_1)\bar{\Theta}_1\bigr\}ds+\Lambda dW,
\end{aligned}
\end{equation}
and $P^1(T)=G_1$. Thus, $(P^1(\cdot),\Lambda(\cdot))$ is an adapted solution to a BSDE with deterministic coefficients. Hence, $P^1(\cdot)$ is deterministic and $\Lambda(\cdot)=0$ which means
\begin{equation}\label{Pi}
\Pi^1=P^1(C+D_1\bar{\Theta}_1).
\end{equation}
Therefore,
\begin{equation}\label{follower RE---}
\dot{P}^1+A^\top P^1+P^1A+C^\top P^1C+(P^1B_1+C^\top P^1D_1+S^{1\top}_1)\bar{\Theta}_1+Q^1=0.
\end{equation}
Using (\ref{Pi}), (\ref{Theta condition}) can be written as
\begin{equation}\label{Theta_1}
\begin{aligned}
0&=B_1^\top P^1+D^\top_1P^1(C+D_1\bar{\Theta}_1)+S^1_1+R^1_{11}\bar{\Theta}_1\\
 &=B_1^\top P^1+D^\top_1P^1C+S^1_1+(R^1_{11}+D^\top_1P^1D_1)\bar{\Theta}_1,\qquad a.e.
\end{aligned}
\end{equation}
This implies
\begin{equation}
\mathcal{R}(B_1^\top P^1+D^\top_1P^1C+S^1_1) \subseteq \mathcal{R}(R^1_{11}+D^\top_1P^1D_1),\qquad a.e.
\end{equation}
Using (\ref{Theta_1}), (\ref{follower RE---}) can be written as
\begin{equation}\nonumber
\begin{aligned}
0&=\dot{P}^1+A^\top P^1+P^1A+C^\top P^1C+(C^\top   P^1D_1+S^{1\top}_1+P^1B_1)\bar{\Theta}_1\\
 &\qquad  +Q^1+\bar{\Theta}^\top_1\big[B_1^\top P^1+D^\top_1P^1C+S^1_1+(R^1_{11}+D^\top_1P^1D_1)\bar{\Theta}_1\big]\\
 &=\dot{P}^1+(A+B_1\bar{\Theta}_1)^\top P^1+P^1(A+B_1\bar{\Theta}_1)+(C+D_1\bar{\Theta}_1)^\top P^1(C+D_1\bar{\Theta}_1)\\
 &\qquad+S^{1\top}_1\bar{\Theta}_1+\bar{\Theta}^\top_1S^1_1+ +Q^1+\bar{\Theta}^\top_1R^1_{11}\bar{\Theta}_1,\qquad a.e.
\end{aligned}
\end{equation}

Since $P^1(T)=G^1 \in \mathbb{S}^n$ and $Q^1(\cdot)$, $R^1_{11}(\cdot)$ are symmetric, by uniqueness, we must have $P^1(\cdot) \in C([0,T];\mathbb{S}^n)$. Denoting $\hat{R}_{11}\equiv R^1_{11}+D^\top_1P^1D_1$, since
\begin{equation}\nonumber
\hat{R}^\dagger_{11}(B_1^\top P^1+D^\top_1P^1C+S^1_1)=-\hat{R}^\dagger_{11}\hat{R}_{11}\bar{\Theta}_{11}
\end{equation}
and $\hat{R}^\dagger_{11}\hat{R}_{11}$ is an orthogonal projection, we see that (\ref{range1}) holds and
\begin{equation}\nonumber
\bar{\Theta}_1=-\hat{R}^\dagger_{11}(B_1^\top P^1+D^\top_1P^1C+S^1_1)+(I-\hat{R}^\dagger_{11}\hat{R}_{11})\theta
\end{equation}
for some $\theta(\cdot) \in L^2(0,T;\mathbb{R}^{m_1 \times n})$. Consequently,
\begin{equation}
\begin{aligned}
&(P^1B_1+C^\top P^1D_1+S^{1\top}_1)\bar{\Theta}_1=\bar{\Theta}^\top_1\hat{R}_{11}\hat{R}^\dagger_{11}(B_1^\top P^1+D^\top_1P^1C+S^1_1)\\
&\qquad =-(P^1B_1+C^\top P^1D_1+S^{1\top}_1)\hat{R}^\dagger_{11}(B_1^\top P^1+D^\top_1P^1C+S^1_1).
\end{aligned}
\end{equation}
Plug the above into (\ref{follower RE---}), we obtain Riccati equation in (\ref{follower Riccati}). To determine $\bar{v}_1(\cdot)$, we define
\begin{equation}\nonumber
\begin{cases}
\eta^1=\bar{y}^{u_2}-P^1\bar{x}^{u_2},\\
\zeta^1=\bar{z}^{u_2}-P^1(C+D_1\bar{\Theta}_1)\bar{x}^{u_2}-P^1D_1\bar{v}_1-P^1(D_2u_2+\sigma).
\end{cases}
\end{equation}
Then (noting (\ref{follower simple FBSDE}), (\ref{follower RE---}) and (\ref{Theta_1}))
\begin{equation}\nonumber
\begin{aligned}
d\eta^1&=d\bar{y}^{u_2}-\dot{P}^1\bar{x}^{u_2}ds-P^1d\bar{x}^{u_2}\\
       &=\big\{-A^\top\bar{y}^{u_2}-C^\top\bar{z}^{u_2}-(Q^1+S^{1\top}_1\bar{\Theta}_1)\bar{x}^{u_2}-S^{1\top}_1\bar{v}_1-S^{1\top}_2u_2-q^1\\
       &\qquad +A^\top P^1\bar{x}^{u_2}+P^1A\bar{x}^{u_2}+C^\top P^1C\bar{x}^{u_2}+(C^\top P^1D_1+S^{1\top}_1+P^1B_1)\bar{\Theta}_1\bar{x}^{u_2}\\
       &\qquad +Q^1\bar{x}^{u_2}-P^1A\bar{x}^{u_2}-P^1B_1\bar{\Theta}_1\bar{x}^{u_2}-P^1B_1\bar{v}_1-P^1B_2u_2-P^1b\bigr\}ds\\
       &\quad +\big\{\bar{z}^{u_2}-P^1C\bar{x}^{u_2}-P^1D_1\bar{\Theta}_1\bar{x}^{u_2}-P^1D_1\bar{v}_1-P^1D_2u_2-P^1\sigma\big\}dW\\
       &=\big\{-A^\top(\eta^1+P^1\bar{x}^{u_2})-C^\top\big[\zeta^1+P^1(C+D_1\bar{\Theta}_1)\bar{x}^{u_2}+P^1D_1\bar{v}_1+P^1(D_2u_2+\sigma)\big]\\
       &\qquad -(Q^1+S^{1\top}_1\bar{\Theta}_1)\bar{x}^{u_2}-S^{1\top}_1\bar{v}_1-S^{1\top}_2u_2-q^1+A^\top P^1\bar{x}^{u_2}+P^1A\bar{x}^{u_2}\\
       &\qquad +C^\top P^1C\bar{x}^{u_2}+(C^\top P^1D_1+S^{1\top}_1+P^1B_1)\bar{\Theta}_1\bar{x}^{u_2}+Q^1\bar{x}^{u_2}\\
       &\qquad -P^1A\bar{x}^{u_2}-P^1B_1\bar{\Theta}_1\bar{x}^{u_2}-P^1B_1\bar{v}_1-P^1B_2u_2-P^1b\bigr\}ds+\zeta^1dW\\
       &=-\big\{A^\top\eta^1+C^\top\zeta^1+(C^\top P^1D_1+P^1B_1+S^{1\top}_1)\bar{v}_1+(C^\top P^1D_2+S^{1\top}_2+P^1B_2)u_2\\
       &\qquad +C^\top P^1\sigma+q^1+P^1b\bigr\}ds+\zeta^1dW.
\end{aligned}
\end{equation}
According to (\ref{ follower stationarity condition}), we have
\begin{equation}\label{v_1 condition}
\begin{aligned}
0&=B_1^\top \bar{y}^{u_2}+D^\top_1\bar{z}^{u_2}+(S^1_1+R^1_{11}\bar{\Theta}_1)\bar{x}^{u_2}+R^1_{11}\bar{v}_1+R^1_{12}u_2+\rho^1_1\\
 &=B_1^\top (\eta^1+P^1\bar{x}^{u_2})+D^\top_1\big[\zeta^1+P^1(C+D_1\bar{\Theta}_1)\bar{x}^{u_2}+P^1D_1\bar{v}_1+P^1(D_2u_2+\sigma)\big]\\
 &\qquad +(S^1_1+R^1_{11}\bar{\Theta}_1)\bar{x}^{u_2}+R^1_{11}\bar{v}_1+R^1_{12}u_2+\rho^1_1\\
 &=B_1^\top\eta^1+D^\top_1\zeta^1+(R^1_{11}+D^\top_1P^1D_1)\bar{v}_1+(R^1_{12}+D^\top_1P^1D_2)u_2+D^\top_1P^1\sigma+\rho^1_1\\
 &\qquad +\big[B_1^\top P^1+D^\top_1P^1C+(R^1_{11}+D^\top_1P^1D_1)\bar{\Theta}_1+S^1_1\big]\bar{x}^{u_2}\\
 &=B_1^\top\eta^1+D^\top_1\zeta^1+(R^1_{11}+D^\top_1P^1D_1)\bar{v}_1+(R^1_{12}+D^\top_1P^1D_2)u_2+D^\top_1P^1\sigma+\rho^1_1.
\end{aligned}
\end{equation}
Hence,
\begin{equation}
B_1^\top\eta^1+D^\top_1\zeta^1+(R^1_{12}+D^\top_1P^1D_2)u_2+D^\top_1P^1\sigma+\rho^1_1 \in \mathcal{R}(R^1_{11}+D^\top_1P^1D_1),\quad a.e.,\,\,a.s.
\end{equation}
Since
\begin{equation}\nonumber
\hat{R}^\dagger_{11}[B_1^\top\eta^1+D^\top_1\zeta^1+(R^1_{12}+D^\top_1P^1D_2)u_2+D^\top_1P^1\sigma+\rho^1_1]=-\hat{R}^\dagger_{11}\hat{R}_{11}\bar{v}_1
\end{equation}
and $\hat{R}^\dagger_{11}\hat{R}_{11}$ is an orthogonal projection, we see that (\ref{range2}) holds and
\begin{equation}\nonumber
\bar{v}_1=-\hat{R}^\dagger_{11}[B_1^\top\eta^1+D^\top_1\zeta^1+(R^1_{12}+D^\top_1P^1D_2)u_2+D^\top_1P^1\sigma+\rho^1_1]+(I-\hat{R}^\dagger_{11}\hat{R}_{11})v
\end{equation}
for some $v(\cdot) \in L^2_\mathbb{F}(0,T;\mathbb{R}^{m_1})$.
Consequently,
\begin{equation}\nonumber
\begin{aligned}
&(C^\top P^1D_1+P^1B_1+S^{1\top}_1)\bar{v}_1\\
&=-(C^\top P^1D_1+P^1B_1+S^{1\top}_1)\hat{R}^\dagger_{11}\big[B_1^\top\eta^1+D^\top_1\zeta^1+(R^1_{12}+D^\top_1P^1D_2)u_2 +D^\top_1P^1\sigma+\rho^1_1\big]\\
&\qquad \qquad+(C^\top P^1D_1+P^1B_1+S^{1\top}_1)(I-\hat{R}^\dagger_{11}\hat{R}_{11})v\\
&=-(C^\top P^1D_1+P^1B_1+S^{1\top}_1)\hat{R}^\dagger_{11}\big[B_1^\top\eta^1+D^\top_1\zeta^1+(R^1_{12}+D^\top_1P^1D_2)u_2 +D^\top_1P^1\sigma+\rho^1_1\big].
\end{aligned}
\end{equation}
Therefore, $(\eta^1,\zeta^1)$ is the adapted solution to the following BSDE:
\begin{equation}\nonumber
\begin{cases}
d\eta^1=-\Big\{ \big[A^\top-(P^1B_1+C^\top P^1D_1+S^{1\top}_1)(R^1_{11}+D^\top_1P^1D_1)^{\dagger}B_1^\top\big]\eta^1\\
\qquad \qquad\quad +\big[C^\top-(P^1B_1+C^\top P^1D_1+S^{1\top}_1)(R^1_{11}+D^\top_1P^1D_1)^{\dagger}D_1^\top\big]\zeta^1\\
\qquad \qquad\quad +\big[C^\top P^1D_2+S^{1\top}_2+P^1B_2\\
\quad\qquad\quad \qquad -(P^1B_1+C^\top P^1D_1+S^{1\top}_1)(R^1_{11}+D^\top_1P^1D_1)^{\dagger}(R^1_{12}+D^\top_1P^1 D_2)\big]u_2\\
\qquad \qquad\quad +\big[C^\top-(P^1B_1+C^\top P^1D_1+S^{1\top}_1)(R^1_{11}+D^\top_1P^1D_1)^{\dagger}D_1^\top\big]P^1\sigma\\
\qquad \qquad\quad -(P^1B_1+C^\top P^1D_1+S^{1\top}_1)(R^1_{11}+D^\top_1P^1D_1)^{\dagger}\rho^1_1+q^1+P^1 b \Big\}ds+\zeta^1dW,\\
\eta^1(T)=g^1.
\end{cases}
\end{equation}

To prove $R^1_{11}+D^\top_1P^1D_1 \geqslant 0$, as well as the sufficiency, we take any $u_1(\cdot) \in \mathcal{U}_1[0,T]$, and let $x(\cdot)\equiv x(\cdot;x,u_1(\cdot),u_2(\cdot)),\,\bar{x}^{u_2}(\cdot)\equiv x(\cdot;x,\bar{\Theta}_1(\cdot),\bar{v}_1(\cdot),u_2(\cdot)) $ be the corresponding state processes. Then, by It\^o's formula, we have
\begin{equation}
\begin{aligned}
&J_1(x;u_1(\cdot),u_2(\cdot))=\mathbb{E}\biggl\{ \int_0^T\Big[\langle Q^1x,x \rangle+2\langle S^1_1x,u_1 \rangle +\langle R^1_{11}u_1,u_1 \rangle+2\langle R^1_{12}u_2+\rho^1_1,u_1 \rangle\\
&\quad +2\langle S^{1\top}_2u_2+q^1,x \rangle+2\langle\rho^1_2,u_2 \rangle+\langle R^1_{22}u_2,u_2 \rangle\Big] ds +\langle G^1x(T),x(T)\rangle+2\langle g^1,x(T) \rangle\biggr\}\\
&=J_1(x;\bar{\Theta}_1(\cdot)\bar{x}^{u_2}(\cdot)+\bar{v}_1(\cdot),u_2(\cdot))\\
&\quad +\mathbb{E}\int_0^T\big\langle (R^1_{11}+D^\top_1P^1D_1)(u_1-\bar{\Theta}_1x-\bar{v}_1),u_1-\bar{\Theta}_1x-\bar{v}_1\big\rangle ds.
\end{aligned}
\end{equation}
Hence,
\begin{equation}\nonumber
J_1(x;\bar{\Theta}_1(\cdot)\bar{x}^{u_2}(\cdot)+\bar{v}_1(\cdot),u_2(\cdot)) \leqslant J_1(x;\bar{\Theta}_1(\cdot)x(\cdot)+v_1(\cdot),u_2(\cdot)),\quad \forall v_1(\cdot) \in \mathcal{U}_1[0,T],
\end{equation}
if and only if
\begin{equation}\nonumber
R^1_{11}+D^\top_1P^1D_1 \geqslant 0,\quad a.e.
\end{equation}
In this case,
\begin{equation}\nonumber
\begin{split}
&V_1(x;u_2(\cdot))=\mathbb{E}\biggl\{ \big\langle P^1(0)x,x\big\rangle +2\big\langle \eta^1(0),x\big\rangle +\int_0^T\Big[ \big\langle (R^1_{22}+D^\top_2P^1D_2)u_2,u_2 \big\rangle\\
      &\qquad +\langle P^1\sigma,\sigma\rangle+2\langle\eta^1,b\rangle+2\langle\zeta^1,\sigma\rangle+2\big\langle B^\top_2\eta^1+D^\top_2\zeta^1+D^\top_2P^1\sigma+\rho^1_2,u_2\big\rangle\\
      &\qquad -\big|\big[(R^1_{11}+D^\top_1P^1D_1)^\dagger\big]^{\frac{1}{2}}\big[B_1^\top\eta^1+D^\top_1\zeta^1+(R^1_{12}+D^\top_1P^1D_2)u_2+D^\top_1P^1\sigma+\rho^1_1\big]\big|^2ds \biggr\}
\end{split}
\end{equation}
which is (\ref{value-1}). This completes the proof. $\qquad\Box$

\section{LQ problem of the leader}

Now, let Problem (SLQ)$_f$ be uniquely closed-loop solvable for given $(x,u_2(\cdot)) \in \mathbb{R}^n \times \mathcal{U}_2[0,T]$. For the sake of brevity, we assume that $(R^1_{11}+D^\top_1P^1D_1)^{-1}$ exists in the following. If $R^1_{11}+D^\top_1P^1D_1$ is not invertible, the concept of generalized inverse could be introduced as before, and the expression becomes more complicated, but there is no other difficulty. Then by (\ref{follower closed-loop optimal}), the follower takes his/her following optimal control:
\begin{equation}\label{follower optimal control}
\begin{aligned}
\bar{u}_1(t)&=\bar{\Theta}(t)\bar{x}^{u_2}(t)+\bar{v}_1(t)=-(R^1_{11}+D^\top_1P^1D_1)^{-1}\big[\big(B_1^\top P^1 +D^\top_1P^1C+S^1_1\big)\bar{x}^{u_2}(t)\\
            &\quad  +B^\top_1\eta^1+D^\top_1\zeta^1+(R^1_{12}+D^\top_1P^1D_2)u_2+D^\top_1P^1\sigma+\rho^1_1\big],\quad t \in [0,T],
\end{aligned}
\end{equation}
where the process triple $(\bar{x}^{u_2}(\cdot),\eta^{1,u_2}(\cdot),\zeta^{1,u_2}(\cdot))\in L^2_{\mathbb{F}}(0,T;\mathbb{R}^n)\times L^2_{\mathbb{F}}(0,T;\mathbb{R}^n)\times L^2_{\mathbb{F}}(0,T;\mathbb{R}^n)$ satisfies the following FBSDE, which, now, is the ``state" equation of the leader:
\begin{equation}\nonumber
\left\{
\begin{aligned}
  d\bar{x}^{u_2}&=\Big\{\big[A-B_1(R^1_{11}+D^\top_1P^1D_1)^{-1}(B_1^\top P^1 +D^\top_1P^1C+S^1_1)\big]\bar{x}^{u_2}\\
                &\qquad+\big[B_2-B_1(R^1_{11}+D^\top_1P^1D_1)^{-1}(R^1_{12}+D^\top_1P^1D_2)\big]u_2\\
                &\qquad -B_1(R^1_{11}+D^\top_1P^1D_1)^{-1}(B^\top_1\eta^{1,u_2}+D^\top_1\zeta^{1,u_2}+D^\top_1P^1\sigma+\rho^1_1)+b \Big\}ds\\
                &\quad +\Big\{\big[C-D_1(R^1_{11}+D^\top_1P^1D_1)^{-1}(B_1^\top P^1 +D^\top_1P^1C+S^1_1)\big]\bar{x}^{u_2}\\
                &\qquad +\big[D_2-D_1(R^1_{11}+D^\top_1P^1D_1)^{-1}(R^1_{12}+D^\top_1P^1D_2)\big]u_2\\
                &\qquad -D_1(R^1_{11}+D^\top_1P^1D_1)^{-1}(B^\top_1\eta^{1,u_2}+D^\top_1\zeta^{1,u_2}+D^\top_1P^1\sigma+\rho^1_1)+\sigma \Big\}dW,\\
   d\eta^{1,u_2}&=-\Big\{ \big[A^\top-(P^1B_1+C^\top P^1D_1+S^{1\top}_1)(R^1_{11}+D^\top_1P^1D_1)^{-1}B_1^\top\big]\eta^{1,u_2}\\
                &\qquad+\big[C^\top-(P^1B_1+C^\top P^1D_1+S^{1\top}_1)(R^1_{11}+D^\top_1P^1D_1)^{-1}D_1^\top\big]\zeta^{1,u_2}\\
                &\qquad+\big[C^\top P^1D_2+S^{1\top}_2+P^1B_2-(P^1B_1+C^\top  P^1D_1+S^{1\top}_1)\\
                &\qquad \times(R^1_{11}+D^\top_1P^1D_1)^{-1}(R^1_{12}+D^\top_1P^1 D_2)\big]u_2\\
                &\qquad +\big[C^\top-(P^1B_1+C^\top P^1D_1+S^{1\top}_1)(R^1_{11}+D^\top_1P^1D_1)^{-1}D_1^\top\big]P^1\sigma\\
                &\qquad-(P^1B_1+C^\top P^1D_1+S^{1\top}_1)(R^1_{11}+D^\top_1P^1D_1)^{-1}\rho^1_1+P^1 b+q^1 \Big\}ds+\zeta^{1,u_2}dW,\\
\bar{x}^{u_2}(0)&=x,\,\,\,\,\,\eta^{1,u_2}(T)=g^1.
\end{aligned}
\right.
\end{equation}
In the above, if we denote
\begin{equation*}
\begin{cases}
\hat{R}^1_{11}=R^1_{11}+D^\top_1P^1D_1,\quad\hat{R}^1_{12}=R^1_{12}+D^\top_1P^1D_2,\quad\hat{R}^1_{21}=R^1_{21}+D^\top_2P^1D_1,\\ \hat{\rho}^1_1=\rho_1^1+D^\top_1P^1\sigma,
\quad \hat{S}^1_1=B_1^\top P^1 +D^\top_1P^1C+S^1_1,\quad \hat{S}^1_2=B_2^\top P^1 +D^\top_2P^1C+S^1_2,
\\\hat{A}=A-B_1(\hat{R}^{1}_{11})^{-1}\hat{S}_1^1,\quad \hat{C}=C-D_1(\hat{R}^{1}_{11})^{-1}\hat{S}^1_1,\quad \hat{D}_1=-D_1(\hat{R}^{1}_{11})^{-1}D^\top_1,\\
\hat{F}_1=-B_1(\hat{R}^{1}_{11})^{-1}B_1^\top,\quad \hat{B}_1=-B_1(\hat{R}^{1}_{11})^{-1}D_1^\top,\quad \hat{B}_2=B_2-B_1(\hat{R}^{1}_{11})^{-1}\hat{R}^1_{12},\quad\\
\hat{D}_2=D_2-D_1(\hat{R}^{1}_{11})^{-1}\hat{R}^1_{12},\quad
\hat{b}=b-B_1(\hat{R}^{1}_{11})^{-1}\hat{\rho}^1_1,\quad \hat{\sigma}=\sigma-D_1(\hat{R}^{1}_{11})^{-1}\hat{\rho}^1_1,\quad\\
\hat{F}_2=\hat{S}^1_2-\hat{R}^1_{21}(\hat{R}^{1}_{11})^{-1}\hat{S}^1_1,\quad
\hat{\beta}=\hat{C}^\top P^1\sigma+P^1b+q^1-(\hat{S}^1_1)\top(\hat{R}^{1}_{11})^{-1}\rho_1^1,
\end{cases}
\end{equation*}
then
\begin{equation}\label{leader state}
\left\{
\begin{aligned}
  d\bar{x}^{u_2}&=\big[ \hat{A}\bar{x}^{u_2}+\hat{F}_1\eta^{1,u_2}+\hat{B}_1\zeta^{1,u_2}++\hat{B}_2u_2+\hat{b} \big]ds\\
                &\qquad+\big[\hat{C}\bar{x}^{u_2}+\hat{B}^\top_1\eta^{1,u_2}+\hat{D}_1\zeta^{1,u_2}+\hat{D}_2u_2+\hat{\sigma}\big]dW\\
   d\eta^{1,u_2}&=-\bigl\{\hat{A}^\top\eta^{1,u_2}+\hat{C}^\top\zeta^{1,u_2}+\hat{F}_2^\top u_2+\hat{\beta} \bigr\}ds+\zeta^{1,u_2}dW,\\
\bar{x}^{u_2}(0)&=x,\,\,\,\,\,\eta^{1,u_2}(T)=g^1
\end{aligned}\right.
\end{equation}
Knowing that the follower has chosen a closed-loop optimal strategy
$$(\bar{\Theta}_1(\cdot),\bar{v}_1(\cdot))\equiv(\bar{\Theta}_1[u_2](\cdot),\bar{v}_1[u_2](\cdot))$$
such that its outcome $\bar{u}_1(\cdot)\equiv\bar{u}_1[u_2](\cdot)$ is of the form (\ref{follower optimal control}):
\begin{equation}
\bar{u}_1=-\hat{R}^1_{11}\hat{S}^1_1\bar{x}^{u_2}-\hat{R}^1_{11}B^\top_1\eta^{1,u_2}-\hat{R}^1_{11}D^\top_1\zeta^{1,u_2}-\hat{R}^1_{11}\hat{R}^1_{12}u_2-\hat{R}^1_{11}\hat{\rho}^1_1,
\end{equation}
the leader wish to choose an optimal control $\bar{u}_2(\cdot) \in \mathcal{U}_2[0,T]$ such that his/her cost functional
\begin{equation}
\begin{aligned}
&\hat{J}_2(x;u_2(\cdot)) \triangleq J_2(x;\bar{u}_1(\cdot),u_2(\cdot))\\
&=\mathbb{E} \Bigg\{\int_0^T \bigg[\bigg\langle
\left( \begin{array}{cccc}
\hat{Q}_{11} & \hat{Q}^\top_{12} & \hat{Q}^\top_{13} &\hat{K}^\top_1\\
\hat{Q}_{12} & \hat{Q}_{22}      & \hat{Q}^\top_{23} &\hat{K}^\top_2\\
\hat{Q}_{13} & \hat{Q}_{23}      & \hat{Q}_{33}      &\hat{K}^\top_3\\
\hat{K}_{1}  & \hat{K}_2         & \hat{K}_3         &\hat{R}_2 \end{array} \right)
\left( \begin{array}{c} \bar{x}^{u_2} \\
                        \eta^{1,u_2}  \\
                        \zeta^{1,u_2} \\
                             u_2     \end{array} \right),
\left( \begin{array}{c} \bar{x}^{u_2} \\
                          \eta^{1,u_2}\\
                         \zeta^{1,u_2}\\
                               u_2   \end{array}\right)\bigg\rangle\\
&\qquad +2\bigg\langle
\left( \begin{array}{c} \hat{q}_1 \\
                         \hat{q}_2\\
                         \hat{q}_3\\
                         \hat{\rho} \end{array} \right),
\left( \begin{array}{c} \bar{x}^{u_2}\\
                         \eta^{1,u_2}\\
                        \zeta^{1,u_2}\\
                             u_2 \end{array} \right)\bigg\rangle+\hat{l}\bigg]ds
+\big\langle G^2\bar{x}^{u_2}(T),\bar{x}^{u_2}(T)\big\rangle+2\big\langle g^2 ,\bar{x}^{u_2}(T)\big\rangle \Bigg\}
\end{aligned}
\end{equation}
is minimized, where
\begin{equation*}
\begin{cases}
\hat{Q}_{11}=Q^2-S^{2\top}_1\hat{R}^1_{11}\hat{S}^1_1-\hat{S}^{1\top}_1(\hat{R}^1_{11})^{-1}S^2_1+\hat{S}^{1\top}_1(\hat{R}^1_{11})^{-1}R^2_{11}\hat{R}^1_{11}\hat{S}^1_1,\\
\hat{Q}_{12}=B_1(\hat{R}^1_{11})^{-1}[R^2_{11}(\hat{R}^1_{11})^{-1}\hat{S}^1_1-S^2_1],\quad
\hat{Q}_{13}=D_1(\hat{R}^1_{11})^{-1}[R^2_{11}(\hat{R}^1_{11})^{-1}\hat{S}^1_1-S^2_1],\\
\hat{Q}_{22}=B_1(\hat{R}^1_{11})^{-1}R^2_{11}(\hat{R}^1_{11})^{-1}B^\top_1,\quad
\hat{Q}_{23}=D_1(\hat{R}^1_{11})^{-1}R^2_{11}(\hat{R}^1_{11})^{-1}B^\top_1,
\\
\hat{Q}_{33}=D_1(\hat{R}^1_{11})^{-1}R^2_{11}(\hat{R}^1_{11})^{-1}D^\top_1,\\
\hat{K}_1=S^2_2-R^2_{21}(\hat{R}^1_{11})^{-1}\hat{S}^1_1-\hat{R}^{1\top}_{12}(\hat{R}^1_{11})^{-1}S^2_1+\hat{R}^{1\top}_{12}\hat{R}^1_{11}R^2_{11}(\hat{R}^1_{11})^{-1}\hat{S}^1_1,\\
\hat{R}_2=R^2_{22}+\hat{R}^{1\top}_{12}(\hat{R}^1_{11})^{-1}R^2_{11}(\hat{R}^1_{11})^{-1}\hat{R}^1_{12}-\hat{R}^{1\top}_{12}(\hat{R}^1_{11})^{-1}R^2_{12}-R^2_{21}(\hat{R}^1_{11})^{-1}\hat{R}^1_{12},\\
\hat{K}_2=[\hat{R}^{1\top}_{12}(\hat{R}^1_{11})^{-1}R^2_{11}-R^2_{21}](\hat{R}^1_{11})^{-1}B^\top_1,\quad \hat{K}_3=[\hat{R}^{1\top}_{12}(\hat{R}^1_{11})^{-1}R^2_{11}-R^2_{21}](\hat{R}^1_{11})^{-1}D^\top_1,\\
\hat{q}_1=q^2+\hat{S}^{1\top}_1(\hat{R}^1_{11})^{-1}\rho^2_1+[S^{2\top}_1-\hat{S}^{1\top}_1(\hat{R}^1_{11})^{-1}R^2_{11}](\hat{R}^1_{11})^{-1}\hat{\rho}^1_1,\\
\hat{q}_2=B_1(\hat{R}^1_{11})^{-1}[R^2_{11}(\hat{R}^1_{11})^{-1}\hat{\rho}^1_1-\rho^2_1],\quad
\hat{q}_3=D_1(\hat{R}^1_{11})^{-1}[R^2_{11}(\hat{R}^1_{11})^{-1}\hat{\rho}^1_1-\rho^2_1],\\
\hat{\rho}=\rho^2_2+[\hat{R}^{1\top}_{12}(\hat{R}^1_{11})^{-1}R^2_{11}-R^2_{21}](\hat{R}^1_{11})^{-1}\hat{\rho}^1_1-\hat{R}^{1\top}_{12}(\hat{R}^1_{11})^{-1}\rho^2_1,\\
\hat{l}=\langle R^2_{11}(\hat{R}^1_{11})^{-1}\hat{\rho}^1_1-2\rho^2_1,(\hat{R}^1_{11})^{-1}\hat{\rho}^1_1  \rangle.
\end{cases}
\end{equation*}
The LQ problem of the leader can be stated as follows.

\vspace{1mm}

\textbf{Problem (SLQ)$_l$}. For given $x \in \mathbb{R}^n$, find a $\bar{u}_2(\cdot) \in \mathcal{U}_2[0,T]$ such that
\begin{equation}\label{LLQ}
\hat{J}_2(x;\bar{u}_2(\cdot))=\underset{u_2(\cdot)\in \mathcal{U}_2[0,T]} {\min}\hat{J}_2(x;u_2(\cdot))\equiv V_2(x).
\end{equation}
The above Problem (SLQ)$_l$ is an LQ problem of FBSDE. Any $\bar{u}_2(\cdot) \in \mathcal{U}_2[0,T]$ satisfying (\ref{LLQ}) is called \textit{an open-loop optimal control} of Problem (SLQ)$_l$ for $x$, the corresponding $(\bar{x}(\cdot),\bar{\eta}^1(\cdot),\bar{\zeta}^1(\cdot))\\\equiv (\bar{x}^{\bar{u}_2}(\cdot),\bar{\eta}^{1,\bar{u}_2}(\cdot),\bar{\zeta}^{1,\bar{u}_2}(\cdot))$ is called \textit{an open-loop optimal state process triple} and $(\bar{x}(\cdot),\bar{\eta}^1(\cdot),\bar{\zeta}^1(\cdot),\\\bar{u}_2(\cdot))$ is called \textit{an open-loop optimal quadruple}. The map $V_2(\cdot)$ is called \textit{the value function} of Problem (SLQ)$_l$. The (unique) {\it open-loop solvability} of Problem (SLQ)$_l$ can be similarly defined as Definition \ref{def2.1}.

%leader open-loop optimal
\begin{mythm}\label{lsc}
Let (H1)-(H2) hold. For a given $x \in \mathbb{R}^n$, a quadruple $(\bar{x}(\cdot),\bar{\eta}^1(\cdot),\bar{\zeta}^1(\cdot),\bar{u}_2(\cdot))$ is an open-loop optimal quadruple of Problem (SLQ)$_l$ if and only if the following stationarity condition holds:
\begin{equation}\label{leader open stationarity}
	\hat{F}_2p^{2,\bar{u}_2}+\hat{B}^\top_2q^{2,\bar{u}_2}+\hat{D}_2^\top k^{2,\bar{u}_2}+\hat{K}_1\bar{x}+\hat{K}_2\bar{\eta}^1+\hat{K}_3\bar{\zeta}^1+\hat{R}_2\bar{u}_2+\hat{\rho}=0,\quad a.e.,\, \mathbb{P}\mbox{-}a.s.,
	\end{equation}
where $(p^{2,\bar{u}_2}(\cdot),q^{2,\bar{u}_2}(\cdot),k^{2,\bar{u}_2}(\cdot))\in L^2_{\mathbb{F}}(0,T;\mathbb{R}^n)\times L^2_{\mathbb{F}}(0,T;\mathbb{R}^n)\times L^2_{\mathbb{F}}(0,T;\mathbb{R}^n)$ is the adapted solution to the following FBSDE:
\begin{equation}\label{leader adjoint equation}
\left\{
\begin{aligned}
dp^{2,\bar{u}_2}&=\big[\hat{A}p^{2,\bar{u}_2}+\hat{F}^\top_1 q^{2,\bar{u}_2}+\hat{B}_1k^{2,\bar{u}_2}+\hat{Q}_{12}\bar{x}+\hat{Q}_{22}\bar{\eta}^1+\hat{Q}^\top_{23}\bar{\zeta}^1+\hat{K}^\top_2\bar{u}_2+\hat{q}_2\big]ds\\
	            &\qquad +\big[\hat{C}p^{2,\bar{u}_2}+\hat{B}^\top_1 q^{2,\bar{u}_2}+\hat{D}^\top_1k^{2,\bar{u}_2}+\hat{Q}_{13}\bar{x}+\hat{Q}_{23}\bar{\eta}^1+\hat{Q}_{33}\bar{\zeta}^1+\hat{K}^\top_3\bar{u}_2+\hat{q}_3\big]dW,\\
dq^{2,\bar{u}_2}&=-\big[\hat{A}^\top q^{2,\bar{u}_2}+\hat{C}^\top k^{2,\bar{u}_2}+\hat{Q}_{11}\bar{x}+\hat{Q}^\top_{12}\bar{\eta}^1+\hat{Q}^\top_{13}\bar{\zeta}^1+\hat{K}_1^\top\bar{u}_2+\hat{q}_1\big]ds+k^{2,\bar{u}_2}dW,\\
	p^{2,\bar{u}_2}(0)&=0,\,\,\,q^{2,\bar{u}_2}(T)=G^2\bar{x}(T)+g^2,
\end{aligned}
\right.
\end{equation}
and the following convexity condition holds:
\begin{equation}\label{leader convexity condition}
\begin{aligned}
	&\mathbb{E} \bigg\{\int_0^T\Big[ \big\langle \hat{Q}_{11}x_{0l},x_{0l} \big\rangle+2\big\langle \hat{Q}_{12}x_{0l},\eta^0\big\rangle+2\big\langle \hat{Q}_{13}x_{0l},\zeta^0\big\rangle+2\big\langle \hat{Q}_{23}\eta^0,\zeta^0\big\rangle+\big\langle \hat{Q}_{22}\eta^0,\eta^0\big\rangle\\
	&\qquad +\big\langle \hat{Q}_{22}\zeta^0,\zeta^0\big\rangle+2\big\langle \hat{K}_1x_{0l},u_2\big\rangle+2\big\langle \hat{K}_2\eta^0,u_2\big\rangle+2\big\langle \hat{K}_3\zeta^0,u_2\big\rangle+\big\langle \hat{R}_2u_2,u_2\big\rangle \Big]ds\\
	&\qquad\qquad\qquad\qquad +\big\langle G^2x_{0l}(T),x_{0l}(T)\big\rangle\bigg\} \geqslant 0,\qquad \forall u_2 \in \mathcal{U}_2[0,T],
\end{aligned}
\end{equation}
where $(x_{0l}(\cdot),\eta^0(\cdot),\zeta^0(\cdot))\in L^2_{\mathbb{F}}(0,T;\mathbb{R}^n)\times L^2_{\mathbb{F}}(0,T;\mathbb{R}^n)\times L^2_{\mathbb{F}}(0,T;\mathbb{R}^n)$ is the solution to the following:
\begin{equation}\label{leader x_0}
\left\{
\begin{aligned}
  dx_{0l}&=\big[\hat{A}x_{0l}+\hat{F}_1\eta^0+\hat{B}_1\zeta^0+\hat{B}_2u_2 \big]ds+\big[\hat{C}x_{0l}+\hat{B}_1^\top\eta^0+\hat{D}_1\zeta^0+\hat{D}_2u_2 \big]dW,\\
  d\eta^0&=-\big[\hat{A}^\top\eta^0+\hat{C}^\top\zeta^0+\hat{F}_2^\top u_2 \big]ds+\zeta^0dW,\\
x_{0l}(0)&=0,\,\,\,\,\,\eta^0(T)=0.
\end{aligned}
\right.
\end{equation}
\end{mythm}

\textit{Proof}. Suppose $(\bar{x}(\cdot),\bar{\eta}^{1}(\cdot),\bar{\zeta}^{1}(\cdot),\bar{u}_2(\cdot))$ is a quadruple corresponding to the given $x \in\mathbb{R}^n$. For any $u_2 \in \mathcal{U}_2[0,T]$ and $\epsilon \in \mathbb{R}$, let $u_2^\epsilon(\cdot)=\bar{u}_2(\cdot)+\epsilon u_2(\cdot)$ and $(\bar{x}^{\epsilon}(\cdot)\equiv x(\cdot;\bar{\Theta}_1,\bar{v}_1,\bar{u}_2+\epsilon u_2(\cdot)),\eta^{1,\epsilon}(\cdot),\zeta^{1,\epsilon}(\cdot))$ be the corresponding state. Then $(\bar{x}^\epsilon(\cdot),\eta^{1,\epsilon}(\cdot),\zeta^{1,\epsilon}(\cdot))$ satisfies
\begin{equation}\nonumber
\left\{
\begin{aligned}
 d\bar{x}^\epsilon&=\big[ \hat{A}\bar{x}^\epsilon+\hat{F}_1\eta^{1,\epsilon}+\hat{B}_1\zeta^{1,\epsilon}+\hat{B}_2(\bar{u}_2+\epsilon u_2)+\hat{b}\big]ds\\
                  &\qquad +\big[ \hat{C}\bar{x}^\epsilon+\hat{B}^\top_1\eta^{1,\epsilon}+\hat{D}_1\zeta^{1,\epsilon}+\hat{D}_2(\bar{u}_2+\epsilon u_2)+\hat{\sigma} \big] dW,\\
d\eta^{1,\epsilon}&=-\big[ \hat{A}^\top\eta^{1,\epsilon}+\hat{C}^\top\zeta^{1,\epsilon}+\hat{F}_2^\top(\bar{u}_2+\epsilon u_2)+\hat{\beta} \big]ds+\zeta^{1,\epsilon} dW,\\
\bar{x}^\epsilon(0)&=x,\,\,\,\,\,\eta^{1,\epsilon}(T)=g^1.
\end{aligned}
\right.
\end{equation}
Thus, $x_{0l}(\cdot) \equiv \frac{\bar{x}^{\epsilon}(\cdot)-\bar{x}(\cdot)}{\epsilon}$ is independent of $\epsilon$ and satisfies (\ref{leader x_0}). Then we get
\begin{equation}\nonumber
\begin{aligned}
&\hat{J}_2(x;\bar{u}_2(\cdot)+\epsilon u_2(\cdot))-\hat{J}_2(x;\bar{u}_2(\cdot))\\
&=2\epsilon\mathbb{E}\biggl\{ \int_0^T\Big[ \big\langle \hat{Q}_{11}\bar{x}^{\bar{u}_2},x_{0l} \big\rangle +\big\langle \hat{Q}_{12}\bar{x}^{\bar{u}_2},\eta^0 \big\rangle
+\big\langle \hat{Q}^\top_{12}\bar{\eta}^1,x_{0l} \big\rangle +\big\langle \hat{Q}_{13}\bar{x}^{\bar{u}_2},\zeta^0 \big\rangle +\big\langle \hat{Q}^\top_{13}\bar{\zeta}^1,x_{0l} \big\rangle\\
&\qquad +\big\langle \hat{Q}_{22}\bar{\eta}^1,\eta^0 \big\rangle +\big\langle \hat{Q}_{23}\bar{\eta}^1,\zeta^0 \big\rangle +\big\langle \hat{Q}^\top_{23}\bar{\zeta}^1,\eta^0\big\rangle +\big\langle \hat{Q}_{33}\bar{\zeta}^1,\zeta^0 \big\rangle +\big\langle \hat{K}_1\bar{x},u_2 \big\rangle +\big\langle \hat{K}^\top_1\bar{u}_2,x_{0l} \big\rangle\\
&\qquad +\big\langle \hat{K}_2\bar{\eta},u_2 \big\rangle +\big\langle \hat{K}^\top_2\bar{u}_2,\eta^{0} \big\rangle +\big\langle \hat{K}_3\bar{\zeta},u_2 \big\rangle +\big\langle \hat{K}^\top_3\bar{u}_2,\zeta^{0} \big\rangle
+\big\langle \hat{R}_2\bar{u_2},u_2 \big\rangle +\big\langle \hat{q}_1,x_{0l}\big\rangle\\
&\qquad +\big\langle\hat{q}_2,\eta^0\big\rangle+\big\langle \hat{q}_3,\zeta^0 \big\rangle+\big\langle \hat{\rho},u_2 \big\rangle \Big]ds +\big\langle G^2\bar{x}(T),x_{0l}(T)\big\rangle
+\big\langle g^2,x_{0l}(T)\big\rangle\biggr\}\\
&\quad +\epsilon^2\mathbb{E}\biggl\{ \int_0^T\Big[ \big\langle \hat{Q}_{11}x_{0l},x_{0l} \big\rangle+2\big\langle \hat{Q}_{12}x_{0l},\eta^0\big\rangle+2\big\langle \hat{Q}_{13}x_{0l},\zeta^0\big\rangle+\big\langle \hat{Q}_{22}\eta^0,\eta^0\big\rangle+2\big\langle \hat{Q}_{23}\eta^0,\zeta^0\big\rangle\\
&\qquad +\big\langle \hat{Q}_{33}\zeta^0,\zeta^0\big\rangle +2\big\langle \hat{K}_1x_{0l},u_2\big\rangle +2\big\langle \hat{K}_2\eta^0,u_2\big\rangle +2\big\langle \hat{K}_3\zeta^0,u_2\big\rangle
+\big\langle \hat{R}_2u_2,u_2\big\rangle \Big]ds\\
&\qquad +\big\langle G^2x_{0l}(T),x_{0l}(T)\rangle \biggr\}.
\end{aligned}
\end{equation}
Applying It\^o's formula to $\big\langle q^{2,\bar{u}_2}(\cdot),x_{0l}(\cdot) \big\rangle$ and $-\big\langle p^{2,\bar{u}_2}(\cdot),\eta^0(\cdot) \big\rangle$, we get
\begin{equation*}
\begin{aligned}
&\hat{J}_2(x;\bar{u}_2(\cdot)+\epsilon u_2(\cdot))-\hat{J}_2(x;\bar{u}_2(\cdot))\\
&=2\epsilon\mathbb{E}\biggl\{  \int_0^T \big\langle \hat{F}_2p^{2,\bar{u}_2}+\hat{B}^\top_2q^{2,\bar{u}_2}+\hat{D}_2^\top k^{2,\bar{u}_2}
+\hat{K}_1\bar{x}+\hat{K}_2\bar{\eta}^1+\hat{K}_3\bar{\zeta}^1+\hat{R}_2\bar{u}_2+\hat{\rho},u_2 \big\rangle ds \biggr\}\\
&\quad +\epsilon^2\mathbb{E}\biggl\{ \int_0^T\Big[ \big\langle \hat{Q}_{11}x_{0l},x_{0l} \big\rangle +2\big\langle \hat{Q}_{12}x_{0l},\eta^0\big\rangle +2\big\langle \hat{Q}_{13}x_{0l},\zeta^0\big\rangle
+\big\langle \hat{Q}_{22}\eta^0,\eta^0\big\rangle +2\big\langle \hat{Q}_{23}\eta^0,\zeta^0\big\rangle\\
&\qquad\qquad  +\big\langle \hat{Q}_{33}\zeta^0,\zeta^0\big\rangle +2\big\langle \hat{K}_1x_{0l},u_2\big\rangle +2\big\langle \hat{K}_2\eta^0,u_2\big\rangle +2\big\langle \hat{K}_3\zeta^0,u_2\big\rangle
+\big\langle \hat{R}_2u_2,u_2\big\rangle \Big]ds\\
&\qquad\qquad +\big\langle G^2x_{0l}(T),x_{0l}(T)\rangle \biggr\}.
\end{aligned}
\end{equation*}
Therefore, $(\bar{x}(\cdot),\bar{\eta}^1(\cdot),\bar{\zeta}^1(\cdot),\bar{u}_2(\cdot))$ is an open-loop optimal quadruple of Problem (SLQ)$_{l}$ if and only if (\ref{leader open stationarity}) and (\ref{leader convexity condition}) hold. The proof is complete. $\qquad\Box$

As in Definition \ref{def2.3}, next, we take $\Theta_2(\cdot) \in \mathcal{Q}_2[0,T],\,\,\check{\Theta}_2(\cdot) \in \mathcal{Q}_2[0,T]$ and $v_2(\cdot)\in \mathcal{U}_2[0,T]$. For any $x \in \mathbb{R}^n$, let us consider the following FBSDE on $[0,T]$:
\begin{equation}\label{leader closedloop system}
\left\{
\begin{aligned}
d\bar{x}^{\Theta_2,\check{\Theta}_2,v_2}&=\big[ (\hat{A}+\hat{B}_2\Theta_2)\bar{x}^{\Theta_2,\check{\Theta}_2,v_2}+(\hat{F}_1+\hat{B}_2\check{\Theta}_2)\eta^{1,\Theta_2,\check{\Theta}_2,v_2}\\
                                        &\qquad +\hat{B}_1\zeta^{1,\Theta_2,\check{\Theta}_2,v_2}+\hat{B}_2v_2+\hat{b} \big]ds\\
                                        &\quad +\big[ (\hat{C}+\hat{D}_2\Theta_2)\bar{x}^{\Theta_2,\check{\Theta}_2,v_2}+(\hat{B}^\top_1+\hat{D}_2\check{\Theta}_2)\eta^{1,\Theta_2,\check{\Theta}_2,v_2}\\
                                        &\qquad +\hat{D}_1\zeta^{1,\Theta_2,\check{\Theta}_2,v_2}+\hat{D}_2v_2+\hat{\sigma} \big]dW,\\
 d\eta^{1,\Theta_2,\check{\Theta}_2,v_2}&=-\big[(\hat{A}+\check{\Theta}_2^\top\hat{F}_2)^\top\eta^{1,\Theta_2,\check{\Theta}_2,v_2}+C^\top\zeta^{1,\Theta_2,\check{\Theta}_2,v_2}
                                         +\hat{F}^\top_2\Theta_2\bar{x}^{\Theta_2,\check{\Theta}_2,v_2}\\
                                        &\qquad +\hat{F}_2^\top v_2+\hat{\beta}\big]ds+\zeta^{1,\Theta_2,\check{\Theta}_2,v_2}dW,\\
\bar{x}^{\Theta_2,\check{\Theta}_2,v_2}(0)&=x,\,\,\,\,\,\eta^{1,\Theta_2,\check{\Theta}_2,v_2}(T)=g^1.
\end{aligned}
\right.
\end{equation}
This is a fully coupled FBSDE which admits a unique solution $(\bar{x}^{\Theta_2,\check{\Theta}_2,v_2}(\cdot),\eta^{1,\Theta_2,\check{\Theta}_2,v_2}(\cdot),\\\zeta^{1,\Theta_2,\check{\Theta}_2,v_2}(\cdot))$, depending on $\Theta_2(\cdot),\,\check{\Theta}_2(\cdot)$ and $v_2(\cdot)$. (\ref{leader closedloop system}) is called the {\it closed-loop system} of the original state equation (\ref{leader state}) under the closed-loop strategy $(\Theta_2(\cdot),\check{\Theta}_2(\cdot),v_2(\cdot))$ of the leader. Similarly, we point out that $(\Theta_2(\cdot),\check{\Theta}_2(\cdot),v_2(\cdot))$ is independent of the initial state $x$. For the sake of simplicity, we denote the above $(\bar{x}^{\Theta_2,\check{\Theta}_2,v_2}(\cdot),\eta^{1,\Theta_2,\check{\Theta}_2,v_2}(\cdot),\zeta^{1,\Theta_2,\check{\Theta}_2,v_2}(\cdot))$ as $(\bar{\check{x}}(\cdot),\check{\eta}^1(\cdot),\check{\zeta}^1(\cdot))$, depending on $\Theta_2(\cdot),\,\check{\Theta}_2(\cdot),\,v_2(\cdot)$, and define
\begin{equation}\label{lclf}
\begin{aligned}
&\hat{J}_2(x;\Theta_2\bar{\check{x}}+\check{\Theta}_2\check{\eta}^1+v_2)\\
&=\mathbb{E} \Bigg\{
\int_0^T \bigg[\bigg\langle
\left( \begin{array}{cccc}
\hat{Q}_{11} & \hat{Q}^\top_{12} & \hat{Q}^\top_{13} &\hat{K}^\top_1\\
\hat{Q}_{12} & \hat{Q}_{22}      & \hat{Q}^\top_{23} &\hat{K}^\top_2\\
\hat{Q}_{13} & \hat{Q}_{23}      & \hat{Q}_{33}      &\hat{K}^\top_3\\
\hat{K}_{1}  & \hat{K}_2         & \hat{K}_3         &\hat{R}_2 \end{array} \right)
\left( \begin{array}{c} \bar{\check{x}} \\
\check{\eta}^1  \\
\check{\zeta}^1 \\
\Theta_2\bar{\check{x}}+\check{\Theta}_2\check{\eta}^1+v_2     \end{array} \right),
\left( \begin{array}{c} \bar{\check{x}} \\
\check{\eta}^1\\
\check{\zeta}^1\\
\Theta_2\bar{\check{x}}+\check{\Theta}_2\check{\eta}^1+v_2   \end{array}\right)\bigg\rangle\\
&\qquad +2\bigg\langle
\left( \begin{array}{c} \hat{q}_1 \\
\hat{q}_2\\
\hat{q}_3\\
\hat{\rho} \end{array} \right),
\left( \begin{array}{c} \bar{\check{x}}\\
\check{\eta}^1\\
\check{\zeta}^1\\
\Theta_2\bar{\check{x}}+\check{\Theta}_2\check{\eta}^1+v_2 \end{array} \right)\bigg\rangle+\hat{l}\bigg]ds +\big\langle G^2\bar{\check{x}}(T),\bar{\check{x}}(T)\big\rangle
+2\big\langle g^2,\bar{\check{x}}(T)\big\rangle \Bigg\}.
\end{aligned}
\end{equation}

\begin{mydef}\label{def4.1}
A triple $(\bar{\Theta}_2(\cdot),\bar{\check{\Theta}}_2,\bar{v}_2(\cdot)) \in \mathcal{Q}_2[0,T] \times \mathcal{Q}_2[0,T] \times \mathcal{U}_2[0,T]$ is called a \textit{closed-loop optimal strategy} of Problem (SLQ)$_l$ on $[0,T]$ if
\begin{equation}
\begin{split}
	&\hat{J}_2(x;\bar{\Theta}_2(\cdot)\bar{x}(\cdot)+\bar{\check{\Theta}}_2(\cdot)\bar{\eta}^1(\cdot)+\bar{v}_2(\cdot))\\
    &\leqslant \hat{J}_2(x;\Theta_2(\cdot)\bar{x}^{\Theta_2,\check{\Theta}_2,v_2(\cdot)}+\check{\Theta}_2(\cdot)\eta^{1,\Theta_2,\check{\Theta}_2,v_2(\cdot)}(\cdot)+v_2(\cdot)),\\
	&\qquad \forall x \in \mathbb{R}^n,\,\,\forall  (\Theta_2(\cdot),\check{\Theta}_2,v_2(\cdot)) \in \mathcal{Q}_2[0,T] \times \mathcal{Q}_2[0,T] \times \mathcal{U}_2[0,T],
\end{split}
\end{equation}
where $\bar{x}(\cdot)\equiv \bar{x}^{\bar{\Theta}_2,\bar{\tilde{\Theta}}_2,\bar{v}_2}(\cdot)$ together with $\bar{\eta}^1(\cdot)$ and $\bar{\zeta}^1(\cdot)$ satisfying (\ref{leader closedloop system}).
\end{mydef}

The following result is similar to Proposition 3.3 of \cite{SY2014}, and the detailed proof is omitted.
\begin{mypro}\label{pro4.1}
Let (H1)-(H2) hold. Then the following are equivalent:
\par (\romannumeral 1) $(\bar{\Theta}_2(\cdot),\bar{\check{\Theta}}_2(\cdot),\bar{v}_2(\cdot)) \in \mathcal{Q}_2[0,T] \times \mathcal{Q}_2[0,T] \times \mathcal{U}_2[0,T]$ is a closed-loop optimal strategy of Problem (SLQ)$_l$.
\par (\romannumeral 2) The following holds:
\begin{equation*}
\begin{aligned}
	\hat{J}_2(x;\bar{\Theta}_2(\cdot)\bar{x}(\cdot)+\bar{\check{\Theta}}_2(\cdot)\bar{\eta}^1(\cdot)+\bar{v}_2(\cdot)) \leqslant \hat{J}_2(x;\bar{\Theta}_2(\cdot)\bar{x}^{\bar{\Theta}_2,\bar{\check{\Theta}}_2,v_2}(\cdot)&+\bar{\check{\Theta}}_2(\cdot)\eta^{1,\bar{\Theta}_2,\bar{\check{\Theta}},v_2}(\cdot)+v_2(\cdot))\\
	                  &\forall x \in \mathbb{R}^n,\,\, \forall v_2(\cdot) \in \mathcal{U}_2[0,T].
\end{aligned}
\end{equation*}

\par (\romannumeral 3) The following holds:
\begin{equation}\label{open and close equivalent}
\begin{aligned}
	\hat{J}_2(x;\bar{\Theta}_2(\cdot)\bar{x}(\cdot)+\bar{\check{\Theta}}_2(\cdot)\bar\eta^1(\cdot)+\bar{v}_2(\cdot)) \leqslant \hat{J}_2(x;u_2(\cdot)),\quad
	\forall x \in \mathbb{R}^n,\,\, \forall u_2(\cdot) \in \mathcal{U}_2[0,T].
\end{aligned}
\end{equation}
\end{mypro}
From (\ref{open and close equivalent}), we can see that for a fixed $x \in \mathbb{R}^n$, the {\it outcome}
\begin{equation}
\bar{u}_2(\cdot) \equiv \bar{\Theta}_2(\cdot)\bar{x}(\cdot)+\bar{\check{\Theta}}_2(\cdot)\bar{\eta}^1(\cdot)+\bar{v}_2(\cdot) \in \mathcal{U}_2[0,T]
\end{equation}
of the closed-loop optimal strategy $(\bar{\Theta}_2(\cdot),\bar{\check{\Theta}}_2(\cdot),\bar{v}_2(\cdot))$ is an open-loop optimal control of Problem (SLQ)$_l$, but the existence of an open-loop optimal control cannot guarantee the existence of a closed-loop optimal strategy. At the end of this paper, we will give some examples to prove that open-loop solvability is weaker than the closed-loop solvability. But if open-loop optimal control and closed-loop optimal strategy both exist, the feedback representation of the open-loop optimal control and the outcome of the closed-loop optimal strategy are consistent.

On the other hand, we can also see that if $(\bar{\Theta}_2(\cdot),\bar{\check{\Theta}}_2(\cdot),\bar{v}_2(\cdot))$ is a closed-loop optimal strategy of Problem (SLQ)$_l$, then $\bar{v}_2(\cdot)$ is an open-loop optimal control of the LQ problem (\ref{leader closedloop system})-(\ref{lclf}), with $\Theta_2(\cdot)=\bar{\Theta}_2(\cdot),\,\,\check{\Theta}_2(\cdot)=\bar{\check{\Theta}}_2(\cdot)$, which we denote it by {\bf Problem (SLQ)$_{ll}$}.

Similar to the above conclusion, we can give the necessary and sufficient conditions for the open-loop solvability of Problem (SLQ)$_{ll}$.

\begin{mypro}
Let (H1)-(H2) hold. For a given $x \in \mathbb{R}^n$, $(\bar{x}^{\bar{v}_2}(\cdot),\bar{\eta}^{1,\bar{v}_2}(\cdot),\bar{\zeta}^{1,\bar{v}_2}(\cdot),\bar{v}_2(\cdot))$ is an open-loop optimal quadruple of Problem (SLQ)$_{ll}$ if and only if the following stationarity condition holds:
\begin{equation}\label{leader open stationarity---}
\begin{aligned}
&\hat{F}_2p^{2,\bar{v}_2}+\hat{B}^\top_2q^{2,\bar{v}_2}+\hat{D}^\top_2k^{2,\bar{v}_2}+(\hat{K}_1+\hat{R}_2\bar{\Theta}_2)\bar{x}^{\bar{v}_2}+(\hat{K}_2+\hat{R}_2\bar{\check{\Theta}}_2)\bar{\eta}^1\\
&\quad +\hat{K}_3\bar{\zeta}^1+\hat{R}_2\bar{v}_2+\hat{\rho}=0,\quad a.e.,\, \mathbb{P}\mbox{-}a.s.,
\end{aligned}
\end{equation}
where $\bar{x}^{\bar{v}_2}\equiv x(\cdot;\bar{\Theta}_1,\bar{v}_1,\bar{\Theta}_2,\bar{\check{\Theta}}_2,\bar{v}_2)$ and $(p^{2,\bar{v}_2}(\cdot),q^{2,\bar{v}_2}(\cdot),k^{2,\bar{v}_2}(\cdot))\in L^2_{\mathbb{F}}(0,T;\mathbb{R}^n)\times L^2_{\mathbb{F}}(0,T;\mathbb{R}^n)\times L^2_{\mathbb{F}}(0,T;\mathbb{R}^n)$ is the adapted solution to the following FBSDE:
\begin{equation}\label{leader adjoint equation---}
\left\{
\begin{aligned}
   dp^{2,\bar{v}_2}&=\bigl[ (\hat{A}+\bar{\check{\Theta}}^\top_2\hat{F}_2)p^{2,\bar{v}_2}+(\hat{F}_1+\hat{B}_2\bar{\check{\Theta}}_2)^\top q^{2,\bar{v}_2}+(\hat{B}_1+\bar{\check{\Theta}}^\top_2\hat{D}^\top_2)k^{2,\bar{v}_2}\\
                   &\qquad +(\hat{Q}_{12}+\bar{\check{\Theta}}^\top_2\hat{K}_1+\hat{K}^\top_2\bar{\Theta}_2+\bar{\check{\Theta}}^\top_2\hat{R}_2\bar{\Theta}_2)\bar{x}^{\bar{v}_2}+(\hat{Q}_{23}
                    +\hat{K}^\top_3\bar{\check{\Theta}}_2)^\top\bar{\zeta}^{1,\bar{v}_2}\\
                   &\qquad +(\hat{Q}_{22}+\hat{K}^\top_2\bar{\check{\Theta}}_2+\bar{\check{\Theta}}^\top_2\hat{K}_2+\bar{\check{\Theta}}^\top_2\hat{R}_2\bar{\check{\Theta}}_2)\bar{\eta}^{1,\bar{v}_2}
                    +(\hat{K}^\top_2+\bar{\check{\Theta}}^\top_2\hat{R}_2)\bar{v}_2+\hat{q}_2+\bar{\check{\Theta}}^\top_2\hat{\rho} \bigr]ds\\
                   &\quad +\bigl[ \hat{C}p^{2,\bar{v}_2}+\hat{B}^\top_1 q^{2,\bar{v}_2}+\hat{D}^\top_1k^{2,\bar{v}_2}+(\hat{Q}_{13}+\hat{K}^\top_3\bar{\Theta}_2)\bar{x}^{\bar{v}_2}\\
                   &\qquad +(\hat{Q}_{23}+\hat{K}^\top_3\bar{\check{\Theta}}_2)\bar{\eta}^{1,\bar{v}_2}+\hat{Q}_{33}\bar{\zeta}^1+\hat{K}^\top_3\bar{v}_2+\hat{q}_3 \big]dW,\\
   dq^{2,\bar{v}_2}&=-\bigl[\bar{\Theta}^\top_2\hat{F}_2 p^{2,\bar{v}_2} + (\hat{A}+\hat{B}_2\bar{\Theta}_2)^\top q^{2,\bar{v}_2}+(\hat{C}+\hat{D}_2\bar{\Theta}_2)^\top k^{2,\bar{v}_2}\\
                   &\qquad +(\hat{Q}_{11}+\bar{\Theta}^\top_2\hat{K}_1+\hat{K}^\top_1\bar{\Theta}_2+\bar{\Theta}^\top_2\hat{R}_2\bar{\Theta}_2)\bar{x}^{\bar{v}_2}\\
                   &\qquad +(\hat{Q}^\top_{12}+\hat{K}^\top_1\bar{\check{\Theta}}_2+\bar{\Theta}_2^\top\hat{K}_2+\bar{\Theta}^\top_2\hat{R}_2\bar{\check{\Theta}}_2)\bar{\eta}^{1,\bar{v}_2}\\
                   &\qquad +(\hat{Q}^\top_{13}+\bar{\Theta}^\top_2\hat{K}_3)\bar{\zeta}^{1,\bar{v}_2}+(\hat{K}^\top_1+\bar{\Theta}^\top_2\hat{R}_2)\bar{v}_2 +\hat{q}_1+\bar{\Theta}^\top_2\hat{\rho}\bigr]ds+k^{2,\bar{v}_2}dW,\\
 p^{2,\bar{v}_2}(0)&=0,\,\,\,q^{2,\bar{v}_2}(T)=G^2\bar{x}^{\bar{v}_2}(T)+g^2,
\end{aligned}\right.
\end{equation}
and the following convexity condition holds:
\begin{equation}\label{leader convexity condition---}
\begin{aligned}
&\mathbb{E} \bigg\{ \int_0^T\Big[\big\langle\big[\hat{Q}_{11}+\bar{\Theta}^\top_2\hat{K}_1+\hat{K}^\top_1\bar{\Theta}_2+\bar{\Theta}_2^\top \hat{R}_2\bar{\Theta}_2\big]x^{v_2}_{0l},x^{v_2}_{0l}\big\rangle\\
&\quad +2\big\langle (\hat{Q}_{12}+\bar{\check{\Theta}}^\top_2\hat{K}_1+\hat{K}^\top_2\bar{\Theta}_2+\bar{\check{\Theta}}^\top_2\hat{R}_2\bar{\Theta}_2)x^{v_2}_{0l},\eta^{0,v_2} \big\rangle
+2\big\langle (\hat{Q}_{13}+\hat{K}^\top_3\bar{\Theta}_2)x^{v_2}_{0l},\zeta^{0,v_2} \big\rangle\\
&\quad +2\big\langle (\hat{Q}_{23}+\hat{K}^\top_3\bar{\check{\Theta}}_2)\eta^{0,v_2},\zeta^{0,v_2} \big\rangle+\big\langle (\hat{Q}_{22}
+\hat{K}^\top_2\bar{\check{\Theta}}_2+\bar{\check{\Theta}}^\top_2\hat{K}_2+\bar{\check{\Theta}}_2^\top \hat{R}_2 \bar{\check{\Theta}}_2)\eta^{0,v_2},\eta^{0,v_2}\big\rangle\\
&\quad +\big\langle \hat{Q}_{33}\zeta^{0,v_2},\zeta^{0,v_2} \big\rangle+2\big\langle (\hat{K}_1+\hat{R}_2\bar{\Theta}_2)x^{v_2}_{0l},v_2 \big\rangle+2\big\langle (\hat{K}_2+\hat{R}_2\bar{\check{\Theta}}_2)\eta^{0,v_2},v_2 \big\rangle\\
&\quad +2\big\langle \hat{K}_3\zeta^{0,v_2},v_2 \big\rangle+\big\langle \hat{R}_2v_2,v_2\big\rangle \Big]ds +\big\langle G^2x^{v_2}_{0l}(T),x^{v_2}_{0l}(T)\big\rangle \bigg\} \geqslant 0,
\end{aligned}
\end{equation}
where $(x_{0l}^{v_2}(\cdot),\eta^{0,v_2}(\cdot),\zeta^{0,v_2}(\cdot))\in L^2_{\mathbb{F}}(0,T;\mathbb{R}^n)\times L^2_{\mathbb{F}}(0,T;\mathbb{R}^n)\times L^2_{\mathbb{F}}(0,T;\mathbb{R}^n)$ is the solution to the following:
\begin{equation}\label{leader x_0---}
\left\{ \begin{aligned}
   dx_{0l}^{v_2}&=\big[ (\hat{A}+\hat{B}_2\bar{\Theta}_2)x_{0l}^{v_2}+(\hat{F}_1+\hat{B}_2\bar{\check{\Theta}}_2)\eta^{0,v_2}+\hat{B}_1\zeta^{0,v_2}+\hat{B}_2v_2\big]ds\\
               &\quad +\big[ (\hat{C}+\hat{D}_2\bar{\Theta}_2)x_{0l}^{v_2}+(\hat{B}^\top_1+\hat{D}_2\bar{\check{\Theta}})\eta^{0,v_2}+\hat{D}_1\zeta^{0,v_2}+\hat{D}_2v_2 \big]dW,\\
  d\eta^{0,v_2}&=-\big[(\hat{A}+\bar{\check{\Theta}}_2^\top\hat{F}_2)^\top\eta^{0,v_2}+C^\top\zeta^{0,v_2}+\hat{F}^\top_2\bar{\Theta}_2x_{0l}^{v_2}
+\hat{F}_2^\top v_2\big]ds+\zeta^{0,v_2}dW,\\
x_{0l}^{v_2}(0)&=0,\quad \eta^{0,v_2}(T)=0.
\end{aligned}\right.
\end{equation}
\end{mypro}

\begin{Remark}
If we consider the outcome of the closed-loop strategy in (\ref{leader closedloop system})-(\ref{lclf}) as follows: $u_2=\Theta_2\bar{x}^{\Theta_2,\check{\Theta}_2,v_2}+\check{\Theta}_2\eta^{1,\Theta_2,\check{\Theta}_2,v_2}+v_2$, which is anticipating. Therefore, inspired by \cite{Yong2002}, we consider a non-anticipating closed-loop system as follows instead:
\begin{equation}\label{nonanticipating closed-loop state}
\left\{
\begin{aligned}
d\bar{x}^{\Theta_2,\tilde{\Theta}_2,v_2}
&=\big[ (\hat{A}+\hat{B}_2\Theta_2)\bar{x}^{\Theta_2,\tilde{\Theta}_2,v_2}+\hat{F}_1\eta^{1,\Theta_2,\tilde{\Theta}_2,v_2}+\hat{B}_1\zeta^{1,\Theta_2,\tilde{\Theta}_2,v_2}+\hat{B}_2\tilde{\Theta}_2p^2\\
&\qquad +\hat{B}_2v_2+\hat{b} \big]ds\\
&\quad +\big[ (\hat{C}+\hat{D}_2\Theta_2)\bar{x}^{\Theta_2,\tilde{\Theta}_2,v_2}+\hat{B}^\top_1\eta^{1,\Theta_2,\tilde{\Theta}_2,v_2}+\hat{D}_1\zeta^{1,\Theta_2,\tilde{\Theta}_2,v_2}\\
&\qquad +\hat{D}_2\tilde{\Theta}_2p^2+\hat{D}_2v_2+\hat{\sigma} \big]dW,\\
d\eta^{1,\Theta_2,\tilde{\Theta}_2,v_2}
&=-\big[\hat{A}^\top\eta^{1,\Theta_2,\tilde{\Theta}_2,v_2}+C^\top\zeta^{1,\Theta_2,\tilde{\Theta}_2,v_2}+\hat{F}^\top_2\Theta_2\bar{x}^{\Theta_2,\tilde{\Theta}_2,v_2}\\
&\qquad +\hat{F}_2^\top\tilde{\Theta}_2p^2+\hat{F}_2^\top v_2+\hat{\beta}\big]ds+\zeta^{1,\Theta_2,\tilde{\Theta}_2,v_2}dW,\\
\bar{x}^{\Theta_2,\tilde{\Theta}_2,v_2}(0)&=x,\quad \eta^{1,\Theta_2,\tilde{\Theta}_2,v_2}(T)=g^1,
\end{aligned}
\right.
\end{equation}
with the cost functional
\begin{equation}\label{nonanticipating cost functional}
\begin{aligned}
&\hat{J}(x,\Theta_2\bar{\tilde{x}}+\tilde{\Theta}_2p^2+v_2)
=\mathbb{E}\bigg\{ \big\langle G^2\bar{\tilde{x}}(T),\bar{\tilde{x}}(T) \big\rangle+2\big\langle g^2,\bar{\tilde{x}}(T) \big\rangle\\
&+\int_0^T\Big[ \big\langle (\hat{Q}_{11}+\Theta^\top_2\hat{K}_1+\hat{K}^\top_1\Theta_2+\Theta_2^\top\hat{R}_2\Theta_2)\bar{\tilde{x}},\bar{\tilde{x}}\big\rangle+2\big\langle (\hat{Q}_{12}+\hat{K}^\top_2\Theta_2)\bar{\tilde{x}},\tilde{\eta}^1\big\rangle\\
&\quad +2\big\langle (\hat{Q}_{13}+\hat{K}^\top_3\Theta_2)\bar{\tilde{x}},\tilde{\zeta}^1\big\rangle+\big\langle \hat{Q}_{22}\tilde{\eta}^1,\tilde{\eta}^1 \big\rangle +2\big\langle \hat{Q}_{23}\tilde{\eta}^1,\tilde{\zeta}^1\big\rangle+\big\langle \hat{Q}_{33}\tilde{\zeta}^1,\tilde{\zeta}^1 \big\rangle\\
&\quad +2\big\langle (\hat{K}_1+\hat{R}_2\Theta_2)\bar{\tilde{x}},v_2 \big\rangle+2\big\langle \hat{K}_2\tilde{\eta}^1,v_2 \big\rangle+2\big\langle \hat{K}_3\tilde{\zeta}^1,v_2 \big\rangle+2\big\langle \hat{q}_3+\hat{K}^\top_3\tilde{\Theta}_2p^2,\tilde{\zeta}^1 \big\rangle\\
&\quad +2\big\langle \hat{q}_1+\hat{K}^\top_1\tilde{\Theta}_2p^2+\Theta^\top_2\hat{R}_2\tilde{\Theta}_2p^2+\Theta^\top_2\hat{\rho},\bar{\tilde{x}} \big\rangle+2\big\langle \hat{q}_2+\hat{K}^\top_2\tilde{\Theta}_2p^2,\tilde{\eta}^1 \big\rangle+\big\langle  \hat{R}_2v_2,v_2\big\rangle\\
&\quad +2\big\langle \hat{R}_2\tilde{\Theta}_2p^2+\hat{\rho},v_2 \big\rangle+\big\langle \tilde{\Theta}^\top_2\hat{R}_2\tilde{\Theta}_2p^2,p^2\big\rangle+2\big\langle\tilde{\Theta}^\top_2\hat{\rho},p^2 \big\rangle+\hat{l}\,\Big]ds\bigg\},
\end{aligned}
\end{equation}
where $p^2(\cdot)$ together with $(q^2(\cdot),k^2(\cdot))$ is the adapted solution to the following FBSDE:
\begin{equation}\label{nonanticipating adjoint equation}
\left\{\begin{aligned}
dp^2&=\bigl[ \hat{A}p^2+\hat{F}_1^\top q^2+\hat{B}_1k^2+(\hat{Q}_{12}+\hat{K}^\top_2\bar{\Theta}_2)\bar{x}+\hat{Q}^\top_{23}\bar{\zeta}^1+\hat{Q}_{22}\bar{\eta}^1+\hat{K}^\top_2\bar{\tilde{\Theta}}^\top_2p^2\\
    &\qquad +\hat{K}^\top_2\bar{v}_2+\hat{q}_2\big]ds\\
    &\quad +\bigl[\hat{C}p^2+\hat{B}^\top_1  q^2+\hat{D}^\top_1k^2+(\hat{Q}_{13}+\hat{K}^\top_3\bar{\Theta}_2)\bar{x}+\hat{Q}_{23}\bar{\eta}^1+\hat{Q}_{33}\bar{\zeta}^1+\hat{K}^\top_3\bar{\tilde{\Theta}}^\top_2p^2\\
    &\qquad +\hat{K}^\top_3\bar{v}_2+\hat{q}_3 \big]dW,\\
dq^2&=-\big[(\hat{F}^\top_2\bar{\Theta}_2+\bar{\tilde{\Theta}}^\top_2\hat{K}_1+\bar{\tilde{\Theta}}^\top_2\hat{R}_2\bar{\Theta}_2)^\top p^2 + (\hat{A}+\hat{B}_2\bar{\Theta}_2)^\top q^2+(\hat{C}+\hat{D}_2\bar{\Theta}_2)^\top k^2\\
      &\qquad +(\hat{Q}_{11}+\bar{\Theta}^\top_2\hat{K}_1+\hat{K}^\top_1\bar{\Theta}_2+\bar{\Theta}^\top_2\hat{R}_2\bar{\Theta}_2)\bar{x}+(\hat{Q}^\top_{12}+\bar{\Theta}_2^\top\hat{K}_2)\bar{\eta}^1\\
      &\qquad +(\hat{Q}^\top_{13}+\bar{\Theta}^\top_2\hat{K}_3)\bar{\zeta}^1+(\hat{K}^\top_1+\bar{\Theta}^\top_2\hat{R}_2)\bar{v}_2 +\hat{q}_1+\bar{\Theta}^\top_2\hat{\rho}\big]ds+k^2dW,\\
p^2(0)&=0,\,\,\,q^2(T)=G^2\bar{x}(T)+g^2.
\end{aligned}\right.
\end{equation}
And the following stationary condition holds:
\begin{equation}\label{nonanticipating stationary condition}
\begin{aligned}
&(\hat{F}_2+\hat{R}_2\bar{\tilde{\Theta}}_2)p^2+\hat{B}^\top_2q^2+\hat{D}^\top_2k^2+(\hat{K}_1+\hat{R}_2\bar{\Theta}_2)\bar{x}+\hat{K}_2\bar{\eta}^1\\
&\quad +\hat{K}_3\bar{\zeta}^1+\hat{R}_2\bar{v}_2+\hat{\rho}=0,\quad a.e.,\, \mathbb{P}\mbox{-}a.s.,
\end{aligned}
\end{equation}
where $(\bar{x}(\cdot),\bar{\eta}^1(\cdot),\bar{\zeta}^1(\cdot)) \equiv (\bar{x}^{\bar{\Theta}_2,\bar{\tilde{\Theta}}_2,\bar{v}_2}(\cdot),\bar{\eta}^{1,\bar{\Theta}_2,\bar{\tilde{\Theta}}_2,\bar{v}_2}(\cdot),\bar{\zeta}^{1,\bar{\Theta}_2,\bar{\tilde{\Theta}}_2,\bar{v}_2}(\cdot))$ is the optimal triple of the closed-loop system (\ref{nonanticipating closed-loop state})-(\ref{nonanticipating cost functional}).
\end{Remark}

Making use of the stationary condition in (\ref{nonanticipating stationary condition}), we rewrite the BSDE in (\ref{nonanticipating adjoint equation}) as follows:
\begin{equation}
\begin{aligned}
 dq^2&=-\big[(\hat{F}^\top_2\bar{\Theta}_2+\bar{\tilde{\Theta}}^\top_2\hat{K}_1+\bar{\tilde{\Theta}}^\top_2\hat{R}_2\bar{\Theta}_2)^\top p^2 + (\hat{A}+\hat{B}_2\bar{\Theta}_2)^\top q^2
      +(\hat{C}+\hat{D}_2\bar{\Theta}_2)^\top k^2\\
 &\qquad +(\hat{Q}_{11}+\bar{\Theta}^\top_2\hat{K}_1+\hat{K}^\top_1\bar{\Theta}_2+\bar{\Theta}^\top_2\hat{R}_2\bar{\Theta}_2)\bar{x}+(\hat{Q}^\top_{12}+\bar{\Theta}_2^\top\hat{K}_2)\bar{\eta}^1\\
 &\qquad +(\hat{Q}^\top_{13}+\bar{\Theta}^\top_2\hat{K}_3)\bar{\zeta}^1+(\hat{K}^\top_1+\bar{\Theta}^\top_2\hat{R}_2)\bar{v}_2
  +\hat{q}_1+\bar{\Theta}^\top_2\hat{\rho}\big]ds+k^2dW,\\
 &=-\big\{ \bar{\Theta}^\top_2\big[(\hat{F}_2+\hat{R}_2\bar{\tilde{\Theta}}_2)p^2+\hat{B}^\top_2q^2+\hat{D}^\top_2k^2+(\hat{K}_1+\hat{R}_2\bar{\Theta}_2)\bar{x}+\hat{K}_2\bar{\eta}^1+\hat{K}_3\bar{\zeta}^1\\
 &\qquad +\hat{R}_2\bar{v}_2+\hat{\rho}\big]+\hat{A}^\top q^2+\hat{C}^\top k^2+\hat{K}^\top_1\bar{\tilde{\Theta}}_2p^2+(\hat{Q}_{11}+\hat{K}^\top_1\bar{\Theta}_2)\bar{x}\\
 &\qquad +\hat{Q}^\top_{12}\bar{\eta}^1+\hat{Q}^\top_{13}\bar{\zeta}^1+\hat{K}^\top_1\bar{v}_2+\hat{q}_1 \big\}ds +k^2dW\\
 &=-\big\{\hat{A}^\top q^2+\hat{C}^\top k^2+\hat{K}^\top_1\bar{\tilde{\Theta}}_2p^2+(\hat{Q}_{11}+\hat{K}^\top_1\bar{\Theta}_2)\bar{x}\\
 &\qquad +\hat{Q}^\top_{12}\bar{\eta}^1+\hat{Q}^\top_{13}\bar{\zeta}^1+\hat{K}^\top_1\bar{v}_2+\hat{q}_1 \big\}ds +k^2dW.
\end{aligned}
\end{equation}
For convenience, we write the state equation and adjoint equation together, and obtain
\begin{equation}\label{leader closedloop optimal system simple}
\left\{
\begin{aligned}
d\bar{x}&=\big[ (\hat{A}+\hat{B}_2\bar{\Theta}_2)\bar{x}+\hat{F}_1\bar{\eta}^1+\hat{B}_1\bar{\zeta}^1+\hat{B}_2\bar{\tilde{\Theta}}_2p^2+\hat{B}_2\bar{v}_2+\hat{b} \big]ds\\
        &\quad +\big[ (\hat{C}+\hat{D}_2\bar{\Theta}_2)\bar{x}+\hat{B}^\top_1\bar{\eta}^1+\hat{D}_1\bar{\zeta}^1+\hat{D}_2\bar{\tilde{\Theta}}_2p^2+\hat{D}_2\bar{v}_2+\hat{\sigma} \big]dW,\\
d\bar{\eta}^1&=-\big[\hat{A}^\top\bar{\eta}^1+\hat{C}^\top\bar{\zeta}^1+\hat{F}^\top_2\bar{\Theta}_2\bar{x}+\hat{F}_2^\top\bar{\tilde{\Theta}}_2p^2+\hat{F}_2^\top \bar{v}_2+\hat{\beta}\big]ds+\bar{\zeta}^1dW,\\
dp^2&=\bigl[ \hat{A}p^2+\hat{F}_1^\top q^2+\hat{B}_1k^2+(\hat{Q}_{12}+\hat{K}^\top_2\bar{\Theta}_2)\bar{x}+\hat{Q}^\top_{23}\bar{\zeta}^1+\hat{Q}_{22}\bar{\eta}^1\\
    &\qquad +\hat{K}^\top_2\bar{\tilde{\Theta}}^\top_2p^2+\hat{K}^\top_2\bar{v}_2+\hat{q}_2\big]ds\\
    &\quad +\bigl[\hat{C}p^2+\hat{B}^\top_1  q^2+\hat{D}^\top_1k^2+(\hat{Q}_{13}+\hat{K}^\top_3\bar{\Theta}_2)\bar{x}+\hat{Q}_{23}\bar{\eta}^1\\
    &\qquad +\hat{Q}_{33}\bar{\zeta}^1+\hat{K}^\top_3\bar{\tilde{\Theta}}^\top_2p^2+\hat{K}^\top_3\bar{v}_2+\hat{q}_3 \big]dW,\\
dq^2&=-\big\{\hat{A}^\top q^2+\hat{C}^\top k^2+\hat{K}^\top_1\bar{\tilde{\Theta}}_2p^2+(\hat{Q}_{11}+\hat{K}^\top_1\bar{\Theta}_2)\bar{x}\\
    &\qquad +\hat{Q}^\top_{12}\bar{\eta}^1+\hat{Q}^\top_{13}\bar{\zeta}^1+\hat{K}^\top_1\bar{v}_2+\hat{q}_1 \big\}ds+k^2dW,\\
\bar{x}(0)&=x,\quad \bar{\eta}^1(T)=g^1,\quad p^2(0)=0,\quad q^2(T)=G^2\bar{x}(T)+g^2.
\end{aligned}
\right.
\end{equation}
Note that the above is a coupled FBSDEs system which is further coupled through the \eqref{nonanticipating stationary condition}. Next, let us set
\begin{equation}
X =\left(\begin{matrix} \bar{x}  \\ p^2  \end{matrix}\right),\,\,\,Y =\left(\begin{matrix} q^2  \\ \bar{\eta}^1  \end{matrix}\right),\,\,\,Z =\left(\begin{matrix} k^2  \\ \bar{\zeta}^1  \end{matrix}\right),\,\,\,
X_0 =\left(\begin{matrix} x  \\ 0  \end{matrix}\right),\,\,\,
\bar{\boldsymbol{\Theta}}_2=\left(\begin{matrix} \bar{\Theta}_2 & \bar{\tilde{\Theta}}_2 \end{matrix}\right),
\end{equation}
and
\begin{equation*}
\begin{cases}
\mathcal{A}=\left(\begin{matrix} \hat{A} & 0 \\ \hat{Q}_{12} & \hat{A} \end{matrix}\right),\quad
\mathcal{B}_2=\left(\begin{matrix} \hat{B}_2 \\ \hat{K}^\top_2 \end{matrix}\right),\quad
\mathcal{F}_1=\left(\begin{matrix} 0 & \hat{F}_1 \\ \hat{F}_1^\top& \hat{Q}_{22} \end{matrix}\right),\quad
\mathcal{B}_1=\left(\begin{matrix} 0 & \hat{B}_1 \\ \hat{B}_1^\top& \hat{Q}^\top_{23} \end{matrix}\right),\quad
\boldsymbol{b}=\left(\begin{matrix} \hat{b} \\ \hat{q}_2 \end{matrix}\right),\quad\\
\mathcal{C}=\left(\begin{matrix} \hat{C} & 0  \\ \hat{Q}_{13} & \hat{C} \end{matrix}\right),\quad
\mathcal{D}_2=\left(\begin{matrix} \hat{D}_2 \\ \hat{K}_3^\top \end{matrix}\right),\quad
\mathcal{D}_1=\left(\begin{matrix} 0 & \hat{D}_1 \\ \hat{D}_1^\top& \hat{Q}_{33} \end{matrix}\right),\quad
\boldsymbol{\sigma}=\left(\begin{matrix} \hat{\sigma} \\ \hat{q}_3 \end{matrix}\right),\quad
\mathcal{Q}_2=\left(\begin{matrix} \hat{Q}_{11}& 0 \\ 0 & 0 \end{matrix}\right),\quad\\
\mathcal{F}_2=\left(\begin{matrix} \hat{K}_1& \hat{F}_2  \end{matrix}\right),\quad
\boldsymbol{\beta}=\left(\begin{matrix} \hat{q}_1 \\ \hat{\beta} \end{matrix}\right),\quad
\mathcal{G}_2=\left(\begin{matrix} G_2 & 0\\ 0& 0 \end{matrix}\right),\quad
\boldsymbol{g}=\left(\begin{matrix} g^1 \\ g^2  \end{matrix}\right).
\end{cases}
\end{equation*}
Then (\ref{leader closedloop optimal system simple}) is equivalent to the following FBSDE:
\begin{equation}\label{High dimension optimal system simple}
\left\{
\begin{aligned}
dX&=\big[ (\mathcal{A}+\mathcal{B}_2\bar{\boldsymbol{\Theta}}_2)X+\mathcal{F}_1Y+\mathcal{B}_1Z+\mathcal{B}_2\bar{v}_2+\boldsymbol{b}  \big]ds\\
&\qquad +\big[ (\mathcal{C}+\mathcal{D}_2\bar{\boldsymbol{\Theta}}_2)X+\mathcal{B}^\top_1Y+\mathcal{D}_1Z+\mathcal{D}_2\bar{v}_2+\boldsymbol{\sigma}  \big]dW,\\
dY&=-\big[ (\mathcal{Q}_2+\mathcal{F}^\top_2\bar{\boldsymbol{\Theta}}_2)X+\mathcal{A}^\top Y+\mathcal{C}^\top Z+\mathcal{F}^\top_2\bar{v}_2+\boldsymbol{\beta} \big]ds+ZdW,\\
X(0)&=X_0,\,\,\,Y(T)=\mathcal{G}_2X(T)+\boldsymbol{g},
\end{aligned}
\right.
\end{equation}
whose adapted solution is $(X(\cdot),Y(\cdot),Z(\cdot))\in L^2_{\mathbb{F}}(0,T;\mathbb{R}^{2n})\times L^2_{\mathbb{F}}(0,T;\mathbb{R}^{2n})\times L^2_{\mathbb{F}}(0,T;\mathbb{R}^{2n})$, with:
\begin{equation}\label{High dimension stationary}
(\hat{R}_2\bar{\boldsymbol{\Theta}}_2+\mathcal{F}_2)X+\mathcal{B}^\top_2Y+\mathcal{D}^\top_2Z+\hat{R}_2\bar{v}_2+\hat{\rho}=0,\quad a.e.,\, \mathbb{P}\mbox{-}a.s.
\end{equation}

For the closed-loop optimal strategies of the leader, we have the following result.
\begin{mythm}\label{Th-cl-l}
Let (H1)-(H2) hold, if Problem (SLQ)$_l$ admits a closed-loop optimal strategy $(\bar{\boldsymbol{\Theta}}_2(\cdot),\bar{v}_2(\cdot))\in \mathcal{Q}_2[0,T] \times \mathcal{Q}_2[0,T] \times \mathcal{U}_2[0,T]$, then it admits the following representation:
\begin{equation}\label{close-loop leader}
\left\{
\begin{aligned}
\bar{\boldsymbol{\Theta}}_2&=-\big[\hat{R}_2+\mathcal{D}^\top_2(I-P\mathcal{D}_1)^{-1}P\mathcal{D}_2\big]^{-1}
                            \big[\mathcal{B}^\top_2P+\mathcal{F}_2+\mathcal{D}^\top_2(I-P\mathcal{D}_1)^{-1}P(\mathcal{C}+\mathcal{B}^\top_1P)\big],\\
	              \bar{v}_2&=-\big[\hat{R}_2+\mathcal{D}^\top_2(I-P\mathcal{D}_1)^{-1}P\mathcal{D}_2\big]^{-1}\Big\{\big[\mathcal{B}^\top_2+\mathcal{D}^\top_2(I-P\mathcal{D}_1)^{-1}P\mathcal{B}^\top_1\big]\eta\\
                       	   &\qquad +\mathcal{D}^\top_2(I-P\mathcal{D}_1)^{-1}\zeta+\mathcal{D}^\top_2(I-P\mathcal{D}_1)^{-1}P\boldsymbol{\sigma}+\hat{\rho} \Big\},\qquad a.e.,\, \mathbb{P}\mbox{-}a.s.
\end{aligned}
\right.
\end{equation}
where $P(\cdot)\equiv\left(\begin{matrix} P_1(\cdot)&P_2(\cdot)  \\ P_3(\cdot)&P_4(\cdot)  \end{matrix}\right)\in C([0,T];\mathbb{R}^{2n \times 2n})$ is the solution to the Riccati equation:
\begin{equation}\label{leader RE}\left\{
\begin{aligned}
	0&=\dot{P}+\mathcal{A}^\top P+P\mathcal{A}+P\mathcal{F}_1P+\mathcal{Q}_2-(P\mathcal{B}_2+\mathcal{F}^\top_2)\big[\hat{R}_2+\mathcal{D}^\top_2(I-P\mathcal{D}_1)^{-1}P\mathcal{D}_2\big]^{-1}\\
	&\quad \times\big[\mathcal{B}^\top_2P+\mathcal{F}_2+\mathcal{D}^\top_2(I-P\mathcal{D}_1)^{-1}P(\mathcal{C}+\mathcal{B}^\top_1P)\big]\\
	&\quad +(\mathcal{C}^\top+P\mathcal{B}_1)(I-P\mathcal{D}_1)^{-1}P\Big\{ \mathcal{C}+\mathcal{B}^\top_1P-\mathcal{D}_2\big[\hat{R}_2+\mathcal{D}^\top_2(I-P\mathcal{D}_1)^{-1}P\mathcal{D}_2\big]^{-1}\\
	&\quad \times\big[\mathcal{B}^\top_2P+\mathcal{F}_2+\mathcal{D}^\top_2(I-P\mathcal{D}_1)^{-1}P(\mathcal{C}+\mathcal{B}^\top_1P)\big] \Big\},\\
&\hspace{-1mm}P(T)=\mathcal{G}_2,
\end{aligned}\right.
\end{equation}
and $(\eta(\cdot),\zeta(\cdot))$ is the adapted solution to the following BSDE:
\begin{equation}\nonumber\left\{
\begin{aligned}
d\eta&=-\bigg\{\Big[\mathcal{A}^\top+P\mathcal{F}_1+(\mathcal{C}^\top+P\mathcal{B}_1)(I-P\mathcal{D}_1)^{-1}P\mathcal{B}^\top_1\\
&\qquad\quad -\big[(\mathcal{C}^\top+P\mathcal{B}_1)(I-P\mathcal{D}_1)^{-1}P\mathcal{D}_2+\mathcal{F}^\top_2+P\mathcal{B}_2\big]\\
&\qquad\quad \times\big[\hat{R}_2+\mathcal{D}^\top_2(I-P\mathcal{D}_1)^{-1}P\mathcal{D}_2\big]^{-1}\big[\mathcal{B}^\top_2+\mathcal{D}^\top_2(I-P\mathcal{D}_1)^{-1}P\mathcal{B}^\top_1\big]\Big]\eta\\
&\qquad +\Big[(\mathcal{C}^\top+P\mathcal{B}_1)(I-P\mathcal{D}_1)^{-1}-\big[(\mathcal{C}^\top+P\mathcal{B}_1)(I-P\mathcal{D}_1)^{-1}P\mathcal{D}_2+\mathcal{F}^\top_2+P\mathcal{B}_2\big]\\
&\qquad\quad \times\big[\hat{R}_2+\mathcal{D}^\top_2(I-P\mathcal{D}_1)^{-1}P\mathcal{D}_2\big]^{-1}\mathcal{D}^\top_2(I-P\mathcal{D}_1)^{-1}\Big]\zeta\\
&\qquad -\big[(\mathcal{C}^\top+P\mathcal{B}_1)(I-P\mathcal{D}_1)^{-1}P\mathcal{D}_2+\mathcal{F}^\top_2+P\mathcal{B}_2\big]\big[\mathcal{D}^\top_2(I-P\mathcal{D}_1)^{-1}P\boldsymbol{\sigma}+\hat{\rho}\big]\\
&\qquad +(\mathcal{C}^\top+P\mathcal{B}_1)(I-P\mathcal{D}_1)^{-1}P\boldsymbol{\sigma}+\boldsymbol{\beta}+P\boldsymbol{b}\bigg\}ds+\zeta dW,\\
\eta(T)&=\boldsymbol{g}.
\end{aligned}\right.
\end{equation}
In this case, the optimal control of the leader is $\bar{u}_2(\cdot)=\boldsymbol\Theta_2(\cdot)X(\cdot)+\bar{v}_2(\cdot)$, where $X(\cdot)$ is the solution to the following SDE:
\begin{equation}\label{nonanticipating X}\left\{
\begin{aligned}
  dX&=\bigg\{\Big[\mathcal{A}+\mathcal{F}_1P+\mathcal{B}_1(I-P\mathcal{D}_1)^{-1}P(\mathcal{C}+\mathcal{B}_1^\top P)-\big[\mathcal{B}_1(I-P\mathcal{D}_1)^{-1}P\mathcal{D}_2+\mathcal{B}_2\big]\\
	&\qquad \times\big[\hat{R}_2+\mathcal{D}_2^\top(I-P\mathcal{D}_1)^{-1}P\mathcal{D}_2\big]^{-1}\big[\mathcal{B}_2^\top P+\mathcal{F}_2+\mathcal{D}_2^\top(I-P\mathcal{D}_1)^{-1}P(\mathcal{C}+\mathcal{B}_1^\top P)\big]\Big]X\\
	&\quad +\Big[\mathcal{F}_1+\mathcal{B}_1(I-P\mathcal{D}_1)^{-1}P\mathcal{B}_1^\top-\big[\mathcal{B}_1(I-P\mathcal{D}_1)^{-1}P\mathcal{D}_2+\mathcal{B}_2\big]\\
	&\qquad \times\big[\hat{R}_2+\mathcal{D}_2^\top(I-P\mathcal{D}_1)^{-1}P\mathcal{D}_2]^{-1}\big[\mathcal{B}_2^\top+\mathcal{D}_2^\top(I-P\mathcal{D}_1)^{-1}P\mathcal{B}_1^\top\big]\Big]\eta\\
	&\quad +\Big[\mathcal{B}_1(I-P\mathcal{D}_1)^{-1}-\big[\mathcal{B}_1(I-P\mathcal{D}_1)^{-1}P\mathcal{D}_2+\mathcal{B}_2\big]\\
	&\qquad \times\big[\hat{R}_2+\mathcal{D}_2^\top(I-P\mathcal{D}_1)^{-1}P\mathcal{D}_2\big]^{-1}\mathcal{D}_2^\top(I-P\mathcal{D}_1)^{-1}\Big]\zeta\\
	&\quad +\Big[\mathcal{B}_1(I-P\mathcal{D}_1)^{-1}P-\big[\mathcal{B}_1(I-P\mathcal{D}_1)^{-1}P\mathcal{D}_2+\mathcal{B}_2\big]\\
	&\qquad \times\big[\hat{R}_2+\mathcal{D}_2^\top(I-P\mathcal{D}_1)^{-1}P\mathcal{D}_2\big]^{-1}\mathcal{D}_2^\top(I-P\mathcal{D}_1)^{-1}P\Big]\boldsymbol{\sigma}\\
	&\qquad +\boldsymbol{b}-\big[\mathcal{B}_1(I-P\mathcal{D}_1)^{-1}P\mathcal{D}_2+\mathcal{B}_2\big]\big[\hat{R}_2+\mathcal{D}_2^\top(I-P\mathcal{D}_1)^{-1}P\mathcal{D}_2\big]^{-1}\hat{\rho}\bigg\}ds\\
	&+\bigg\{\Big[\mathcal{C}+\mathcal{B}_1^\top P+\mathcal{D}_1(I-P\mathcal{D}_1)^{-1}P(\mathcal{C}+\mathcal{B}_1^\top P)-\big[\mathcal{D}_1(I-P\mathcal{D}_1)^{-1}P\mathcal{D}_2+\mathcal{D}_2\big]\\
	&\qquad \times\big[\hat{R}_2+\mathcal{D}_2^\top(I-P\mathcal{D}_1)^{-1}P\mathcal{D}_2\big]^{-1}\big[\mathcal{B}_2^\top P+\mathcal{F}_2+\mathcal{D}_2^\top(I-P\mathcal{D}_1)^{-1}P(\mathcal{C}+\mathcal{B}_1^\top P)\big]\Big]X\\
	&\quad +\Big[\mathcal{B}^\top_1+\mathcal{D}_1(I-P\mathcal{D}_1)^{-1}P\mathcal{B}_1^\top-\big[\mathcal{D}_1(I-P\mathcal{D}_1)^{-1}P\mathcal{D}_2+\mathcal{D}_2\big]\\
	&\qquad \times\big[\hat{R}_2+\mathcal{D}_2^\top(I-P\mathcal{D}_1)^{-1}P\mathcal{D}_2\big]^{-1}\big[\mathcal{B}_2^\top+\mathcal{D}_2^\top(I-P\mathcal{D}_1)^{-1}P\mathcal{B}_1^\top\big]\Big]\eta\\
	&\quad +\Big[\mathcal{D}_1(I-P\mathcal{D}_1)^{-1}-\big[\mathcal{D}_1(I-P\mathcal{D}_1)^{-1}P\mathcal{D}_2+\mathcal{D}_2\big]\\
	&\qquad \times\big[\hat{R}_2+\mathcal{D}_2^\top(I-P\mathcal{D}_1)^{-1}P\mathcal{D}_2\big]^{-1}\mathcal{D}_2^\top(I-P\mathcal{D}_1)^{-1}\Big]\zeta\\
	&\quad +\Big[I+\mathcal{D}_1(I-P\mathcal{D}_1)^{-1}P-\big[\mathcal{D}_1(I-P\mathcal{D}_1)^{-1}P\mathcal{D}_2+\mathcal{D}_2\big]\\
	&\qquad \times\big[\hat{R}_2+\mathcal{D}_2^\top(I-P\mathcal{D}_1)^{-1}P\mathcal{D}_2\big]^{-1}\mathcal{D}_2^\top(I-P\mathcal{D}_1)^{-1}P\Big]\boldsymbol{\sigma}\\
	&\qquad-\big[\mathcal{D}_1(I-P\mathcal{D}_1)^{-1}P\mathcal{D}_2+\mathcal{D}_2\big]\big[\hat{R}_2+\mathcal{D}_2^\top(I-P\mathcal{D}_1)^{-1}P\mathcal{D}_2\big]^{-1}\hat{\rho}\bigg\}dW,\\
	X(0)&=X_0.
\end{aligned}\right.
\end{equation}
\end{mythm}

\textit{Proof.} Let $(\bar{\boldsymbol{\Theta}}_2(\cdot),\bar{v}_2(\cdot))$ be a closed-loop optimal strategy of Problem (SLQ)$_l$ over $[0,T]$. Since (\ref{High dimension optimal system simple}) admits a solution for each $X_0 \in \mathbb{R}^{2n}$, and $(\bar{\boldsymbol{\Theta}}_2(\cdot),\bar{v}_2(\cdot))$ is independent of $X_0$, by substracting solutions corresponding to $X_0$ and $0$, the later from the former, we see that for any $X_0 \in \mathbb{R}^{2n}$, the following FBSDE admits an adapted solution $(\tilde{X}(\cdot),\tilde{Y}(\cdot),\tilde{Z}(\cdot))$:
\begin{equation}
\left\{
\begin{aligned}
d\tilde{X}&=\big[ (\mathcal{A}+\mathcal{B}_2\bar{\boldsymbol{\Theta}}_2)\tilde{X}+\mathcal{F}_1\tilde{Y}+\mathcal{B}_1\tilde{Z}\big]ds\\
&\qquad +\big[ (\mathcal{C}+\mathcal{D}_2\bar{\boldsymbol{\Theta}}_2)\tilde{X}+\mathcal{B}^\top_1\tilde{Y}+\mathcal{D}_1\tilde{Z}\big]dW,\\
d\tilde{Y}&=-\big[ (\mathcal{Q}_2+\mathcal{F}^\top_2\bar{\boldsymbol{\Theta}}_2)\tilde{X}+\mathcal{A}^\top\tilde{Y}+\mathcal{C}^\top \tilde{Z}\big]ds+\tilde{Z}dW,\\
\tilde{X}(0)&=X_0,\,\,\,\tilde{Y}(T)=\mathcal{G}_2\tilde{X}(T).\\
&\hspace{-8mm}(\hat{R}_2\bar{\boldsymbol{\Theta}}_2+\mathcal{F}_2)\tilde{X}+\mathcal{B}^\top_2\tilde{Y}+\mathcal{D}^\top_2\tilde{Z}=0,\quad a.e.,\, \mathbb{P}\mbox{-}a.s.
\end{aligned}\right.
\end{equation}

Now, we let
\begin{equation}
\left\{
\begin{aligned}
d\mathbb{X}&=\big[ (\mathcal{A}+\mathcal{B}_2\bar{\boldsymbol{\Theta}}_2)\mathbb{X}+\mathcal{F}_1\mathbb{Y}+\mathcal{B}_1\mathbb{Z}\big]ds\\
&\qquad +\big[ (\mathcal{C}+\mathcal{D}_2\bar{\boldsymbol{\Theta}}_2)\mathbb{X}+\mathcal{B}^\top_1\mathbb{Y}+\mathcal{D}_1\mathbb{Z}\big]dW,\\
d\mathbb{Y}&=-\big[ (\mathcal{Q}_2+\mathcal{F}^\top_2\bar{\boldsymbol{\Theta}}_2)\mathbb{X}+\mathcal{A}^\top\mathbb{Y}+\mathcal{C}^\top \mathbb{Z}\big]ds+\mathbb{Z}dW,\\
\mathbb{X}(0)&=I_{2n \times 2n},\,\,\,\mathbb{Y}(T)=\mathcal{G}_2\mathbb{X}(T).\\
\end{aligned}
\right.
\end{equation}
Clearly, $\mathbb{X}(\cdot),\mathbb{Y}(\cdot),\mathbb{Z}(\cdot)$ are all well-defined $(2n\times 2n)$-matrix valued processes. Further,
\begin{equation}\label{simply stationary}
(\hat{R}_2\bar{\boldsymbol{\Theta}}_2+\mathcal{F}_2)\mathbb{X}+\mathcal{B}^\top_2\mathbb{Y}+\mathcal{D}^\top_2\mathbb{Z}=0,\quad a.e.,\, \mathbb{P}\mbox{-}a.s.
\end{equation}
Drawing on the method of Yong \cite{Yong2006}, we can check that $\mathbb{X}(\cdot)^{-1}$ exists and satisfies the SDE:
\begin{equation}
\left\{
\begin{aligned}
  d\mathbb{X}^{-1}&=\Big\{-\mathbb{X}^{-1}\big[(\mathcal{A}+\mathcal{B}_2\bar{\boldsymbol{\Theta}}_2)\mathbb{X}+\mathcal{F}_1\mathbb{Y}+\mathcal{B}_1\mathbb{Z}\big]\mathbb{X}^{-1}
                   +\mathbb{X}^{-1}\big[(\mathcal{C}+\mathcal{D}_2\bar{\boldsymbol{\Theta}}_2)\mathbb{X}\\
                  &\qquad +\mathcal{B}^\top_1\mathbb{Y}+\mathcal{D}_1\mathbb{Z}\big]
                   \mathbb{X}^{-1}\big[(\mathcal{C}+\mathcal{D}_2\bar{\boldsymbol{\Theta}}_2)\mathbb{X}+\mathcal{B}^\top_1\mathbb{Y}+\mathcal{D}_1\mathbb{Z}\big]\mathbb{X}^{-1} \Big\}ds\\
                  &\quad -\mathbb{X}^{-1}[(\mathcal{C}+\mathcal{D}_2\bar{\boldsymbol{\Theta}}_2)\mathbb{X}+\mathcal{B}^\top_1\mathbb{Y}+\mathcal{D}_1\mathbb{Z}]\mathbb{X}^{-1}dW,\\
\mathbb{X}(0)^{-1}&=I_{2n\times 2n}.
\end{aligned}
\right.
\end{equation}
We define
\begin{equation}\label{PX-1}
P(\cdot)=\mathbb{Y}(\cdot)\mathbb{X}(\cdot)^{-1},\qquad \Pi(\cdot)=\mathbb{Z}(\cdot)\mathbb{X}(\cdot)^{-1}.
\end{equation}
By It\^o's formula, we obtain
\begin{equation}\nonumber
\begin{aligned}
dP&=\Big\{ -(\mathcal{Q}_2+\mathcal{F}^\top_2\bar{\boldsymbol{\Theta}}_2)-\mathcal{A}^\top P-\mathcal{C}^\top\Pi
   -P(\mathcal{A}+\mathcal{B}_2\bar{\boldsymbol{\Theta}}_2)-P\mathcal{F}_1P-P\mathcal{B}_1\Pi\\
  &\qquad +P\big[\mathcal{C}+\mathcal{D}_2\bar{\boldsymbol{\Theta}}_2+\mathcal{B}^\top_1P+\mathcal{D}_1\Pi\big]\big[\mathcal{C}+\mathcal{D}_2\bar{\boldsymbol{\Theta}}_2+\mathcal{B}^\top_1P+\mathcal{D}_1\Pi\big]\\
  &\qquad -\Pi\big[\mathcal{C}+\mathcal{D}_2\bar{\boldsymbol{\Theta}}_2+\mathcal{B}^\top_1P+\mathcal{D}_1\Pi\big]\Big\}ds
   +\Big\{ \Pi-P\big[\mathcal{C}+\mathcal{D}_2\bar{\boldsymbol{\Theta}}_2+\mathcal{B}^\top_1P+\mathcal{D}_1\Pi\big] \Big\}dW.
\end{aligned}
\end{equation}
Let
\begin{equation}\nonumber
\Lambda \triangleq (I-P\mathcal{D}_1)\Pi-P\big[\mathcal{C}+\mathcal{D}_2\bar{\boldsymbol{\Theta}}_2+\mathcal{B}^\top_1P\big],
\end{equation}
which leads to
\begin{equation}\nonumber
\begin{aligned}
dP&=-\Big\{ (\mathcal{Q}_2+\mathcal{F}^\top_2\bar{\boldsymbol{\Theta}}_2)+\mathcal{A}^\top P+\mathcal{C}^\top\Pi+P(\mathcal{A}+\mathcal{B}_2\bar{\boldsymbol{\Theta}}_2)+P\mathcal{F}_1P+P\mathcal{B}_1\Pi\\
  &\qquad +(\Lambda-\Pi)\big[\mathcal{C}+\mathcal{D}_2\bar{\boldsymbol{\Theta}}_2+\mathcal{B}^\top_1P+\mathcal{D}_1\Pi\big]
   +\Pi\big[\mathcal{C}+\mathcal{D}_2\bar{\boldsymbol{\Theta}}_2+\mathcal{B}^\top_1P+\mathcal{D}_1\Pi\big]\Big\}ds+\Lambda dW\\
  &=-\Big\{ \mathcal{Q}_2+\mathcal{F}^\top_2\bar{\boldsymbol{\Theta}}_2+\mathcal{A}^\top P+\mathcal{C}^\top\Pi
   +P(\mathcal{A}+\mathcal{B}_2\bar{\boldsymbol{\Theta}}_2)+P\mathcal{F}_1P+P\mathcal{B}_1\Pi\\
  &\qquad  +\Lambda\big[\mathcal{C}+\mathcal{D}_2\bar{\boldsymbol{\Theta}}_2+\mathcal{B}^\top_1P+\mathcal{D}_1\Pi\big]\Big\}ds+\Lambda dW,
\end{aligned}
\end{equation}
and $P(T)=\mathcal{G}_2$. Thus, $(P(\cdot),\Lambda(\cdot))$ is the adapted solution to a BSDE with deterministic coefficients. Hence, $P(\cdot)$ is deterministic and $\Lambda(\cdot)=0$ which means
\begin{equation}\nonumber
\Pi=(I-P\mathcal{D}_1)^{-1}P\big[\mathcal{C}+\mathcal{D}_2\bar{\boldsymbol{\Theta}}_2+\mathcal{B}^\top_1P\big].
\end{equation}
Therefore, we get
\begin{equation}\label{leader riccati}
\begin{aligned}
&\dot{P}+\mathcal{A}^\top P+P\mathcal{A}+P\mathcal{F}_1P+\mathcal{Q}_2+(\mathcal{C}^\top+P\mathcal{B}_1)(I-P\mathcal{D}_1)^{-1}P(\mathcal{C}+\mathcal{B}^\top_1P)\\
&\quad +[P\mathcal{B}_2+\mathcal{F}^\top_2+(\mathcal{C}^\top+P\mathcal{B}_1)(I-P\mathcal{D}_1)^{-1}P\mathcal{D}_2]\bar{\boldsymbol{\Theta}}_2=0.
\end{aligned}
\end{equation}
Moreover, (\ref{simply stationary}) and (\ref{PX-1}) imply
\begin{equation}\label{leader Theta condition}
\begin{aligned}
0&=\hat{R}_2\bar{\boldsymbol{\Theta}}_2+\mathcal{F}_2+\mathcal{B}^\top_2P+\mathcal{D}^\top_2\Pi\\
 &=\hat{R}_2\bar{\boldsymbol{\Theta}}_2+\mathcal{F}_2+\mathcal{B}^\top_2P+\mathcal{D}^\top_2(I-P\mathcal{D}_1)^{-1}P\big[\mathcal{C}+\mathcal{D}_2\bar{\boldsymbol{\Theta}}_2+\mathcal{B}^\top_1P\big]\\
 &=\mathcal{B}^\top_2P+\mathcal{F}_2+\mathcal{D}^\top_2(I-P\mathcal{D}_1)^{-1}P(\mathcal{C}+\mathcal{B}^\top_1P)+\big[\hat{R}_2+\mathcal{D}^\top_2(I-P\mathcal{D}_1)^{-1}P\mathcal{D}_2\big]\bar{\boldsymbol{\Theta}}_2.
\end{aligned}
\end{equation}
Thus
\begin{equation}\label{Theta2}
\bar{\boldsymbol{\Theta}}_2=-\big[\hat{R}_2+\mathcal{D}^\top_2(I-P\mathcal{D}_1)^{-1}P\mathcal{D}_2\big]^{-1}\big[\mathcal{B}^\top_2P+\mathcal{F}_2+\mathcal{D}^\top_2(I-P\mathcal{D}_1)^{-1}P(\mathcal{C}+\mathcal{B}^\top_1P)\big].
\end{equation}
Plugging the above into (\ref{leader riccati}), we obtain Riccati equation in (\ref{leader RE}). To determine $\bar{v}_2(\cdot)$, we define
\begin{equation}\label{relationship}\left\{
\begin{aligned}
  \eta&\triangleq Y-PX,\\
 \zeta&\triangleq Z-P\big[(\mathcal{C}+\mathcal{D}_2\bar{\boldsymbol{\Theta}}_2)X+\mathcal{B}^\top_1Y+\mathcal{D}_1Z+\mathcal{D}_2\bar{v}_2+\boldsymbol{\sigma}\big].
\end{aligned}\right.
\end{equation}
Consequently,
\begin{equation}\label{eta}
\begin{aligned}
d\eta&=\Big\{ -(\mathcal{Q}_2+\mathcal{F}^\top_2\bar{\boldsymbol{\Theta}}_2)X-\mathcal{A}^\top Y-\mathcal{C}^\top Z-\mathcal{F}^\top_2\bar{v}_2-\boldsymbol{\beta}+\mathcal{A}^\top PX+P\mathcal{A}X\\
     &\qquad+P\mathcal{F}_1PX+\mathcal{Q}_2X+(\mathcal{C}^\top+P\mathcal{B}_1)(I-P\mathcal{D}_1)^{-1}P(\mathcal{C}+\mathcal{B}^\top_1P)X\\
&\qquad +[P\mathcal{B}_2+\mathcal{F}^\top_2+(\mathcal{C}^\top+P\mathcal{B}_1)(I-P\mathcal{D}_1)^{-1}P\mathcal{D}_2]\bar{\boldsymbol{\Theta}}_2X\\
&\qquad -P\mathcal{A}X-P\mathcal{B}_2\bar{\boldsymbol{\Theta}}_2X-P\mathcal{F}_1Y-P\mathcal{B}_1Z-P\mathcal{B}_2\bar{v}_2-P\boldsymbol{b}\Big\}ds\\
&\qquad +\big\{Z-P\big[(\mathcal{C}+\mathcal{D}_2\bar{\boldsymbol{\Theta}}_2)X+\mathcal{B}^\top_1Y+\mathcal{D}_1Z+\mathcal{D}_2\bar{v}_2+\boldsymbol{\sigma}\big]\big\}dW\\
&=-\Big\{\mathcal{A}^\top(\eta+PX)+\mathcal{C}^\top(I-P\mathcal{D}_1)^{-1}\big[P(\mathcal{C}+\mathcal{D}_2\bar{\boldsymbol{\Theta}}_2)X+P\mathcal{B}^\top_1\eta+P\mathcal{B}^\top_1PX\\
&\qquad +P\mathcal{D}_2\bar{v}_2+P\boldsymbol{\sigma}+\zeta\big]+\mathcal{F}^\top_2\bar{v}_2+\boldsymbol{\beta}-\mathcal{A}^\top PX-P\mathcal{F}_1PX\\
&\qquad -(\mathcal{C}^\top+P\mathcal{B}_1)(I-P\mathcal{D}_1)^{-1}P(\mathcal{C}+\mathcal{B}^\top_1P)X-(\mathcal{C}^\top+P\mathcal{B}_1)(I-P\mathcal{D}_1)^{-1}P\mathcal{D}_2\bar{\boldsymbol{\Theta}}_2X\\
&\qquad +P\mathcal{F}_1\eta+P\mathcal{F}_1PX+P\mathcal{B}_2\bar{v}_2+P\boldsymbol{b}+P\mathcal{B}_1(I-P\mathcal{D}_1)^{-1}\big[P(\mathcal{C}+\mathcal{D}_2\bar{\boldsymbol{\Theta}}_2)X\\
&\qquad +P\mathcal{B}^\top_1\eta+P\mathcal{B}^\top_1PX+P\mathcal{D}_2\bar{v}_2+P\boldsymbol{\sigma}+\zeta\big]\Big\}ds+\zeta dW\\
&=-\Big\{\big[\mathcal{A}^\top+(\mathcal{C}^\top+P\mathcal{B}_1)(I-P\mathcal{D}_1)^{-1}P\mathcal{B}^\top_1+P\mathcal{F}_1\big]\eta+(\mathcal{C}^\top+P\mathcal{B}_1)(I-P\mathcal{D}_1)^{-1}\zeta\\
&\qquad +\big[(\mathcal{C}^\top+P\mathcal{B}_1)(I-P\mathcal{D}_1)^{-1}P\mathcal{D}_2+\mathcal{F}^\top_2+P\mathcal{B}_2  \big]\bar{v}_2\\
&\qquad +(\mathcal{C}^\top+P\mathcal{B}_1)(I-P\mathcal{D}_1)^{-1}P\boldsymbol{\sigma}+\boldsymbol{\beta}+P\boldsymbol{b}\Big\}ds+\zeta dW,
\end{aligned}
\end{equation}
and $\eta(T)=\boldsymbol{g}$. According to (\ref{High dimension stationary}) and (\ref{leader Theta condition}), we have
\begin{equation}\nonumber
\begin{aligned}
0&=(\hat{R}_2\bar{\boldsymbol{\Theta}}_2+\mathcal{F}_2)X+\mathcal{D}^\top_2(I-P\mathcal{D}_1)^{-1}
  \big[P(\mathcal{C}+\mathcal{D}_2\bar{\boldsymbol{\Theta}}_2)X+P\mathcal{B}^\top_1\eta+P\mathcal{B}^\top_1PX+\zeta+P\boldsymbol{\sigma}\big]\\
&\qquad +\mathcal{B}^\top_2(\eta+PX)+\big[\hat{R}_2+\mathcal{D}^\top_2(I-P\mathcal{D}_1)^{-1}P\mathcal{D}_2\big]\bar{v}_2+\hat{\rho}\\
&=\big[\mathcal{B}^\top_2+\mathcal{D}^\top_2(I-P\mathcal{D}_1)^{-1}P\mathcal{B}^\top_1\big]\eta+\mathcal{D}^\top_2(I-P\mathcal{D}_1)^{-1}\zeta+\mathcal{D}^\top_2(I-P\mathcal{D}_1)^{-1}P\boldsymbol{\sigma}+\hat{\rho}\\
&\qquad +\big[\hat{R}_2+\mathcal{D}^\top_2(I-P\mathcal{D}_1)^{-1}P\mathcal{D}_2\big]\bar{v}_2,\quad a.e.,\, \mathbb{P}\mbox{-}a.s.
\end{aligned}
\end{equation}
Then
\begin{equation}\label{v2}
\begin{aligned}
&\bar{v}_2=-\big[\hat{R}_2+\mathcal{D}^\top_2(I-P\mathcal{D}_1)^{-1}P\mathcal{D}_2\big]^{-1}\Big\{\big[\mathcal{B}^\top_2+\mathcal{D}^\top_2(I-P\mathcal{D}_1)^{-1}P\mathcal{B}^\top_1\big]\eta\\
&\qquad\qquad +\mathcal{D}^\top_2(I-P\mathcal{D}_1)^{-1}\zeta+\mathcal{D}^\top_2(I-P\mathcal{D}_1)^{-1}P\boldsymbol{\sigma}+\hat{\rho}\Big\},\quad a.e.,\, \mathbb{P}\mbox{-}a.s.
\end{aligned}
\end{equation}
Inserting the above into (\ref{eta}), we achieve
\begin{equation}\nonumber\left\{
\begin{aligned}
d\eta&=-\Big\{\Big[\mathcal{A}^\top+P\mathcal{F}_1+(\mathcal{C}^\top+P\mathcal{B}_1)(I-P\mathcal{D}_1)^{-1}P\mathcal{B}^\top_1\\
     &\qquad\quad -\big[(\mathcal{C}^\top+P\mathcal{B}_1)(I-P\mathcal{D}_1)^{-1}P\mathcal{D}_2+\mathcal{F}^\top_2+P\mathcal{B}_2\big]\\
     &\qquad\quad \times\big[\hat{R}_2+\mathcal{D}^\top_2(I-P\mathcal{D}_1)^{-1}P\mathcal{D}_2\big]^{-1}\big[\mathcal{B}^\top_2+\mathcal{D}^\top_2(I-P\mathcal{D}_1)^{-1}P\mathcal{B}^\top_1\big]\Big]\eta\\
&\qquad +\Big[(\mathcal{C}^\top+P\mathcal{B}_1)(I-P\mathcal{D}_1)^{-1}-\big[(\mathcal{C}^\top+P\mathcal{B}_1)(I-P\mathcal{D}_1)^{-1}P\mathcal{D}_2+\mathcal{F}^\top_2+P\mathcal{B}_2\big]\\
&\qquad\quad \times\big[\hat{R}_2+\mathcal{D}^\top_2(I-P\mathcal{D}_1)^{-1}P\mathcal{D}_2\big]^{-1}\mathcal{D}^\top_2(I-P\mathcal{D}_1)^{-1}\Big]\zeta\\
&\qquad -\big[(\mathcal{C}^\top+P\mathcal{B}_1)(I-P\mathcal{D}_1)^{-1}P\mathcal{D}_2+\mathcal{F}^\top_2+P\mathcal{B}_2\big]\big[\mathcal{D}^\top_2(I-P\mathcal{D}_1)^{-1}P\boldsymbol{\sigma}+\hat{\rho}\big]\\
&\qquad +(\mathcal{C}^\top+P\mathcal{B}_1)(I-P\mathcal{D}_1)^{-1}P\boldsymbol{\sigma}+\boldsymbol{\beta}+P\boldsymbol{b}\Big\}ds+\zeta dW\\
\eta(T)&=\boldsymbol{g}.
\end{aligned}\right.
\end{equation}
In this case,
\begin{equation}\label{u2=Theta2X+v2}
\begin{aligned}
&\bar{u}_2=\bar{\boldsymbol{\Theta}}_2X+\bar{v}_2\\
&=-\big[\hat{R}_2+\mathcal{D}^\top_2(I-P\mathcal{D}_1)^{-1}P\mathcal{D}_2\big]^{-1}\Big\{\big[\mathcal{B}^\top_2P+\mathcal{F}_2+\mathcal{D}^\top_2(I-P\mathcal{D}_1)^{-1}P(\mathcal{C}+\mathcal{B}^\top_1P)\big]X\\
&\qquad +\big[\mathcal{B}^\top_2+\mathcal{D}^\top_2(I-P\mathcal{D}_1)^{-1}P\mathcal{B}^\top_1\big]\eta+\mathcal{D}^\top_2(I-P\mathcal{D}_1)^{-1}\zeta
 +\mathcal{D}^\top_2(I-P\mathcal{D}_1)^{-1}P\boldsymbol{\sigma}+\hat{\rho} \Big\}.
\end{aligned}
\end{equation}
Putting (\ref{Theta2}), (\ref{relationship}) and (\ref{v2}) into the equation of $X$ in (\ref{High dimension optimal system simple}), we obtain (\ref{nonanticipating X}). The proof is complete. $\qquad\Box$

Noting that optimal control $\bar{u}_2(\cdot)$ of the leader has a closed-loop representation (\ref{u2=Theta2X+v2}) with the ``state" $X(\cdot)=\left(\begin{matrix} \bar{x}(\cdot)\\ p^2(\cdot)\end{matrix}\right)$ being the solution to (\ref{nonanticipating X}). Likewise, for the follower, the optimal control $\bar{u}_1(\cdot)$ can also be represented in the following way:
\begin{equation}\label{folower non-anticipating}
\begin{aligned}
\bar{u}_1&=-(\hat{R}^1_{11})^{-1}\bigg\{ \Big[\left(\begin{matrix}
                                              \hat{S}^1_1&0
                                                    \end{matrix}\right)+\left(\begin{matrix}
                                                    0&B^\top_1
                                                    \end{matrix}\right)P
                                                    +\left(\begin{matrix}
                                                    0&D^\top_1
                                                    \end{matrix}\right)(I-P\mathcal{D}_1)^{-1}P(\mathcal{C}+\mathcal{B}^\top_1P)\\
&\qquad -\mathcal{R}\big[\mathcal{D}^\top_2(I-P\mathcal{D}_1)^{-1}P(\mathcal{C}+\mathcal{B}^\top_1P)+\mathcal{B}^\top_2P+\mathcal{F}_2\big] \Big]X \\
&\qquad +\Big[\left(\begin{matrix}
                    0&B^\top_1
                    \end{matrix}\right)
                    +\left(\begin{matrix}
                    0&D^\top_1
                    \end{matrix}\right)(I-P\mathcal{D}_1)^{-1}P\mathcal{B}^\top_1-\mathcal{R}\big[\mathcal{D}^\top_2(I-P\mathcal{D}_1)^{-1}P\mathcal{B}^\top_1+\mathcal{B}^\top_2\big]\eta\\
&\qquad +\Big[\left(\begin{matrix}
                    0&D^\top_1
                    \end{matrix}\right)(I-P\mathcal{D}_1)^{-1}-\mathcal{R}\mathcal{D}^\top_2(I-P\mathcal{D}_1)^{-1}\Big]\zeta+\hat{\rho}^1_1-\mathcal{R}\hat{\rho}\\
&\qquad +\Big[\left(\begin{matrix}
                    0&D^\top_1
                    \end{matrix}\right)(I-P\mathcal{D}_1)^{-1}P-\mathcal{R}\mathcal{D}^\top_2(I-P\mathcal{D}_1)^{-1}P\Big]\boldsymbol{\sigma}\bigg\},
\end{aligned}
\end{equation}
where $\mathcal{R}\triangleq\big[\hat{R}^1_{12}+\left(\begin{matrix} 0&D^\top_1\end{matrix}\right)(I-P\mathcal{D}_1)^{-1}P\mathcal{D}_2\big]\big[\hat{R}_2+\mathcal{D}^\top_2(I-P\mathcal{D}_1)^{-1}P\mathcal{D}_2\big]^{-1}$.

\begin{Remark}
When the coefficients of the inhomogeneous terms and the cross terms are zero, $(\eta(\cdot),\zeta(\cdot)) \equiv (0,0)$ and the Riccati equation \eqref{leader RE} is the same as (3.38) of \cite{Yong2002} for the special case when its coefficients are reduced to deterministic functions (that is $\hat{\Lambda}(\cdot)\equiv0$ in (3.38)). So if the open-loop optimal control can be expressed as a feedback form and the closed-loop optimal strategy exists, the feedback representation of the open-loop optimal control is consistent with the outcome of the closed-loop optimal strategy.
\end{Remark}

\begin{Remark}\label{value function}
In \cite{Yong2002}, the homogeneous state equation and cost functional are considered, and the cost functional does not contain neither the cross terms of state and control nor the cross terms of $u_1$ and $u_2$. In this case, $Y=PX$. By using $It\hat{o}$'s formula to $\langle Y(\cdot),X(\cdot) \rangle $, the value function could be obtained. However, in this paper we consider a general model with the inhomogeneous and the cross terms, this method does not work because of the existence of $g^2$ in \eqref{leader closedloop optimal system simple}. We try to decouple it before using the dimensional expansion technique, but we get a very complicated equation, and it doesn't help us to simplify the cost functional.	

When $b,\sigma,g_i,q^i,\rho^i_1,\rho^i_2,i=1,2$ are all equal to 0, we can use the above method to obtain the value function
\begin{equation}
V_2(x)=\mathbb{E}\langle X(0),Y(0)\rangle=\langle P_1(0)x,x\rangle.
\end{equation}

It is well known that, in stochastic optimal control problems, the value function can be obtained by using $It\hat{o}$'s formula to solutions to the Riccati equation and related BSDE. In the leader-follower game of this paper, since Riccati equation's solution $P(\cdot)\equiv\left(\begin{matrix} P_1(\cdot)&P_2(\cdot)  \\ P_3(\cdot)&P_4(\cdot)  \end{matrix}\right)$ is an $2\times2$ matrix-valued equation which $P_1(\cdot),P_2(\cdot),P_3(\cdot), P_4(\cdot)$ are coupled together and BSDE $(\eta(\cdot),\zeta(\cdot))$ is a vector-valued equation which two components $(\eta_1(\cdot),\zeta_1(\cdot))$ and $(\eta_2(\cdot),\zeta_2(\cdot))$ are also coupled each other, we can not show the expressions for $P_1(\cdot)$ and $(\eta_1(\cdot),\zeta_1(\cdot))$ and then the method of classical optimal control problem fails.
\end{Remark}

%\begin{Remark}\label{non-anticipationg}
%Since the system we consider is inhomogeneous, the backward stochastic differential equations $(\eta(\cdot),\zeta(\cdot))$ has the nontrivial solution. Thus from $(\ref{folower non-anticipating})$ and $(\ref{u2=Theta2X+v2})$, %$u_1(\cdot)$ and $u_2(\cdot)$ are related to $(\eta(\cdot),\zeta(\cdot))$, this is in fact a pseudo-nonanticipating form.
%\end{Remark}

\section{Examples}

In this section, we first give an example to show that for our LQ leader-follower stochastic differential game, the open-loop solvability is weaker than the closed-loop solvability. Then we give a practical example to demonstrate the effectiveness of our theoretic results.

The following example shows that the LQ leader-follower stochastic differential game may have only open-loop Stackelberg equilibria.

\begin{Example}
Consider the following state equation:
\begin{equation}\nonumber\left\{
\begin{aligned}
	dX^{u_1,u_2}(s)&=\big[u_1(s)+u_2(s)\big]ds,\quad s\in [0,1],\\
	X^{u_1,u_2}(0)&=x \in \mathbb{R},
\end{aligned}\right.
\end{equation}
and the cost functionals
\begin{equation}\nonumber
\begin{aligned}
J_1(x;u_1,u_2)&=|X^{u_1,u_2}(1)|^2,\\
J_2(x;u_1,u_2)&=|X^{u_1,u_2}(1)|^2+\int_0^1|u_2(s)|^2ds.
\end{aligned}
\end{equation}
We claim that $\bar{u}_1(\cdot)=-x-u_2(\cdot)$ is an open-loop optimal control of Problem $(SLQ)_f$ for the initial pair $(0,x)$. In this case, it is easy to check that for any $u_2(\cdot) \in \mathcal{U}_2[0,1]$,
\begin{equation}\nonumber
V_1(x)=\mathop{\min}\limits_{u_1(\cdot) \in \mathcal{U}_1[0,T]}J_1(x;u_1,u_2)=0.
\end{equation}
Plugging $\bar{u}_1(\cdot)$ into the state equation, we have
\begin{equation}\nonumber\left\{
\begin{aligned}
dX^{u_2}(s)&=-xds,\quad s\in [0,1],\\
X^{u_2}(0)&=x \in \mathbb{R}.
\end{aligned}\right.
\end{equation}
Then
\begin{equation}\nonumber
X^{u_2}(s)=x-x\int_0^sds=x(1-s),\quad s\in [0,1].
\end{equation}
Thus $X^{u_2}(1)=0$, then
\begin{equation}\nonumber
J_2(x;\bar{u}_1,u_2)=\int_0^1|u_2(s)|^2ds.
\end{equation}
It is easy to know when $\bar{u}_2(\cdot)\equiv0$, this cost functional is minimized. However, such a problem is not closed-loop solvable. In fact, if it is closed-loop solvable, then we may let $(\bar{\theta}_1(\cdot),\bar{v}_1(\cdot))$  be a closed-loop optimal strategy of the follower. Therefore, it is necessary that
\begin{equation}\nonumber
0=X^{\bar{\theta}_1,\bar{v}_1,u_2}(1)=e^{-\int_0^1\bar{\theta}_1(s)ds}x+\int_0^1\big[\bar{v}_1(s)+u_2(s)\big]e^{\int_s^1\bar{\theta}_1(r)dr}ds,\quad\forall x\in \mathbb{R},\,\,u_2(\cdot) \in \mathcal{U}_2[0,1].
\end{equation}
This is impossible.
\end{Example}

We consider a simplified model for a dynamic research and development resource allocation problem in Chen and Cruz \cite{CC1972}, but random noise is considered.

\begin{Example}
It is assumed that two firms are competing with each other for a share of the market for a specific consumption goods. One is a start-up with a relatively small size and its correspondingly research and development spending is lower than the other, and is labeled 1. The other is a more larger company and has a certain market share, which has more experience and spends more on research, and has more access to information, we will label it 2. However, each firm's share of the market does depend on its research and development effort. Let $u_1(\cdot)$ and $u_2(\cdot)$ be the amounts of money invested in research and development by firm 1 and firm 2, respectively. The evolution of the technology gap $X(\cdot)$ between firm 1 and firm 2, is modeled by
\begin{equation}\left\{
\begin{aligned}
dX(t)&=\Big[-\frac{\gamma}{2}X(t)+u_2(t)-\alpha u_1(t)\Big]dt+\beta X(t) dW,\\
 X(0)&=x,
\end{aligned}\right.
\end{equation}
where $\gamma>0$ is the discount rate and the multiplying factor $\alpha \geqslant 1$ accounts for the fact that it is easier for a developing firm, firm 1, to catch up than for firm 2, which is technically advanced, to innovate, and $\beta\in\mathbb{R}$ represents some random environmental effect. Thus, the revenues of these two firms for finite horizon are, respectively,
\begin{equation}
J_1=\mathbb{E}\int_0^T\Big[ -\frac{V}{2x^2_0}X^2(t)-u^2_1(t) \Big] dt,
\end{equation}
\begin{equation}
J_2=\mathbb{E}\int_0^T\Big[ \frac{V}{2x^2_0}X^2(t)-u^2_2(t) \Big] dt,
\end{equation}
where $V\in\mathbb{R}$ is the quasi-rent that is assumed to be constant and $x_0$ is some constant such that, when $X(\cdot)$ reaches $x_0$, the market is completely taken over by firm 2. The goals of both firms are to maximize their own revenues. In order to apply what we have obtained in the previous section, we need to minimize $-J_1$ and $-J_2$. The above problem is an LQ leader-follower stochastic differential game, with the firm 1 as the follower and firm 2 as the leader.

We thus can apply the results obtained in the previous section to seek the closed-loop Stackelberg equilibrium. In the follower's problem, the optimal closed-loop strategy $(\bar{\Theta}_1(\cdot),\bar{v}_1(\cdot))$ of firm 1 exists if and only if the following Riccati equation admits a solution $P^1(\cdot) \in C[0,T]$:
\begin{equation}
\begin{cases}
\dot{P}^1+(\beta^2-\gamma)P^1-\alpha^2(P^1)^2 +\frac{V}{2x_0^2}=0,\\
P^1(T)=0,\\
\end{cases}
\end{equation}	
and the following BSDE admits a solution $(\eta^{1,u_2}(\cdot),\zeta^{1,u_2}(\cdot)) \in L^2_{\mathbb{F}}(0,T;\mathbb{R}) \times L^2_{\mathbb{F}}(0,T;\mathbb{R})$:
\begin{equation}\left\{
\begin{aligned}
d\eta^{1,u_2}&=-\Big[ -\Big(\frac{\gamma}{2}+P^1\alpha^2\Big)\eta^{1,u_2}+\beta\zeta^{1,u_2}+P^1u_2 \Big] ds+\zeta^{1,u_2}dW,\\
\eta^{1,u_2}(T)&=0.
\end{aligned}\right.
\end{equation}	
In this case, the closed-loop optimal strategy $(\bar{\Theta}_1(\cdot),\bar{v}_1(\cdot))\in \mathcal{Q}_1[0,T] \times \mathcal{U}_1[0,T]$ of firm 1 admits the following representation:
\begin{equation}
\bar{\Theta}_1=\alpha P^1,\quad
\bar{v}_1=\alpha\eta^{1,u_2},\qquad a.e.,\, \mathbb{P}\mbox{-}a.s.
\end{equation}
Further, the value function $V_1(\cdot,\cdot)$ is given by
\begin{equation}\label{value-11}
\begin{split}
V_1(x;u_2(\cdot))=\mathbb{E}\biggl\{ \big\langle P^1(0)x,x\big\rangle+2\big\langle \eta^{1,u_2}(0),x\big\rangle+\int_0^T\Big[2\big\langle \eta^{1,u_2},u_2\big\rangle-\big|\alpha\eta^{1,u_2}\big|^2\Big] ds \biggr\}.
\end{split}
\end{equation}
Then we will concentrate on the leader's problem. If the leader's problem of firm 2 is closed-loop solvable, then the closed-loop optimal strategy $(\boldsymbol\Theta_2(\cdot),\bar{v}_2(\cdot))\equiv(\bar{\Theta}_2(\cdot),\bar{\tilde{\Theta}}_2(\cdot),\bar{v}_2(\cdot))\in \mathcal{Q}_2[0,T] \times \mathcal{Q}_2[0,T] \times \mathcal{U}_2[0,T]$ admits the following representation:
\begin{equation}\label{close-loop leader--}
\left\{
\begin{aligned}
\bar{\Theta}_2&=-P_1,\\
\bar{\tilde{\Theta}}_2&=-(P_2+P^1),\\
\bar{v}_2&=0,\qquad\qquad a.e.,\, \mathbb{P}\mbox{-}a.s.
\end{aligned}
\right.
\end{equation}
where $P(\cdot)\equiv\left(\begin{matrix} P_1(\cdot)&P_2(\cdot) \\ P_2(\cdot)^\top&P_4(\cdot) \end{matrix}\right) \in C([0,T];\mathbb{S}^{2n \times 2n})$ is the solution to the Riccati equation:
\begin{equation}\label{leader RE--}
\begin{cases}
&\dot{P}+\mathcal{A}^\top P+P\mathcal{A}+\mathcal{C}^\top P\mathcal{C}+P\mathcal{F}_1P+\mathcal{Q}_2-(P\mathcal{B}_2+\mathcal{F}_2)(\mathcal{B}^\top_2P+\mathcal{F}^\top_2)=0,\\
&P(T)=0.
\end{cases}
\end{equation}
where
\begin{equation*}
\begin{cases}
\mathcal{A} \triangleq\left(\begin{matrix} -\big(\frac{\gamma}{2}+P^1\alpha^2\big) & 0 \\ 0 & -\big(\frac{\gamma}{2}+P^1\alpha^2\big) \end{matrix}\right),\quad
\mathcal{F}_1 \triangleq\left(\begin{matrix} 0 & -\alpha^2 \\ -\alpha^2 & 0 \end{matrix}\right),\quad\\
\mathcal{B}_2 \triangleq\left(\begin{matrix} 1 \\ 0 \end{matrix}\right),\quad
\mathcal{C} \triangleq\left(\begin{matrix} \beta & 0  \\ 0 & \beta \end{matrix}\right),\quad
\mathcal{Q}_2 \triangleq\left(\begin{matrix} -\frac{V}{2x_0^2} & 0 \\ 0 & 0 \end{matrix}\right),\quad
\mathcal{F}_2 \triangleq\left(\begin{matrix} 0 \\ P^1 \end{matrix}\right).
\end{cases}
\end{equation*}
Further, as Remark \ref{value function}, the value function of the leader admits the following representation:
\begin{equation}\label{value function-leader}
V_2(x)=\big\langle P_1(0)x,x \big\rangle.
\end{equation}

Consider this problem with the parameters' values for all plots are $\gamma=0.05$, $\alpha^2=2.5$, $\beta^2=0.000065$, $\frac{V}{2x_0^2}=0.8$, the solutions $P^1(\cdot)$ and $P_1(\cdot),\,P_2(\cdot),\,P_4(\cdot)$ to the above Riccati equations are shown in the following figure.
\begin{figure}[htbp]
	\centering
	\includegraphics[scale=0.3]{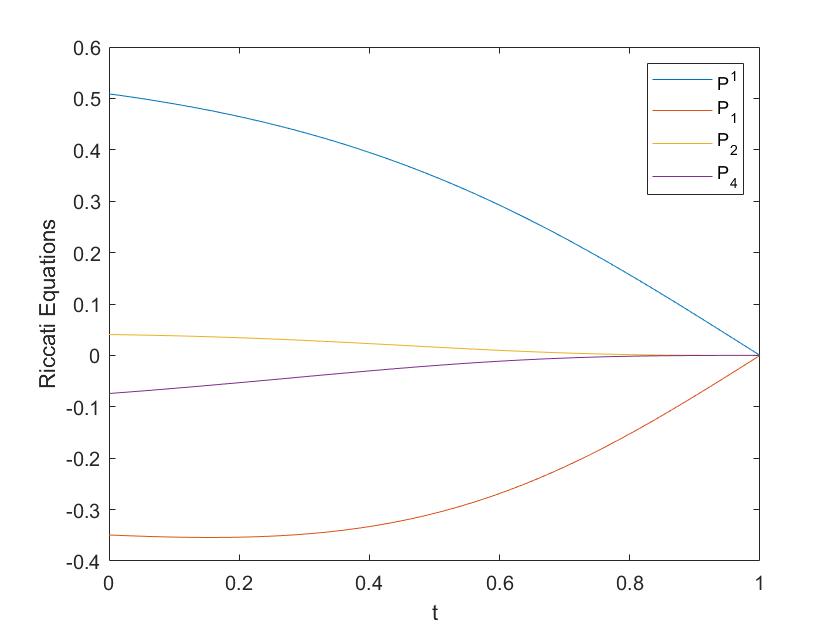}
    \caption{The trajectories of $P^1(\cdot)$ and $P_1(\cdot),\,P_2(\cdot),\,P_4(\cdot)$}
\end{figure}
\end{Example}

\section{Concluding remarks}

In this paper, we have investigated the closed-loop representation of open-loop Stackelberg equilibria and the closed-loop Stackelberg equilibria, for an LQ leader-follower stochastic differential game. We have come to the conclusion that the existence of open-loop Stackelberg equilibria is equivalent the solvability of some systems of FBSDEs together with some convexity conditions (Proposition \ref{fsc} and Theorem \ref{lsc}). For the follower's problem, the existence of closed-loop optimal strategies is equivalent to the solvabilities of some Riccati equation and some BSDE (Theorem \ref{Th-cl-f}); for the leader's problem, necessary conditions for the existence of  closed-loop optimal strategies are given (Theorem \ref{Th-cl-l}). Up to now, we could not obtain sufficient conditions for the existence of closed-loop optimal strategies of the leader, since the state equation of the leader is an fully-coupled FBSDE and the completion-of-square technique is invalid. This is a gap for the closed-loop Stackelberg equilibria of our LQ leader-follower stochastic differential game. We wish to fill in this gap in the future.
\par Problems with random coefficients, and mean-field'type state equations and cost functionals, are interesting and challenging research topics. We will consider them in the near future.


\begin{thebibliography}{0}
	
	\bibitem{BB1981} A. Bagchi, T. Ba\c{s}ar, Stackelberg strategies in linear-quadratic stochastic differential games, \emph{J. Optim. Theory Appl.}, {\bf 35}(3), 443-464, 1981.	
	
	\bibitem{BCS2015} A. Bensoussan, S. K. Chen, and S. P. Sethi,  The maximum principle for global solutions of stochastic Stackelberg differential games, \emph{SIAM J. Control Optim.}, {\bf 53}(4), 1956-1981, 2015.
	
	\bibitem{BO1998} T. Ba\c{s}ar, G. J. Olsder, \emph{Dynamic Noncooperative Game Theory, 2nd Edition}, SIAM, Philadelphia, 1998.
	
	\bibitem{CA1976} D. Castanon, M. Athans, On stochastic dynamic Stackelberg strategies, \emph{Automatica}, {\bf 12}(2), 177-183, 1976.
	
	\bibitem{CC1972} C. I. Chen, J. B. Cruz, Stackelberg solution for two-person games with biased information patterns,  \emph{IEEE Trans. Autom. Control}, {\bf 17}(6), 791-798, 1972.
	
	\bibitem{DW2019} K. Du, Z. Wu, Linear-quadratic Stackelberg game for mean-field backward stochastic differential system and applicationm, \emph{Math. Prob. Eng.}, {\bf 2019}, Article ID 1798585, 17 pages, 2019.
	
    \bibitem{HSW2020}J. H. Huang, K. H. Si, and Z. Wu, Linear-quadratic mixed Stackelberg-Nash stochastic differential game with major-minor agents. \emph{Appl. Math. Optim.}, https://doi.org/10.1007/s00245-020-09713-z.

	\bibitem{LJZ2019} Y. N. Lin, X. S. Jiang, and W. H. Zhang, Open-loop Stackelberg strategy for the linear quadratic mean-field stochastic differential game, \emph{IEEE Tran. Autom. Control}, {\bf 64}(1), 97-110, 2019.
	
	\bibitem{LSY2016} X. Li, J. R. Sun, and J. M. Yong, Mean-field stochastic linear quadratic optimal control problems: closed-loop solvability, \emph{Proba. Uncer. Quan. Risk}, {\bf 1}(1), 24 pages, 2016.
	
	\bibitem{LSY2020} X. Li, J. T. Shi, and J. M. Yong, Mean-field linear-quadratic stochastic differential games in an infinite horizon, 2020. https://arxiv.org/abs/2007.06130
	
	\bibitem{LiYu2018} N. Li, Z. Y. Yu, Forward-backward stochastic differential equations and linear-quadratic generalized Stackelberg games, \emph{SIAM J. Control Optim.}, {\bf 56}(6), 4148-4180, 2018.
	
    \bibitem{LiShi2021} Z. X. Li, J. T. Shi, Linear quadratic Stackelberg stochastic differential games: Closed-loop solvability, accepted by \emph{The 40th Chinese Control Conference}, Shanghai, July 26-28, 2021.

	\bibitem{MB2018} J. Moon, T. Ba\c{s}ar, Linear quadratic mean field Stackelberg differential games, \emph{Automatica}, {\bf 97}, 200-213, 2018.
	
	\bibitem{MX2015} H. Mukaidani, H. Xu, Stackelberg strategies for stochastic systems with multiple followers, \emph{Automatica}, {\bf 53}, 53-59, 2015.
	
	\bibitem{OSU2013} B. \O ksendal, L. Sandal, and  J. Ub\o e,  Stochastic Stackelberg equilibria with applications to time dependent newsvendor models, \emph{J. Econ. Dyna. $\&$ Control.}, {\bf 37}(7), 1284-1299, 2013.
	
	\bibitem{Stackelberg1934} H. von Stackelberg, \emph{Marktform und Gleichgewicht}, Springer, Vienna, 1934. (An English translation appeared in \emph{The Theory of the Market Economy}, Oxford University Press, 1952.)
	
	\bibitem{SC1973-1} M. Simaan, J. B. Cruz Jr., On the Stackelberg game strategy in non-zero games, \emph{J. Optim. Theory Appl.}, {\bf 11}(5), 533-555, 1973.
	
	\bibitem{SC1973-2} M. Simaan, J. B. Cruz Jr., Additional aspects of the Stackelberg strategy in nonzero-sum games, \emph{J. Optim. Theory Appl.}, {\bf 11}(6), 613-626, 1973.
	
	\bibitem{SWX2016} J. T. Shi, G. C. Wang, and J. Xiong, Leader-follower stochastic differential game with asymmetric information and applications, \emph{Automatica}, {\bf 63}, 60-73, 2016.
	
	\bibitem{SWX2017} J. T. Shi, G. C. Wang and J. Xiong,  Linear-quadratic stochastic Stackelberg differential game with asymmetric information, \emph{Sci. China Infor. Sci.}, {\bf 60}, 1-15, 2017.
	
	\bibitem{SWX2020} J. T. Shi, G. C. Wang and J. Xiong,  Stochastic linear-quadratic Stackelberg differential game with overlapping information, \emph{ESAIM: COCV.}, {\bf 26}, Article Number 83, 2020.
	
	\bibitem{SLY2016} J. R. Sun, X. Li, and J. M. Yong, Open-loop and closed-loop solvabilities for stochastic linear quadratic optimal control problems, \emph{SIAM J.Control Optim.}, {\bf 54}(5), 2274-2308, 2016.

	\bibitem{SY2014} J. R. Sun, J. M. Yong, Linear quadratic stocahastic differential games: open-loop and closed-loop saddle points, \emph{SIAM J. Control Optim.}, {\bf 52}(6), 4082-4121, 2014.
	
	\bibitem{SY2018} J. R. Sun, J. M. Yong, Stochastic linear quadratic optimal control problems in infinite horizon, \emph{Appl. Math. Optim.}, {\bf 78}, 145-183, 2018.
	
	\bibitem{SY2019} J. R. Sun,  J. M. Yong, Linear quadratic stocahastic two-person nonzero-sum differential games: open-loop and closed-loop Nash equilibria, \emph{Stoc. Proc. Appl.}, {\bf 129}(2), 381-418, 2019.
	
	\bibitem{SunYong2020} J. R. Sun, J. M. Yong, Stochastic Linear-Quadratic Optimal Control Theory: Open-Loop and Closed-Loop Solutions, \emph{Springer Briefs in Mathematics}, 2020.
	
	\bibitem{WZ2020} G. C. Wang, S. S. Zhang, A Mean-field linear-quadratic stochastic Stackelberg differential game with one leader and two followers, \emph{J.  Syst. Sci. Complex.}, {\bf 33}, 1383-1401, 2020.
	
	\bibitem{XSZ2018} J. J. Xu, J. T. Shi, and H. S. Zhang,  A leader-follower stochastic linear quadratic differential game with time delay, \emph{ Sci. China Infor. Sci.}, {\bf 61}, 112202:1-112202:13, 2018.
	
	\bibitem{XZ2016} J. J. Xu, H. S. Zhang,  Sufficient and necessary open-loop Stackelberg strategy for two-player game with time delay, \emph{IEEE Trans. Cyber.}, {\bf 46}(2), 438-449, 2016.
	
	\bibitem{Yong2002} J. M. Yong, A leader-follower stochastic linear quadratic differential games, \emph{SIAM J. Control Optim.}, {\bf 41}(4), 1015-1041, 2002.
	
	\bibitem{Yong2006} J. M. Yong, Linear forward-backward stochastic differential equations with random coefficients, \emph{ Probab. Theory Relat. Fields}, {\bf 135}(1), 53-83, 2006.
	
	\bibitem{ZS2020}Y. Y. Zheng, J. T. Shi, A Stackelberg game of backward stochastic differential equations with applications, \emph{Dyna. Games Appl.}, {\bf 10}(4), 968-992, 2020.
	
\end{thebibliography}
\end{document}